\documentclass[11pt,letterpaper]{article}
\usepackage[T1]{fontenc}
\usepackage{lmodern}
\usepackage[margin=1in]{geometry}
\usepackage{amsmath,amssymb,amsthm,mathtools,bm}
\usepackage{thmtools}
\usepackage{graphicx,adjustbox,microtype,enumitem,booktabs,placeins,needspace}
\usepackage{xcolor}
\usepackage{tikz}
\usetikzlibrary{calc,arrows.meta,angles,quotes,positioning,patterns}
\usepackage{flafter}
\usepackage[font=small,labelfont=bf,labelsep=period]{caption}
\usepackage[authoryear,round]{natbib}
\bibpunct[, ]{(}{)}{,}{a}{}{,}

\usepackage{xurl}
\usepackage[hidelinks,hypertexnames=false]{hyperref}
\hypersetup{pdftitle={Asymptotic Optimality of Hexagonal Layouts in Euclidean Location-Routing},pdfauthor={Xin Chen, Weijun Xie}}

\allowdisplaybreaks[1]
\newcommand{\doi}[1]{\href{https://doi.org/#1}{\nolinkurl{#1}}}
\newcommand{\R}{\mathbb{R}}
\newcommand{\dd}{\,\mathrm{d}}
\newcommand{\dist}{\operatorname{dist}}
\newcommand{\arcosh}{\operatorname{arcosh}}
\newcommand{\cT}{\mathcal{T}}
\newcommand{\cC}{\mathcal{C}}
\newcommand{\rs}{r_{*}}
\newcommand{\norm}[1]{\lVert #1\rVert}
\newcommand{\abs}[1]{\lvert #1\rvert}

\colorlet{acc}{teal!70!black}
\colorlet{capc}{orange!60}
\colorlet{soft}{teal!13}
\tikzset{dot/.style={circle,fill,inner sep=0pt,minimum size=2.6pt},
odot/.style={circle,draw,fill=white,inner sep=0pt,minimum size=3.2pt},
lab/.style={font=\scriptsize},
>={Stealth[length=1.6mm,width=1.3mm]}}

\newcommand{\eqCsix}{\eqref{eq:Csix}}
\newcommand{\eqabR}{\eqref{eq:abR}}
\newcommand{\lemHex}{Lemma~\ref{lem:hexagon}}
\newcommand{\lemSpacing}{Lemma~\ref{lem:existence}}

\newcommand{\appSide}{Appendix~\ref{app:side}}
\newcommand{\appOverlap}{Appendix~\ref{app:overlap}}
\newcommand{\appV}{Appendix~\ref{app:V}}
\newcommand{\appEta}{Appendix~\ref{app:eta}}

\newcommand{\appBdry}{Appendix~\ref{app:boundary}}

\newcommand{\citeGruber}{\citealp{Gruber1999}}

\newcommand{\citeXO}{\citealp{XieOuyang2015}}
\DeclareMathOperator{\artanh}{artanh}

\title{Asymptotic Optimality of Hexagonal Layouts in Euclidean Location--Routing}
\author{Xin Chen\\
{\small The University of Texas at Dallas}\\
{\small Richardson, Texas, USA}
\and
Weijun Xie\\
{\small City University of Hong Kong}\\
{\small Hong Kong}}
\date{}
\begin{document}
\raggedbottom
\maketitle
\begin{abstract}
We study the joint location of facilities and the design of a closed replenishment tour through them. Customers are uniformly distributed over a square and travel to their nearest facility. The objective is to minimize the sum of total customer travel distance and the weighted length of the replenishment tour, with all distances measured in the Euclidean metric. We consider an asymptotic regime in which the area of the square and the number of facilities increase together at unit facility density. For every fixed routing weight, we prove that the optimal cost per facility converges to that of an optimized staggered-row layout with hexagonal service regions. We obtain this result by constructing a feasible layout and tour in the finite square and deriving a matching lower bound valid for arbitrary facility locations and visiting orders. The upper and lower bounds differ from their common leading term by at most a constant times the square root of the number of facilities and a constant, respectively.
\end{abstract}

\section{Introduction}\label{sec:introduction}

Facility locations jointly determine the cost of serving customers 
and the cost of replenishing facilities. Reducing one cost may increase the other. Dispersing facilities over a service region improves 
customer access but may require a longer replenishment route. 
Conversely, locating facilities close to one another along a route 
reduces replenishment travel but may increase the distance customers 
must travel to reach a facility. This trade-off arises in the design 
of retail networks, customer collection points, and parcel-locker 
systems \citep{Cachon2014,Boysen2021,Janinhoff2024}. A fundamental 
question is how facilities should be arranged in space to minimize 
the combined cost of customer and replenishment travel.

\citet{XieOuyang2015} study this question on an infinite homogeneous 
plane with facility setup, outbound transportation, and inbound 
replenishment costs. Under the Euclidean metric, they propose a 
staggered-row layout in which each facility serves an elongated 
cyclic hexagon. By comparing the cost of this layout with a relaxed 
lower bound, they show that the layout is nearly optimal. Whether 
the hexagonal layout is asymptotically optimal, however, remains 
unresolved. The difficulty lies in deriving a matching lower bound 
that applies to arbitrary facility locations and replenishment 
routes, without imposing restrictions on the shapes or areas of 
individual service regions.

In this paper, we resolve this question in a fixed-density 
setting. We consider the joint placement of a given number of distinct facilities 
in a square whose area equals the facility count and the design of a closed replenishment 
tour visiting all facilities. The objective is to minimize the 
sum of customer travel distances to their nearest facilities and 
the weighted length of the replenishment tour. We derive an upper 
bound using the hexagonal layout of \citet{XieOuyang2015} and a 
lower bound that holds for every feasible layout and tour. The 
two bounds share the same leading term and differ by at most 
a constant times the square root of the number of facilities for any fixed routing weight. Thus, as the service 
region and the number of facilities grow at the same rate, the 
hexagonal layout is asymptotically optimal. The results also 
characterize the optimal spacing within the staggered-row family as a function of 
the relative cost of replenishment travel.

\subsection{Literature review}
Joint location and routing models determine facility locations, 
customer assignments, and vehicle routes. The formulations 
and solution methods are reviewed by 
\citet{NagySalhi2007,ProdhonPrins2014,DrexlSchneider2015,SchneiderDrexl2017}. 
Subsequent reviews classify a broader range of models 
\citep{Mara2021,HosodaIrohara2022,MarianovEiselt2024}, including 
single- and multi-echelon systems, multiple periods, and multiple 
objectives \citep{Cavagnini2026a,Cavagnini2026b}. Other extensions 
consider intermediate stops, hub routing, sustainability, resilience, 
and last-mile distribution 
\citep{Schiffer2019,Wandelt2025,Figueiredo2025,Kahalimoghadam2026}. 
For discrete location--routing problems, 
\citet{Harks2013,Zhao2026} develop approximation algorithms for 
capacitated systems, \citet{Liguori2023} strengthen exact solution 
methods using knapsack cuts, and \citet{Sobhanan2025} combine genetic 
algorithms with neural cost prediction for hierarchical vehicle 
routing. Our focus is on the spatial structure of optimal facility 
layouts and replenishment routes. We derive an additive bound on 
the optimality gap per facility, which is proven to vanish as the service region grows.

An application of our model is out-of-home delivery, where 
customers travel to service points supplied by a carrier. 
\citet{Boysen2021,Janinhoff2024,Mohri2024,ZhangDemir2025} review 
last-mile delivery concepts, alternative delivery points, and 
parcel-locker systems. \citet{Janjevic2021} study the design of 
distribution networks with multiple tiers, service levels, and 
transportation modes. \citet{Arslan2021} formulate a 
location-or-routing problem in which service can be provided through 
either location or routing decisions. More recently, 
\citet{CachonCetinYu2026} study a dynamic vehicle routing problem 
with service location flexibility, allowing each customer to be 
served within a given radius of the customer's address. Under their 
fixed-route policy, the vehicle selects its direction by an 
independent fair coin toss at each departure. Our model considers 
a static setting in which customers use their nearest facility and 
a vehicle replenishes all facilities along a closed tour. In the 
absence of capacity constraints and service times, the facility 
locations and their visiting order determine the two travel costs.

Our analysis builds on the literature on continuous approximation 
and Euclidean location theory. Continuous approximation represents 
discrete customer demand by a spatial demand field 
\citep{Franceschetti2017,Ansari2018}, with applications to 
route-length estimation \citep{Daganzo1984}, terminal system design 
\citep{OuyangDaganzo2006}, and planar hub location--routing 
\citep{Ghaffarinasab2018}. Recent applications include integrated 
truck--drone delivery \citep{Zhang2024}, strategic last-mile 
distribution planning \citep{Kirubakaran2026}, and large-scale 
delivery network design evaluated through agent-based simulation 
\citep{Osorio2026}. For the case without routing costs, 
\citet{HaimovichMagnanti1988} derive extremum and asymptotic 
properties of hexagonal partitions under uniform demand, and 
\citet{Gruber1999} provides a proof of the planar moment theorem 
based on the number of polygon sides. 
\citet{BourneCristoferi2021} derive asymptotic triangular-lattice 
optimality for a different location objective involving quadratic 
transport costs and mass penalization. The Euclidean traveling 
salesman asymptotics of \citet{Beardwood1959} concern randomly 
sampled points, whereas the facility locations in our model are 
decision variables. We extend the geometric side-count argument 
to account for replenishment travel, using Abel integral equations 
\citep{GorenfloVessella1991} to construct a radius-dependent 
Lagrange multiplier and derive the required inequality for 
individual triangles.

The two studies most closely related to ours are 
\citet{CarlssonJia2015} and \citet{XieOuyang2015}. 
\citet{CarlssonJia2015} consider the length of a tour through 
facilities, the integral of customer distance to the nearest 
facility, and facility installation costs. They obtain honeycomb 
and Archimedean-spiral configurations in two limiting regimes 
and develop a constant-factor approximation algorithm. Fixing 
the number of facilities and normalizing the customer-travel 
coefficient in their model recovers our objective. 
\citet{XieOuyang2015} construct layouts of elongated cyclic 
hexagons on an infinite homogeneous plane and derive a lower 
bound by approximating a service region with a disk cut by two 
parallel lines, with equally spaced facilities along a straight 
replenishment route. These idealized regions need not tile the 
plane. Their reported optimality gap within 0.3\% is for an 
objective that includes facility setup costs and optimizes 
facility density. (Section and result numbers for \citet{XieOuyang2015} refer to the author version, arXiv:1401.2415v2; the published article appeared in 2015.)
We adopt their hexagonal construction and extend it to a finite 
square with an exact facility count and a closed replenishment 
tour. Unlike their geometric relaxation, our lower bound uses 
the disk overlaps along the actual tour and the number of sides 
of the Voronoi cells. The resulting upper and lower bounds have 
the same leading coefficient, proving the asymptotic optimality 
of the elongated hexagonal construction for every fixed routing 
weight at unit facility density. Section~\ref{sec:numerics} 
compares our limiting cost with their relaxed lower bound under 
the same normalization.

\subsection{Summary of contributions}

We consider a joint facility location and replenishment routing 
problem in a square whose area equals the number of facilities. 
For every fixed nonnegative routing weight, we derive finite-region 
bounds with a common leading term given by the minimum cost per 
facility of the staggered-row hexagonal layout multiplied by the 
number of facilities. Our main contributions are as follows.
\begin{itemize}
\item  Theorem~\ref{thm:upper} extends the 
staggered-row construction of \citet{XieOuyang2015} to a finite 
square with exactly the required number of facilities and a closed 
replenishment tour. The additional cost of serving boundary regions, 
connecting rows, and completing the facility count is at most 
proportional to the square root of the number of facilities.

\item Theorem~\ref{thm:lower} derives a 
matching lower bound, up to an additive constant independent of 
the number of facilities. We formulate a geometric relaxation 
that links the Voronoi triangles of service regions to disk overlaps 
along the replenishment tour. Using Abel integral equations, we 
construct Lagrange multipliers under which all six triangles of 
the reference hexagon attain the lower bound on individual 
triangles. The resulting bound holds for arbitrary facility 
locations and visiting orders.

\item Corollary~\ref{cor:gap} 
proves the asymptotic optimality of the staggered-row hexagonal 
construction and quantifies its finite-size optimality gap. 
Lemma~\ref{lem:spacing} further characterizes the unique optimal 
spacing within the staggered-row family as a function of the 
routing weight.
\end{itemize}

The remainder of the paper is organized as follows. Section~\ref{sec:model} defines the model and reference construction and gives the main results. Section~\ref{sec:upper} proves the upper bound. Sections~\ref{sec:lb1} and~\ref{sec:lb2} prove the lower bound: Section~\ref{sec:lb1} derives the relaxation, and Section~\ref{sec:lb2} constructs its multipliers. Section~\ref{sec:numerics} reports numerical results. Detailed proofs are in Appendices~\ref{app:upper}--\ref{app:global}.

\paragraph{Notation.}
The Euclidean norm is $\|\cdot\|$, and $\dist(x,X)=\min_{y\in X}\|x-y\|$. We write $B_t(x)=\{y:\|y-x\|\le t\}$, use $|A|$ for the area of a planar set and $\#X$ for the cardinality of a finite set, and put $(z)_+=\max\{z,0\}$. The indicator of $E$ is $\mathbf 1_E$. The relation $f\sim g$ means $f/g\to1$ in the given limit. Throughout the large-$N$ analysis, $r$ is fixed, and constants indexed by $r$ do not depend on $N$.

\section{Model, reference construction, and main results}\label{sec:model}

Let $N\ge3$ be an integer. Place $N$ distinct facilities in a square with unit customer density:
\begin{equation}\label{eq:domain}
Q_N=[0,\sqrt N]^2,\qquad X=\{x_1,\ldots,x_N\}\subset Q_N.
\end{equation}
Each customer uses the nearest facility. The total customer-distance cost is
\begin{equation}\label{eq:service}
M_N(X)=\int_{Q_N}\min_i\|x-x_i\|\dd x.
\end{equation}
A replenishment vehicle visits all facilities and returns to its starting point. For a cyclic order $\sigma$, with $\sigma(N+1)=\sigma(1)$, define
\begin{equation}\label{eq:route}
L(X,\sigma)=\sum_{j=1}^N\|x_{\sigma(j+1)}-x_{\sigma(j)}\|,
\qquad L_N(X)=\min_\sigma L(X,\sigma).
\end{equation}
For a given relative routing cost $r\ge0$, the optimal combined cost is
\begin{equation}\label{eq:opt}
A_N(r)=\inf_{\substack{X\subset Q_N\text{ distinct}\\\#X=N}}\{M_N(X)+rL_N(X)\}.
\end{equation}
The area and facility count both equal $N$, so the average facility density is one. The service regions are clipped Voronoi cells and may have unequal areas. Apart from the conditions above, facility locations and the visiting order are unrestricted. We omit facility setup cost since the number of facilities is fixed for each problem.

For $N$ facilities in a square of area $A>0$, put $\delta=N/A$ and rescale coordinates by $\sqrt\delta$. The rescaled square has area $N$ and facility density one. With customer density one in the original square, write $\widetilde X=\sqrt\delta X$. A change of variables gives
\[
\int_{[0,\sqrt A]^2}\dist(x,X)\dd x+rL(X,\sigma)
=\delta^{-3/2}\{M_N(\widetilde X)+r\delta L(\widetilde X,\sigma)\}.
\]
Hence the routing weight in the normalized problem is $r\delta=rN/A$.

\subsection{Staggered rows}
In the construction of \citet{XieOuyang2015}, facilities in each row are separated by $d$, adjacent rows are separated by $1/d$, and successive rows are offset by $d/2$. Let
\begin{equation}\label{eq:lattice}
d_0=\sqrt{2/\sqrt3},\qquad
\Lambda_d=\left\{\left(d\left(i+\frac j2\right),\frac jd\right):i,j\in\mathbb Z\right\},\qquad 0<d\le d_0.
\end{equation}

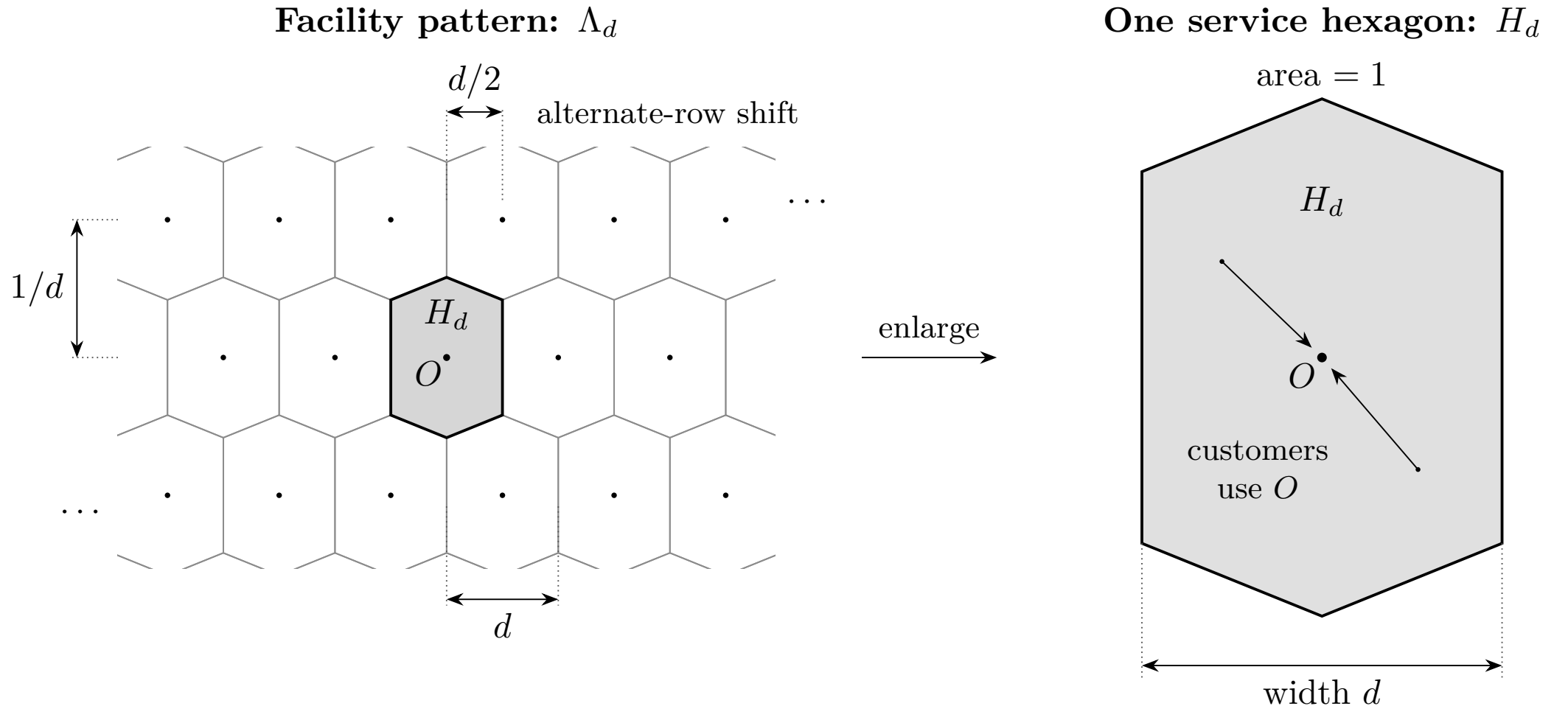
\begin{figure}[!htbp]\centering
\begingroup
\resizebox{0.86\textwidth}{!}{
\begin{tikzpicture}[color=black,>=Stealth,font=\small,
every node/.style={text=black},dim/.style={<->,thin},guide/.style={densely dotted,thin,black!65}]
\pgfmathsetmacro{\fd}{0.9}
\pgfmathsetmacro{\fa}{\fd/2}
\pgfmathsetmacro{\fr}{1/(2*\fd)+\fd^3/8}
\pgfmathsetmacro{\fv}{1/\fd-\fr}
\node[font=\small\bfseries] at (3.75,7.0) {Facility pattern: $\Lambda_d$};
\node[font=\small\bfseries] at (12.5,7.0) {One service hexagon: $H_d$};
\begin{scope}[shift={(3.75,3.65)},scale=1.24]
\begin{scope}
\clip (-2.65,-1.7) rectangle (2.65,1.7);
\foreach \k in {-2,...,2}{
\foreach \j in {-4,...,4}{
\pgfmathsetmacro{\px}{\fd*(\j+\k/2)}
\pgfmathsetmacro{\py}{\k/\fd}
\begin{scope}[shift={(\px,\py)}]
\draw[black!45,thin] (0,\fr)--(\fa,\fv)--(\fa,-\fv)--(0,-\fr)--(-\fa,-\fv)--(-\fa,\fv)--cycle;
\fill[black] (0,0) circle (0.022);
\end{scope}
}
}
\filldraw[fill=black!16,draw=black,thick]
(0,\fr)--(\fa,\fv)--(\fa,-\fv)--(0,-\fr)--(-\fa,-\fv)--(-\fa,\fv)--cycle;
\fill[black] (0,0) circle (0.028);
\node[above,font=\small] at (0,0.12) {$H_d$};
\node[below left,inner sep=1pt] at (0,0) {$O$};
\end{scope}
\draw[guide] (0,1.27)--(0,2.02);
\draw[guide] (\fa,1.27)--(\fa,2.02);
\draw[dim] (0,1.98)--(\fa,1.98) node[midway,above] {$d/2$};
\node[anchor=west,font=\footnotesize] at (0.62,1.98) {alternate-row shift};
\draw[guide] (-2.65,0)--(-3.02,0);
\draw[guide] (-2.65,{1/\fd})--(-3.02,{1/\fd});
\draw[dim] (-2.98,0)--(-2.98,{1/\fd}) node[midway,left] {$1/d$};
\draw[guide] (0,-1.2)--(0,-2.0);
\draw[guide] (\fd,-1.2)--(\fd,-2.0);
\draw[dim] (0,-1.95)--(\fd,-1.95) node[midway,below] {$d$};
\node at (-2.93,-1.25) {$\cdots$};
\node at (2.93,1.25) {$\cdots$};
\end{scope}
\draw[->,thin] (7.9,3.65)--(9.25,3.65) node[midway,above,font=\footnotesize] {enlarge};
\begin{scope}[shift={(12.5,3.65)},scale=4]
\filldraw[fill=black!12,draw=black,thick]
(0,\fr)--(\fa,\fv)--(\fa,-\fv)--(0,-\fr)--(-\fa,-\fv)--(-\fa,\fv)--cycle;
\fill[black] (0,0) circle (0.012);
\node[above,font=\small] at (0,0.32) {$H_d$};
\node[below left,inner sep=1.5pt] at (0,0) {$O$};
\foreach \x/\y in {-0.25/0.24,0.24/-0.28}{
\fill[black] (\x,\y) circle (0.006);
\draw[->,thin,shorten >=4pt] (\x,\y)--(0,0);
}
\node[align=center,font=\footnotesize] at (-0.16,-0.28) {customers\\use $O$};
\node[above,font=\small] at (0,\fr) {area $=1$};
\draw[guide] (-\fa,-\fv)--(-\fa,-0.80);
\draw[guide] (\fa,-\fv)--(\fa,-0.80);
\draw[dim] (-\fa,-0.77)--(\fa,-0.77) node[midway,below,font=\small] {width $d$};
\end{scope}
\end{tikzpicture}
}
\endgroup
\caption{Staggered-row facility pattern $\Lambda_d$ and the Voronoi cell $H_d$ of one facility. Row spacing is $1/d$, within-row spacing is $d$, and each full cell has area one. The drawing uses $d=0.9$.}\label{fig:pattern}
\end{figure}

A fundamental parallelogram of $\Lambda_d$ has base $d$ and height $1/d$, hence area one and one facility per unit area. Within a row, adjacent facilities are distance $d$ apart. The nearest facilities in an adjacent row are displaced by $d/2$ horizontally and $1/d$ vertically, giving distance $\sqrt{d^{-2}+d^2/4}$. At $d=d_0$ these two distances agree, and the cell is regular. For $0<d\le d_0$, the within-row neighbors are among the closest facilities. Decreasing $d$ makes the cells narrower and taller and shortens the within-row tour edges (Figure~\ref{fig:pattern}). Denote the Voronoi cell of the origin by $H_d$ and define

\begin{equation}\label{eq:Csix}
\phi(d)=\int_{H_d}\|x\|\dd x,\qquad C_6(r)=\min_{0<d\le d_0}\{\phi(d)+rd\}.
\end{equation}
Here $\phi(d)$ is the service cost of one cell, and $rd$ is the weighted cost of one within-row tour edge. The finite construction also includes connections between rows.

\subsection{Main results}\label{sec:mainresults}
The following finite-square bounds have the common leading coefficient $C_6(r)$. Lemma~\ref{lem:existence} shows that the minimum defining this coefficient is attained at a positive spacing.

\Needspace{4\baselineskip}
\begin{restatable}[Finite-square upper bound]{theorem}{finiteSquareUpperBound}\label{thm:upper}
Fix $r\ge0$. For every $0<d\le d_0$ and every integer $N\ge\max\{3,\lceil d^{-2}\rceil\}$, there are a layout $X$ of $N$ distinct facilities in $Q_N$ and a closed visiting order $\sigma$ such that
\begin{equation}\label{eq:upperparts}
M_N(X)\le N\phi(d)+4(d^2+d^{-2})\sqrt N,\qquad L(X,\sigma)\le Nd+4\sqrt N.
\end{equation}
In particular, choose a minimizing spacing $d$ in \eqref{eq:Csix}, fixed independently of $N$, and set
\begin{equation}\label{eq:Dr}
D_r=4(d^2+d^{-2}+r).
\end{equation}
For the same size threshold,
\begin{equation}\label{eq:upper}
A_N(r)\le NC_6(r)+D_r\sqrt N.
\end{equation}
\end{restatable}

The lower bound holds for every layout and every closed visiting order.

\Needspace{4\baselineskip}
\begin{restatable}[Universal lower bound]{theorem}{universalLowerBound}\label{thm:lower}
For each fixed $r\ge 0$ there is a finite constant $K_r\ge 0$, independent of $N$, such that for every
integer $N\ge 3$, every set $X$ of $N$ distinct facilities in $Q_N$, and every closed tour $\sigma$ through them,
\begin{equation}\label{eq:thm2}
M_N(X)+r\,L(X,\sigma)\;\ge\;N\,C_6(r)-K_r .
\end{equation}
In particular $A_N(r)\ge N C_6(r)-K_r$. The constant is given in \eqref{eq:final}, and $K_0=0$.
\end{restatable}

The bounds have the common leading term $NC_6(r)$. The correction $D_r\sqrt N$ comes from the finite-square construction, and $K_r$ comes from the disk-union estimate. Combining the bounds gives the following result.

\Needspace{4\baselineskip}
\begin{restatable}[Finite-size gap and limiting cost]{corollary}{finiteSizeGap}\label{cor:gap}
Fix $r\ge0$. Use $D_r$ and the size threshold from Theorem~\ref{thm:upper}, and use $K_r$ from Theorem~\ref{thm:lower}. For every $N$ above that threshold,
\begin{equation}\label{eq:gapboth}
NC_6(r)-K_r\le A_N(r)\le NC_6(r)+D_r\sqrt N.
\end{equation}
In particular, the absolute error per facility satisfies
\begin{equation}\label{eq:gapabs}
\left|\frac{A_N(r)}N-C_6(r)\right|\le\max\left\{\frac{K_r}N,\frac{D_r}{\sqrt N}\right\},
\end{equation}
and
\begin{equation}\label{eq:limit}
\lim_{N\to\infty}\frac{A_N(r)}N=C_6(r).
\end{equation}
\end{restatable}

\begin{proof}
Theorems~\ref{thm:upper} and~\ref{thm:lower} give \eqref{eq:gapboth}. Dividing by $N$ and subtracting $C_6(r)$ gives
\[
-\frac{K_r}{N}\le\frac{A_N(r)}N-C_6(r)\le\frac{D_r}{\sqrt N}.
\]
Taking absolute values yields \eqref{eq:gapabs}. Both constants are finite and independent of $N$ for fixed $r$, so \eqref{eq:limit} follows.
\end{proof}

The same estimates prove asymptotic optimality of the constructed layout and tour. Choose a minimizing spacing $d_r$ and let $(X,\sigma)$ be the layout and visiting order of Theorem~\ref{thm:upper}. Then
\[
0\le M_N(X)+rL(X,\sigma)-A_N(r)\le D_r\sqrt N+K_r.
\]
After division by $N$, the construction's excess cost per facility tends to zero for fixed $r$.

\begin{figure}[!htbp]\centering
\begin{adjustbox}{max width=\linewidth}
\begin{tikzpicture}[box/.style={draw=black!55,rounded corners=1pt,align=center,text width=6.7cm,inner sep=6pt,font=\small},arr/.style={-{Stealth[length=1.5mm]},font=\small}]
\node[font=\small\bfseries] at (0,0) {UPPER BOUND};
\node[font=\small\bfseries] at (7.5,0) {LOWER BOUND};
\node[font=\scriptsize] at (0,-0.38) {Section~\ref{sec:upper}, Appendix~\ref{app:upper}};
\node[font=\scriptsize] at (7.5,-0.38) {Sections~\ref{sec:lb1}--\ref{sec:lb2}, Appendices~\ref{app:local}--\ref{app:global}};
\node[box,anchor=north] (u1) at (0,-0.68) {Cell geometry: Lemma~\ref{lem:hexagon}\\Minimizing spacing: Lemma~\ref{lem:existence}\\Xie--Ouyang construction};
\node[box,anchor=north] (l1) at (7.5,-0.68) {Section~\ref{sec:lb1}: tour and side count\\Relaxation and Lagrange lower bound\\Section~\ref{sec:lb2}: construct and verify multipliers};
\node[box,anchor=north] (u2) at (0,-2.65) {Finite-square construction: Theorem~\ref{thm:upper}\\Place rows, bound service, close the tour.\\Additional cost $\le D_r\sqrt N$.};
\node[box,anchor=north] (l2) at (7.5,-2.65) {Section~\ref{sec:completion}: complete the proof\\Trade-off and boundary extension\\Choose the spacing for the routing weight};
\node[box,fill=black!3,anchor=north] (u3) at (0,-4.60) {\textbf{Theorem~\ref{thm:upper}}\\$A_N(r)\le NC_6(r)+D_r\sqrt N$};
\node[box,fill=black!3,anchor=north] (l3) at (7.5,-4.60) {\textbf{Theorem~\ref{thm:lower}}\\$A_N(r)\ge NC_6(r)-K_r$};
\node[box,text width=13.9cm,anchor=north] (bot) at (3.75,-6.05) {\textbf{Corollary~\ref{cor:gap}: finite-size gap and limit}\quad Section~\ref{sec:mainresults}\\[3pt]
$\displaystyle\left|\frac{A_N(r)}N-C_6(r)\right|\le\max\left\{\frac{K_r}N,\frac{D_r}{\sqrt N}\right\}\longrightarrow0$\\[3pt]$r$ fixed, $\lim_{N\to\infty}A_N(r)/N=C_6(r)$.};
\draw[arr] (u1.south)--(u2.north);\draw[arr] (l1.south)--(l2.north);
\draw[arr] (u2.south)--(u3.north);\draw[arr] (l2.south)--(l3.north);
\draw[arr] (u3.south)--(bot.north -| u3.south);\draw[arr] (l3.south)--(bot.north -| l3.south);
\end{tikzpicture}
\end{adjustbox}
\caption{Outline of the proof. The upper bound uses the hexagonal construction of \citet{XieOuyang2015}. The lower bound applies to arbitrary layouts and tours. The two bounds imply Corollary~\ref{cor:gap}.}\label{fig:dependencies}
\end{figure}
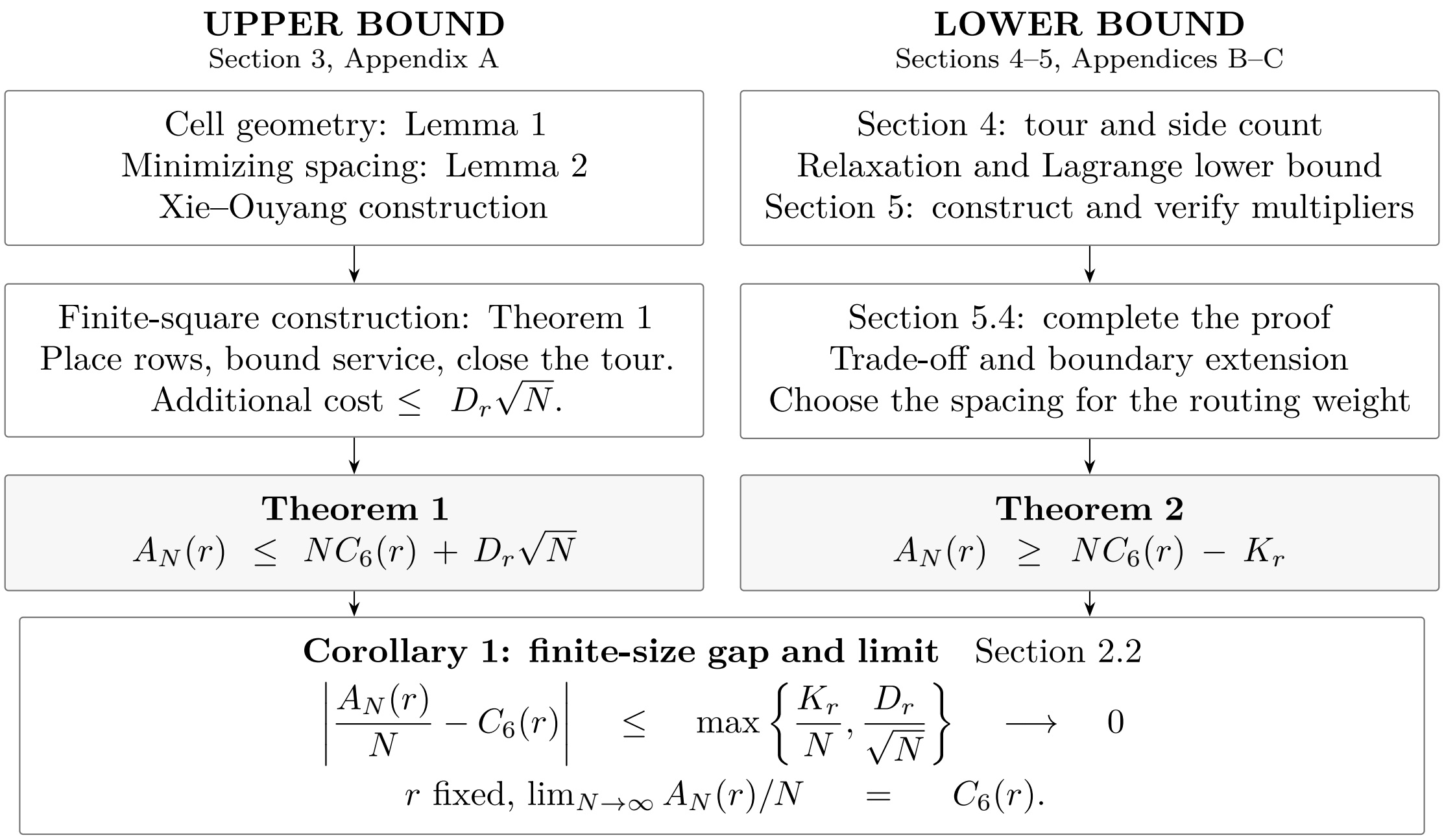
\section{Proof of Theorem 1: the finite-square construction}\label{sec:upper}

We implement the hexagonal pattern in $Q_N$, using the cell formula for its interior service cost. We then bound service near the boundary, connect the rows, and complete the facility count. Each correction costs at most a constant times $\sqrt N$ for fixed $d$. Appendix~\ref{app:upperconstruction} gives the complete construction and constants.

\subsection{Cell geometry}\label{sec:cellgeometry}
Let $a$ and $b$ be the distances from the facility to the two vertical and four sloping sides of $H_d$, respectively, and let $R$ be the common distance to its vertices:
\begin{equation}\label{eq:abR}
a=\frac d2,\qquad b=\frac12\sqrt{d^{-2}+d^2/4},\qquad R=\frac1{2d}+\frac{d^3}{8}.
\end{equation}
For $0<d\le d_0$, we have $0<a\le b<R$. We call the two sides at distance $a$ the \emph{near} sides and the four at distance $b$ the \emph{far} sides. Splitting the cell into six triangles gives its service cost.

\Needspace{4\baselineskip}
\begin{restatable}[The unit-area hexagon]{lemma}{unitAreaHexagon}\label{lem:hexagon}
The service cell is
\begin{equation}\label{eq:hexagon}
H_d=\{(x,y):|x|\le d/2,\ |y|\le R-(d^2/2)|x|\}.
\end{equation}
Its six vertices are at distance $R$ from the facility, two of its sides are at distance $a$, and four are at distance $b$. For $0<h<R$, the integral of $\|x\|$ over the triangle with vertices $(0,0)$ and $(h,\pm\sqrt{R^2-h^2})$ is
\begin{equation}\label{eq:JR}
J_R(h)=\frac{hR\sqrt{R^2-h^2}+h^3\arcosh(R/h)}{3}.
\end{equation}
Moreover,
\begin{equation}\label{eq:phiformula}
|H_d|=1,\qquad \phi(d)=2J_R(a)+4J_R(b).
\end{equation}
\end{restatable}
The proof is in Appendix~\ref{app:hexagon}.

Figure~\ref{fig:referencehex} shows the two near-side and four far-side triangles. We use $H_d$ as the \emph{reference hexagon} in the lower-bound proof.

\begin{figure}[!htbp]\centering
\begingroup
\resizebox{0.79\textwidth}{!}{
\begin{tikzpicture}[color=black,>=Stealth,x=4.8cm,y=4.8cm,font=\small,
every node/.style={text=black}]
\pgfmathsetmacro{\fd}{0.9}
\pgfmathsetmacro{\fa}{\fd/2}
\pgfmathsetmacro{\fr}{1/(2*\fd)+\fd^3/8}
\pgfmathsetmacro{\fv}{1/\fd-\fr}
\draw[black!65,dashed,thin] (0,0) circle (\fr);
\fill[black!16] (0,0)--(\fa,\fv)--(\fa,-\fv)--cycle;
\fill[black!16] (0,0)--(-\fa,\fv)--(-\fa,-\fv)--cycle;
\foreach \x/\y in {0/\fr,\fa/\fv,\fa/-\fv,0/-\fr,-\fa/-\fv,-\fa/\fv}{
\draw[black!55,thin] (0,0)--(\x,\y);
\fill[black] (\x,\y) circle (0.008);
}
\draw[thick] (0,\fr)--(\fa,\fv)--(\fa,-\fv)--(0,-\fr)--(-\fa,-\fv)--(-\fa,\fv)--cycle;
\draw[->,thin] (0,0)--(-\fa,\fv) node[pos=0.55,above right,inner sep=2pt] {$R$};
\draw[->,thin] (0,0)--({\fa/2},{(\fr+\fv)/2}) node[pos=0.55,left,inner sep=2pt] {$b$};
\draw[<->,thin] (0,-0.10)--(\fa,-0.10) node[midway,below,inner sep=2pt] {$a$};
\draw[black!55,thin] (0,0)--(0,-0.14);
\fill[black] (0,0) circle (0.014);
\node[below left,inner sep=2pt] at (0,0) {$O$};
\node[above,inner sep=3pt] at (0,\fr) {$(0,R)$};
\node[anchor=west,inner sep=3pt,font=\footnotesize] at (\fa,\fv) {$(a,1/d-R)$};
\node[below,inner sep=5pt] at (0,-\fr) {$|H_d|=1$};
\node[anchor=west,align=left,text=black,font=\small] at (0.96,0.02)
{Here, $a=\dfrac d2$\\[9pt]
$b=\dfrac12\sqrt{d^{-2}+d^2/4}$\\[9pt]
$R=\dfrac1{2d}+\dfrac{d^3}{8}$};
\end{tikzpicture}
}
\endgroup
\caption{The six triangles obtained by joining the facility to the vertices of $H_d$. Two have a side at distance $a$ from the facility, four have a side at distance $b$, and all vertices of $H_d$ are at distance $R$. Their service costs sum to $\phi(d)=2J_R(a)+4J_R(b)$.}\label{fig:referencehex}
\end{figure}

\subsection{Choice of spacing}\label{sec:existence}
The service cost diverges as $d$ approaches zero, preventing an optimal spacing from collapsing to zero. To see this, compare $H_d$ with the centered rectangle
\[
[-d/2,d/2]\times[-1/(2d),1/(2d)].
\]
Both regions have area one and lie in the strip $|x|\le d/2$. Their set differences therefore have equal area. On the part of the hexagon outside the rectangle, $|y|\ge1/(2d)$. On the part of the rectangle outside the hexagon, $|y|\le1/(2d)$. Their common part contributes equally to both integrals, so the integral of $|y|$ over the hexagon is at least that over the rectangle. Since $\sqrt{x^2+y^2}\ge |y|$, we have
\[
\phi(d)\ge\int_{H_d}|y|\,\dd x\dd y
\ge d\int_{-1/(2d)}^{1/(2d)}|y|\,\dd y=\frac1{4d}.
\]
This lower bound diverges as $d\downarrow0$. Together with continuity of $\phi$, it gives the following result.

\Needspace{4\baselineskip}
\begin{restatable}[Existence of a minimizing spacing]{lemma}{minimizingSpacingExists}\label{lem:existence}
For $0<d\le d_0$, $\phi(d)\ge1/(4d)$. For every fixed $r\ge0$, the minimum in \eqref{eq:Csix} is attained at a positive spacing.
\end{restatable}
The proof is in Appendix~\ref{app:existence}.

Thus a minimizing spacing is chosen once for the given $r$ and is kept fixed as $N$ increases.

\subsection{Finite-square construction}\label{sec:finiteconstruction}

Write $s=\sqrt N$, $m=\lfloor s/d\rfloor$, and $n=\lfloor sd\rfloor$. Place $m$ facilities in each of $n$ staggered rows, and put $q=N-mn$. Then $q\le(d+d^{-1})s$. For each of these $mn$ facilities, take its hexagonal cell in the translated infinite lattice, and let $W$ be the union of these full cells. They have total service integral $mn\phi(d)$, although $W$ may extend outside the square. The region $Q_N\setminus W$ has area $O(s)$, and every customer there is within a distance depending only on $d$ of one of these facilities. Its service cost is therefore $O(s)$.

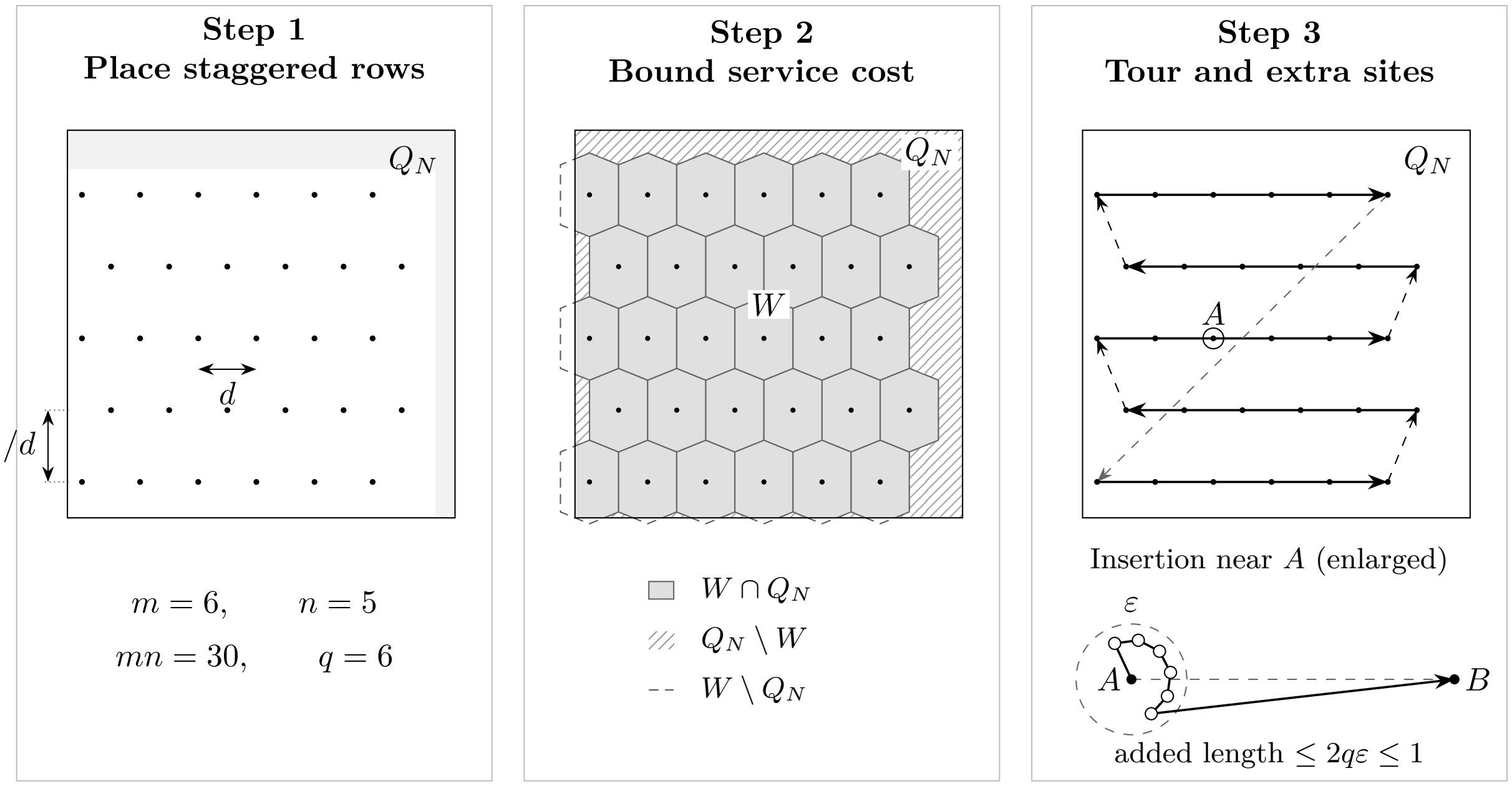
\begin{figure}[!htbp]\centering
\begingroup
\resizebox{\textwidth}{!}{
\begin{tikzpicture}[color=black,>=Stealth,font=\small,
every node/.style={text=black},frame/.style={black!25,thin}]
\pgfmathsetmacro{\fd}{0.9}
\pgfmathsetmacro{\fa}{\fd/2}
\pgfmathsetmacro{\fr}{1/(2*\fd)+\fd^3/8}
\pgfmathsetmacro{\fv}{1/\fd-\fr}

\foreach \p/\title in {0/{Place staggered rows},1/{Bound service cost},2/{Tour and extra sites}}{
\begin{scope}[xshift={5.5*\p cm}]
\draw[frame] (0,0) rectangle (5.15,8.4);
\pgfmathtruncatemacro{\stepno}{\p+1}
\node[align=center,font=\small\bfseries] at (2.575,7.91) {Step \stepno\\\title};
\end{scope}
}

\begin{scope}[shift={(0.55,2.85)},scale=0.70]
\fill[black!5] (0,5.4) rectangle (6,6);
\fill[black!5] (5.7,0) rectangle (6,6);
\draw[thin] (0,0) rectangle (6,6);
\node[anchor=north east,font=\small] at (5.94,5.94) {$Q_N$};
\foreach \k in {0,...,4}{\foreach \j in {0,...,5}{
\pgfmathsetmacro{\px}{\fd*(\j+0.25+mod(\k,2)/2)}
\pgfmathsetmacro{\py}{(\k+0.5)/\fd}
\fill[black] (\px,\py) circle (0.045);
}}
\draw[<->,thin] (2.025,2.30)--(2.925,2.30) node[midway,below] {$d$};
\draw[black!60,densely dotted,thin] (-0.35,{0.5/\fd})--(0.05,{0.5/\fd});
\draw[black!60,densely dotted,thin] (-0.35,{1.5/\fd})--(0.05,{1.5/\fd});
\draw[<->,thin] (-0.30,{0.5/\fd})--(-0.30,{1.5/\fd}) node[midway,left] {$1/d$};
\end{scope}
\node at (2.575,1.90) {$m=6,\qquad n=5$};
\node at (2.575,1.32) {$mn=30,\qquad q=6$};

\begin{scope}[shift={(6.05,2.85)},scale=0.70]
\fill[pattern=north east lines,pattern color=black!35] (0,0) rectangle (6,6);
\foreach \k in {0,...,4}{\foreach \j in {0,...,5}{
\pgfmathsetmacro{\px}{\fd*(\j+0.25+mod(\k,2)/2)}
\pgfmathsetmacro{\py}{(\k+0.5)/\fd}
\begin{scope}[shift={(\px,\py)}]
\draw[densely dashed,thin,black!65]
(0,\fr)--(\fa,\fv)--(\fa,-\fv)--(0,-\fr)--(-\fa,-\fv)--(-\fa,\fv)--cycle;
\end{scope}
}}
\begin{scope}
\clip (0,0) rectangle (6,6);
\foreach \k in {0,...,4}{\foreach \j in {0,...,5}{
\pgfmathsetmacro{\px}{\fd*(\j+0.25+mod(\k,2)/2)}
\pgfmathsetmacro{\py}{(\k+0.5)/\fd}
\begin{scope}[shift={(\px,\py)}]
\filldraw[fill=black!12,draw=black!60,thin]
(0,\fr)--(\fa,\fv)--(\fa,-\fv)--(0,-\fr)--(-\fa,-\fv)--(-\fa,\fv)--cycle;
\end{scope}
}}
\end{scope}
\foreach \k in {0,...,4}{\foreach \j in {0,...,5}{
\pgfmathsetmacro{\px}{\fd*(\j+0.25+mod(\k,2)/2)}
\pgfmathsetmacro{\py}{(\k+0.5)/\fd}
\fill[black] (\px,\py) circle (0.04);
}}
\draw[thin] (0,0) rectangle (6,6);
\node[anchor=north east,font=\small,fill=white,inner sep=1pt] at (5.94,5.94) {$Q_N$};
\node[fill=white,inner sep=1pt] at (3,3.3) {$W$};
\end{scope}
\filldraw[fill=black!12,draw=black!60] (6.85,1.97) rectangle (7.12,2.16);
\node[anchor=west,font=\footnotesize] at (7.28,2.065) {$W\cap Q_N$};
\fill[pattern=north east lines,pattern color=black!40] (6.85,1.43) rectangle (7.12,1.62);
\node[anchor=west,font=\footnotesize] at (7.28,1.525) {$Q_N\setminus W$};
\draw[densely dashed,black!65] (6.85,0.985)--(7.12,0.985);
\node[anchor=west,font=\footnotesize] at (7.28,0.985) {$W\setminus Q_N$};

\begin{scope}[shift={(11.55,2.85)},scale=0.70]
\draw[thin] (0,0) rectangle (6,6);
\node[anchor=north east,font=\small] at (5.94,5.94) {$Q_N$};
\foreach \k in {0,...,4}{\foreach \j in {0,...,5}{
\pgfmathsetmacro{\px}{\fd*(\j+0.25+mod(\k,2)/2)}
\pgfmathsetmacro{\py}{(\k+0.5)/\fd}
\coordinate (p\k-\j) at (\px,\py);
}}
\foreach \k in {0,...,4}{
\pgfmathtruncatemacro{\parity}{mod(\k,2)}
\ifnum\parity=0
\draw[->,line width=0.7pt] (p\k-0)--(p\k-5);
\else
\draw[->,line width=0.7pt] (p\k-5)--(p\k-0);
\fi
}
\draw[->,dashed,thin] (p0-5)--(p1-5);
\draw[->,dashed,thin] (p1-0)--(p2-0);
\draw[->,dashed,thin] (p2-5)--(p3-5);
\draw[->,dashed,thin] (p3-0)--(p4-0);
\draw[->,dashed,thin,black!60] (p4-5)--(p0-0);
\foreach \k in {0,...,4}{\foreach \j in {0,...,5}{\fill[black] (p\k-\j) circle (0.045);}}
\draw[thin] (p2-2) circle (0.16);
\node[above,font=\small,inner sep=4pt] at (p2-2) {$A$};
\end{scope}
\begin{scope}[xshift=11cm]
\node[font=\footnotesize] at (2.575,2.38) {Insertion near $A$ (enlarged)};
\coordinate (ia) at (1.08,1.10);
\coordinate (ib) at (4.58,1.10);
\draw[dashed,thin,black!60] (ia) circle (0.60);
\draw[dashed,thin,black!60] (ia)--(ib);
\node[above,font=\small] at (1.08,1.70) {$\varepsilon$};
\foreach \n/\angle in {1/115,2/80,3/45,4/10,5/-25,6/-60}{
\coordinate (z\n) at ($(ia)+(\angle:0.43)$);
}
\draw[->,line width=0.7pt] (ia)--(z1)--(z2)--(z3)--(z4)--(z5)--(z6)--(ib);
\foreach \n in {1,...,6}{\filldraw[fill=white,draw=black,thin] (z\n) circle (0.065);}
\fill[black] (ia) circle (0.055); \fill[black] (ib) circle (0.055);
\node[left,inner sep=3pt] at (ia) {$A$};
\node[right,inner sep=3pt] at (ib) {$B$};
\node[font=\footnotesize,text=black] at (2.575,0.28) {added length $\le 2q\varepsilon\le1$};
\end{scope}
\end{tikzpicture}
}
\endgroup
\caption{Finite-square construction. Place staggered rows, bound service in the boundary region, and connect the rows into a closed tour. Additional sites are placed near one of the row facilities. The example uses $N=36$, $d=0.9$, and 30 sites in the rows. The insertion detail is enlarged.}\label{fig:finiteconstruction}
\end{figure}

Visit the first row from left to right, then alternate directions
in successive rows, joining rows at the same end. Connect the last
facility to the first to close the tour. There are $O(s)$ row connections, each of length bounded in terms of $d$,
and the closing edge is no longer than the square's diagonal. If $q>0$, choose an interior facility $A$ in these rows and let $B$ be its next stop. If $mn=1$, take $B=A$.
Place the $q$ new facilities at distinct points other than $A$ in a ball
of radius $\varepsilon\le1/(2N)$ around $A$, contained in the square and
free of other row facilities. Replace $AB$ by a path through these
points, adding at most $2q\varepsilon\le1$ to the tour length. No
nearest-facility distance increases, since no row facility moves.
Figure~\ref{fig:finiteconstruction} illustrates the construction.
Appendix~\ref{app:upperconstruction} gives the full proof and constants.

\section{Proof of Theorem 2: the relaxation and lower bound}\label{sec:lb1}

We derive a lower bound from four aggregate constraints on a layout and tour: total area, total angle, side count, and disk overhang. Section~\ref{sec:lb2} constructs Lagrange multipliers for which the bound on each side triangle is attained by the reference triangles.

Throughout Sections~\ref{sec:lb1} and~\ref{sec:lb2}, $X$ is a set of $N\ge 3$ distinct facilities in
$Q_N$, visited by a closed tour. Its tour-edge lengths are $\ell_1,\dots,\ell_N$, and its total length is $L=\sum_i\ell_i$. We
fix a reference spacing $d\in(0,d_0]$ and write $a\le b<R$ for the two side distances and the vertex
radius of the reference hexagon $H_d$, as in \eqabR, with $a=d/2$. The competing layout $X$ and its tour need not follow the reference pattern. Service regions may have different areas and shapes, and $L$ may be smaller or larger than $Nd$. The reference parameters and the multipliers constructed in Section~\ref{sec:lb2} depend only on
$d$, not on $N$ or the competing layout.

We first assume that the facilities lie in the interior of $Q_N$. This gives a positive distance from each facility to every side of its cell and a full angle of $2\pi$ around that facility. Lemma~\ref{lem:boundary} removes this temporary assumption by a limiting argument. Dividing each service cell into triangles from its facility gives at most $6N$ triangles, with total area $N$ and total angle $2\pi N$ at the facilities. A fourth constraint links their total disk overhang---the disk area beyond their opposite sides---to the overlaps of disks centered at consecutive tour stops. Sections~\ref{sec:uncovered} and~\ref{sec:triangles} prove these facts.

\subsection{The inequality to be proved}\label{sec:goal}

In the hexagonal layout with spacing $d$, each facility accounts for service cost $\phi(d)$ and tour
length $d$. A competing layout may have a shorter tour, $L<Nd$. We bound its service-cost difference in terms of the saving in tour length, up to a constant independent of $N$:
\begin{equation}\label{eq:tradeoff}
M_N(X)-N\phi(d)\;\ge\;\rs(d)\,(Nd-L)-K(d).
\end{equation}
Here $\rs(d)$ is the service-cost penalty per unit of tour length saved in this bound, and $K(d)$ is an allowance independent of $N$. Both are nonnegative and are defined in \eqref{eq:rK}. Adding $\rs(d)L$ to both sides removes the tour length from the right-hand side:
\begin{equation}\label{eq:tradeoff-total}
M_N(X)+\rs(d)\,L\;\ge\;N\{\phi(d)+\rs(d)\,d\}-K(d).
\end{equation}
To obtain Theorem~\ref{thm:lower}, choose $d$ to minimize $\phi(d)+rd$. Lemma~\ref{lem:spacing} shows that this choice also gives $\rs(d)=r$ through the identity
\begin{equation}\label{eq:rstar-phi}
\rs(d)=-\phi'(d).
\end{equation}

\subsection{Uncovered area and the tour}\label{sec:uncovered}\label{sec:tour}

\begin{figure}[t]
\centering
\begin{adjustbox}{max width=\linewidth}
\begin{tikzpicture}[scale=1.1]

\def\pts{(0.5,0.55),(1.45,0.4),(2.5,0.6),(0.4,1.5),(1.5,1.55),(2.55,1.45),(0.6,2.5),(1.55,2.6),(2.45,2.4)}
\fill[black!30] (0,0) rectangle (3,3);
\begin{scope}
\clip (0,0) rectangle (3,3);
\foreach \p in \pts {\fill[white] \p circle (0.55);}
\foreach \p in \pts {\draw[acc,thin] \p circle (0.55);}
\end{scope}
\foreach \p in \pts {\node[dot,minimum size=2pt] at \p {};}
\draw (0,0) rectangle (3,3);
\draw[->,thin] (1.5,1.55) -- ++(40:0.55);
\node[lab] at (1.84,1.62) {$t$};
\fill[black!30] (0.3,-0.44) rectangle (0.55,-0.24);
\node[lab,anchor=west] at (0.57,-0.34) {uncovered, area $U_X(t)$};
\node[lab,font=\scriptsize\bfseries] at (1.5,3.3) {(a) A layout in $Q_N$};

\begin{scope}[shift={(5.9,1.5)},scale=2.3]
\fill[soft] (0,0.64668) -- (0.45,0.46443) -- (0.45,-0.46443) -- (0,-0.64668) -- (-0.45,-0.46443) -- (-0.45,0.46443) -- cycle;
\begin{scope}\clip (0.45,-0.7) rectangle (0.75,0.7); \fill[capc] (0,0) circle (0.53);\end{scope}
\begin{scope}\clip (-0.45,-0.7) rectangle (-0.75,0.7); \fill[capc] (0,0) circle (0.53);\end{scope}
\draw[acc] (0,0) circle (0.53);
\draw[thick] (0,0.64668) -- (0.45,0.46443) -- (0.45,-0.46443) -- (0,-0.64668) -- (-0.45,-0.46443) -- (-0.45,0.46443) -- cycle;
\draw[thin] (0,0) -- (0.45,0) node[midway,below=-1pt,lab]{$a$};
\draw[thin] (0,0) -- (0.225,0.55556) node[pos=0.55,left=-1pt,lab]{$b$};
\draw[->,thin] (0,0) -- (225:0.53) node[pos=0.55,above left=-2pt,lab]{$t$};
\node[dot] at (0,0) {};
\node[lab] at (0,-0.84) {two caps, each of area $c_a(t)$};
\node[lab] at (0.4,0.63) {$H_d$};
\end{scope}
\node[lab,font=\scriptsize\bfseries] at (5.9,3.3) {(b) The reference cell, $a\le t\le b$};

\begin{scope}[shift={(10.75,1.5)},scale=2.3]
\fill[soft] (-0.45,0) circle (0.53);
\fill[soft] (0.45,0) circle (0.53);
\begin{scope}\clip (-0.45,0) circle (0.53); \fill[capc] (0.45,0) circle (0.53);\end{scope}
\draw[acc] (-0.45,0) circle (0.53);
\draw[acc] (0.45,0) circle (0.53);
\draw[densely dashed,thin] (0,-0.36) -- (0,0.36);
\draw[thin] (-0.45,0) -- (0.45,0);
\node[dot] at (-0.45,0) {}; \node[dot] at (0.45,0) {};
\node[lab,below] at (-0.45,-0.01) {$x_i$};
\node[lab,below] at (0.47,-0.01) {$x_{i+1}$};
\draw[<->,thin] (-0.45,-0.64) -- (0.45,-0.64) node[midway,below,lab]{$d=2a$};
\node[lab] at (0,-0.98) {lens of area $o_t(d)=2c_a(t)$};
\end{scope}
\node[lab,font=\scriptsize\bfseries] at (10.75,3.3) {(c) Two consecutive stops};
\end{tikzpicture}
\end{adjustbox}
\caption{Uncovered area and overlap. (a) $U_X(t)$ is the area of $Q_N$ outside all disks of radius $t$.
(b) For $a\le t\le b$, only the two near-side caps can lie outside $H_d$. Their areas are zero at $t=a$. (c) Two consecutive disks of the hexagonal layout overlap in a lens made
of the same two caps. Panels (b) and (c) use $d=0.9$ and $t=0.53$.}
\label{fig:uncovered}
\end{figure}
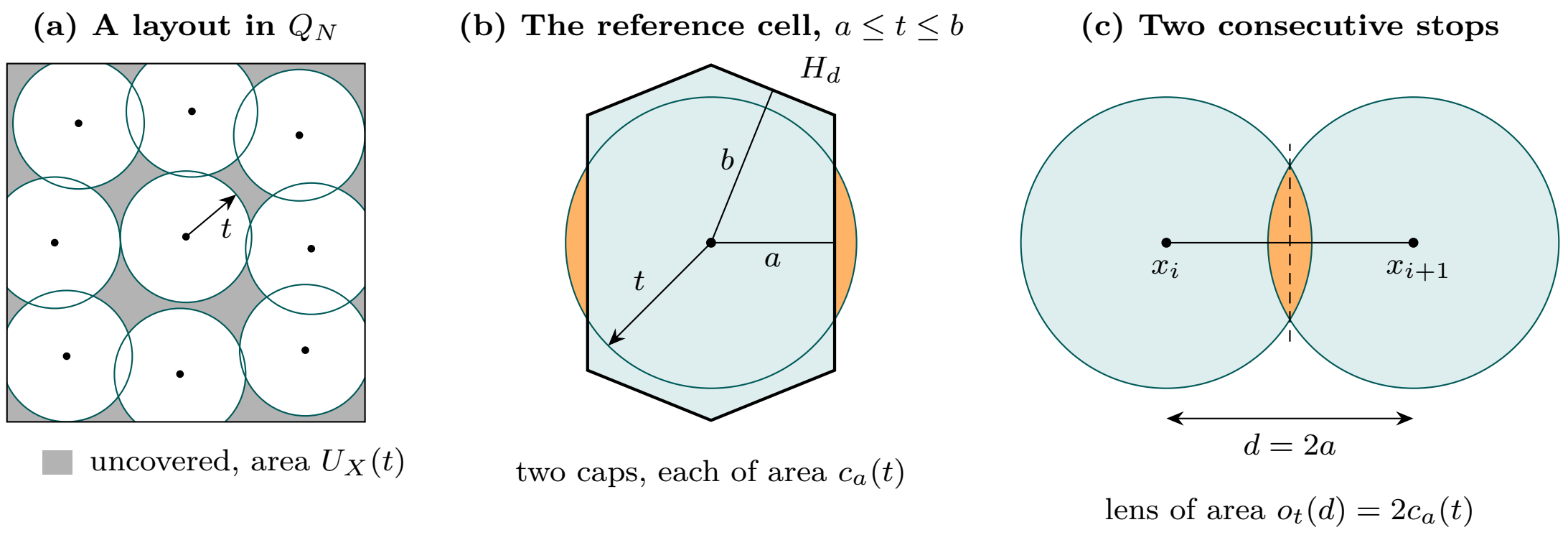

The quantity that connects the tour to the service cost is the area that the facilities leave
uncovered at a given distance. For $t\ge 0$ let $B_t(x)$ be the closed disk of radius $t$ centered at
$x$, and define the \emph{uncovered area}
\begin{equation}\label{eq:U}
U_X(t)=\bigl|\{x\in Q_N:\dist(x,X)>t\}\bigr|
=N-\Bigl|\,Q_N\cap\bigcup_{i=1}^{N}B_t(x_i)\Bigr| .
\end{equation}
Thus $U_X(t)$ measures the part of the square whose customers are farther than $t$ from every facility (Figure~\ref{fig:uncovered}a). We use $t$ as an auxiliary radius. Section~\ref{sec:prices} integrates the coverage constraint over $a\le t\le b$ with a nonnegative weight.

For a customer at $x$, the distance $D(x)=\dist(x,X)$ equals $\int_0^\infty\mathbf 1_{\{t<D(x)\}}\dd t$. The indicator is measurable and nonnegative, so Tonelli's theorem \citep[Theorem~1.7.18]{Tao2011} permits us to reverse the order of integration:
\[
M_N(X)=\int_{Q_N}\int_0^\infty\mathbf 1_{\{t<D(x)\}}\dd t\dd x
=\int_0^\infty U_X(t)\dd t.
\]

For the reference hexagon, with its facility at the origin, write $u_d(t)=\abs{H_d\setminus B_t(0)}$.
Two sides of $H_d$ are at distance $a$ from the
facility, and four are at distance $b\ge a$. As $t$ grows, the disk first lies inside the hexagon, then
crosses the two near sides, and only when $t>b$ crosses the other four. For $a\le t\le b$, the only possible exterior parts are the two near-side caps, whose areas are zero at $t=a$ (Figure~\ref{fig:uncovered}b). Let $c_h(t)$
denote the area of the cap that a line at distance $h$ from the center cuts off a disk of radius $t$:
\begin{equation}\label{eq:cap}
c_h(t)=\begin{cases}
t^2\arccos(h/t)-h\sqrt{t^2-h^2}, & 0\le h<t,\\[2pt]
0, & h\ge t.
\end{cases}
\end{equation}
For $h<t$, the cutting line meets the disk in a chord of length $2\sqrt{t^2-h^2}$. Moving the line outward decreases the cap area at this rate, so
\begin{equation}\label{eq:capderiv}
\partial_h c_h(t)=-2\sqrt{(t^2-h^2)_+},\qquad h>0,\quad t\ge0.
\end{equation}
The derivative is zero for $h\ge t$, including at $h=t$. The calculation is given in \appOverlap.

For $a\le t\le b$, the area inside both the disk and the reference cell is $\pi t^2-2c_a(t)$. These caps also describe overlap along the tour: for two disks of radius $t$ whose centers are $\ell$ apart, the overlap area is the sum of two cap areas (Figure~\ref{fig:uncovered}c):
\begin{equation}\label{eq:overlap}
o_t(\ell)=2\,c_{\ell/2}(t).
\end{equation}
For fixed $t$, the lens area is nonincreasing and convex in the separation $\ell$. Its derivative
$o_t'(\ell)=-\sqrt{(4t^2-\ell^2)_+}$ is nondecreasing (\appOverlap). In the reference layout, consecutive within-row facilities are $d=2a$ apart. Hence $2c_a(t)=o_t(d)$ and
\begin{equation}\label{eq:ud}
u_d(t)=1-\pi t^2+o_t(d),\qquad a\le t\le b .
\end{equation}
Thus, for these radii, the part of a disk outside its own cell has the same area as its overlap with the next disk on the tour.

The sum $N\pi t^2$ overestimates the area covered by overlapping disks. The tour selects $N$ pairs of consecutive disks, whose lens areas are $o_t(\ell_i)$. Subtracting these overlaps gives the following estimate.

\Needspace{4\baselineskip}
\begin{restatable}[Tour constraint]{lemma}{tourConstraint}\label{lem:tour}
For every $t\ge 0$, we have
\begin{equation}\label{eq:tour}
U_X(t)\;\ge\;N(1-\pi t^2)+\sum_{i=1}^{N}o_t(\ell_i)-\pi t^2 .
\end{equation}
\end{restatable}

Only consecutive-stop lenses are subtracted. The allowance $\pi t^2$ accounts for points that lie in all $N$ disks, as the proof below shows.

\begin{figure}[!htbp]
\centering
\begin{adjustbox}{max width=\linewidth}
\begin{tikzpicture}[scale=0.78]
\coordinate (A) at (-1,0); \coordinate (B) at (0,1.3); \coordinate (C) at (1,0);
\coordinate (D) at (1.5,-2.2); \coordinate (E) at (-1.5,-2.2); \coordinate (Y) at (0,0.35);

\foreach \P in {A,B,C} {\fill[soft] (\P) circle (1.25);}
\foreach \P in {A,B,C} {\draw[acc] (\P) circle (1.25);}
\foreach \P in {D,E} {\draw[black!45,densely dashed] (\P) circle (1.25);}
\draw[black!55,thin] (A) -- (B) -- (C) -- (D) -- (E) -- cycle;
\foreach \P/\pos in {A/left,B/above,C/right} {\node[dot] at (\P) {}; \node[lab,\pos] at (\P) {$\P$};}
\foreach \P/\pos in {D/below right,E/below left} {\node[odot] at (\P) {}; \node[lab,\pos] at (\P) {$\P$};}
\draw[thick] ($(Y)+(-0.09,-0.09)$) -- ($(Y)+(0.09,0.09)$) ($(Y)+(-0.09,0.09)$) -- ($(Y)+(0.09,-0.09)$);
\node[lab,right] at ($(Y)+(0.05,-0.1)$) {$y$};
\draw[->,thin] (B) -- ++(150:1.25) node[midway,above=-1pt,lab]{$t$};
\node[lab,font=\scriptsize\bfseries] at (0,3.05) {(a) Disks of radius $t$ around the five stops};
\node[lab] at (0,-4.0) {$y$ lies in the disks around $A$, $B$ and $C$, so $k(y)=3$};

\begin{scope}[shift={(9.2,0)}]
\coordinate (A) at (-1,0); \coordinate (B) at (0,1.3); \coordinate (C) at (1,0);
\coordinate (D) at (1.5,-2.2); \coordinate (E) at (-1.5,-2.2); \coordinate (Y) at (0,0.35);
\draw[black!55,thin,densely dashed] (C) -- (D) -- (E) -- (A);
\draw[acc,line width=1.6pt] (A) -- (B) -- (C);
\foreach \P/\pos in {A/left,B/above,C/right} {\node[dot,minimum size=3.6pt] at (\P) {}; \node[lab,\pos] at (\P) {$\P$};}
\foreach \P/\pos in {D/below right,E/below left} {\node[odot] at (\P) {}; \node[lab,\pos] at (\P) {$\P$};}
\draw[thick] ($(Y)+(-0.09,-0.09)$) -- ($(Y)+(0.09,0.09)$) ($(Y)+(-0.09,0.09)$) -- ($(Y)+(0.09,-0.09)$);
\node[lab,right] at ($(Y)+(0.05,-0.1)$) {$y$};
\node[lab,font=\scriptsize\bfseries] at (0,3.05) {(b) The tour $A\to B\to C\to D\to E\to A$};
\node[lab] at (0,-3.55) {only $AB$ and $BC$ join two marked stops,};
\node[lab] at (0,-4.0) {so $m(y)=2$ and $k(y)-m(y)=1$};
\end{scope}
\end{tikzpicture}
\end{adjustbox}
\caption{The count behind Lemma~\ref{lem:tour}. The point $y$ is counted three times in the total disk
area and twice in the lenses of consecutive stops ($AB$ and $BC$), so once in the difference. It
also lies in the disks around $A$ and $C$, but $AC$ is not a tour edge and that overlap is not
subtracted. Here $k(y)$ counts disks containing $y$, and $m(y)$ counts consecutive-stop lenses containing $y$.}
\label{fig:count}
\end{figure}
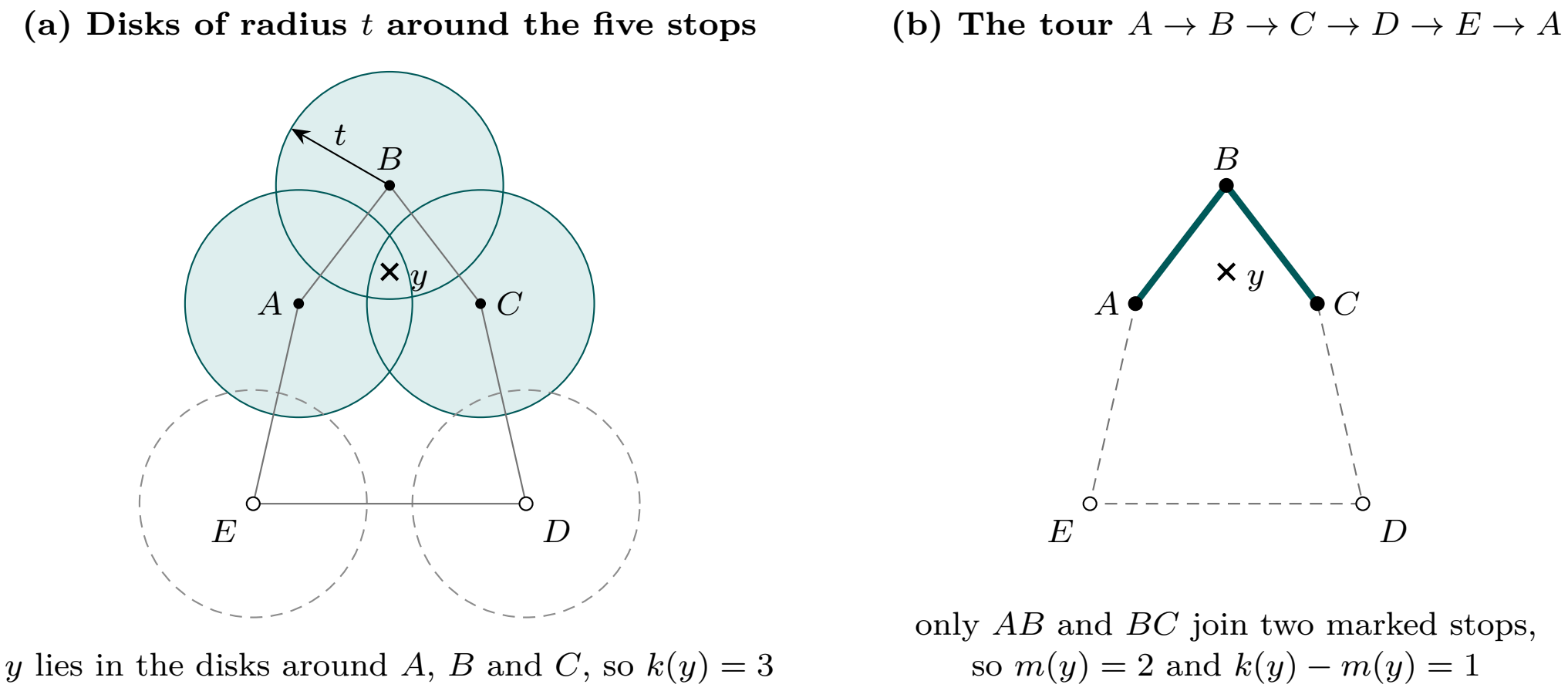

\tourConstraint*
\begin{proof}
Number the facilities in tour order, put $B_i=B_t(x_i)$ and $B_{N+1}=B_1$. For a point $y$ of the
plane, let $k(y)$ be the number of disks $B_i$ that contain $y$, and $m(y)$ the number of lenses
$B_i\cap B_{i+1}$ that contain $y$. We claim that
\begin{equation}\label{eq:count}
k(y)-m(y)\ge 1\quad\text{whenever } 1\le k(y)\le N-1 .
\end{equation}
Mark the stops whose disks contain $y$ (Figure~\ref{fig:count}) and orient the tour in visiting
order. There are $k(y)$ tour edges ending at marked stops. Of these, $m(y)$ also start at marked stops.
When $0<k(y)<N$, start at an unmarked stop and follow the tour to the first marked stop. Its
incoming tour edge starts at an unmarked stop, so it is not counted in $m(y)$. Thus
$m(y)\le k(y)-1$, proving \eqref{eq:count}.

If $k(y)=N$, every stop is marked, all $N$ tour edges are counted, and $k(y)-m(y)=0$. Such points lie in
$\bigcap_i B_i$. If $k(y)=0$ then $m(y)=0$. In all cases
\[
\mathbf 1_{\bigcup_i B_i}(y)\;\le\;k(y)-m(y)+\mathbf 1_{\bigcap_i B_i}(y).
\]
Integrate over the plane. The integral of $k$ is $N\pi t^2$. The integral of $m$ is
$\sum_i\abs{B_i\cap B_{i+1}}=\sum_i o_t(\ell_i)$, where $N\ge 3$ ensures that the $N$ lenses belong to
$N$ different pairs. The common intersection $\bigcap_i B_i$ lies in one disk. Hence
\begin{equation}\label{eq:union}
\Bigl|\bigcup_{i}B_i\Bigr|\;\le\;N\pi t^2-\sum_{i}o_t(\ell_i)+\pi t^2 .
\end{equation}
The part of $Q_N$ covered by the disks is no larger than their union, so \eqref{eq:U} and
\eqref{eq:union} give \eqref{eq:tour}.
\end{proof}

\subsection{Service regions and side triangles}\label{sec:triangles}

Continue to assume that all facilities lie in the interior of $Q_N$. Lemma~\ref{lem:boundary} will remove this assumption.

Each customer uses the nearest facility, so facility $x_i$ serves the region
\begin{equation}\label{eq:cell}
C_i=\bigl\{x\in Q_N:\ \norm{x-x_i}\le\norm{x-x_j}\ \text{for all } j\bigr\},
\end{equation}
its Voronoi cell, clipped to the square. For each $j\ne i$, the condition $\norm{x-x_i}\le\norm{x-x_j}$
describes the half-plane bounded by the perpendicular bisector of $x_i$ and $x_j$. Hence $C_i$ is a
convex polygon, each of whose sides lies on such a bisector or on the boundary of $Q_N$. The cells cover
$Q_N$ and meet only along their boundaries, and $x_i$ is an interior point of $C_i$ since the
facilities are interior and distinct.

\begin{figure}[t]
\centering
\begin{adjustbox}{max width=\linewidth}
\begin{tikzpicture}

\begin{scope}[scale=2.5]
\coordinate (O) at (0,0); \coordinate (A) at (1,-0.2126); \coordinate (B) at (1,1.2799);
\fill[capc] (1,-0.2126) -- (1,0.8307) arc[start angle=39.715,end angle=-12,radius=1.3] -- cycle;
\fill[soft] (O) -- (A) -- (B) -- cycle;
\draw[black!45,densely dashed,thin] (1,-0.62) -- (1,1.45);
\draw[acc,thin] (-48:1.3) arc[start angle=-48,end angle=60,radius=1.3];
\draw[black!50,thin] (A) -- ($(O)!1.42!(A)$);
\draw[black!50,thin] (B) -- ($(O)!1.10!(B)$);
\draw[thick] (O) -- (A) -- (B) -- cycle;
\draw[thin,densely dotted] (O) -- (1,0);
\node[lab,above=-1.5pt] at (0.8,0) {$h$};
\draw[thin] (0.93,0) -- (0.93,0.07) -- (1,0.07);
\draw[->,thin] (O) -- (1,0.5774);
\node[lab,above=-1.5pt,rotate=30] at (0.58,0.335) {$s=h\sec\theta$};
\draw[thin] (0.3,0) arc[start angle=0,end angle=30,radius=0.3];
\node[lab] at (14:0.4) {$\theta$};
\draw[->,thin] (O) -- (-34:1.3);
\node[lab,below left=-2pt] at (-34:0.7) {$t$};
\node[dot] at (O) {}; \node[lab,left] at (O) {$x_i$};
\node[lab,anchor=south west] at ($(O)!1.10!(B)$) {$\theta_+$};
\node[lab,anchor=west] at ($(O)!1.42!(A)$) {$\theta_-$};
\node[lab,anchor=west] at (1.42,0.42) {overhang $e_t(T)$};
\draw[thin] (1.41,0.42) -- (1.15,0.36);
\node[lab] at (0.72,0.72) {$T$};
\end{scope}
\node[lab,font=\scriptsize\bfseries] at (1.7,4.15) {(a) A side triangle};

\begin{scope}[shift={(8.75,1.25)},scale=3.0]
\draw[black!45,densely dashed,thin] (0,0) circle (0.64668);
\begin{scope}\clip (0.45,-0.4) rectangle (0.6,0.4); \fill[capc] (0,0) circle (0.53);\end{scope}
\fill[soft] (0,0) -- (0.45,-0.46443) -- (0.45,0.46443) -- cycle;
\draw[acc,thin] (-52:0.53) arc[start angle=-52,end angle=52,radius=0.53];
\draw[thick] (0,0) -- (0.45,-0.46443) -- (0.45,0.46443) -- cycle;
\draw[thin,densely dotted] (0,0) -- (0.45,0);
\node[lab,below=-1pt] at (0.3,0) {$h$};
\draw[thin] (0.16,0) arc[start angle=0,end angle=45.9,radius=0.16];
\node[lab] at (21:0.23) {$\beta$};
\node[lab,above left=-2pt] at (0.2,0.22) {$R$};
\node[dot] at (0,0) {};
\node[lab,anchor=west] at (0.7,0.0) {cap $c_h(t)$};
\draw[thin] (0.69,0) -- (0.5,0.0);
\node[lab,anchor=north] at (0,-0.67) {dashed: circle of radius $R$};
\end{scope}
\node[lab,font=\scriptsize\bfseries] at (8.75,4.15) {(b) The symmetric family};

\begin{scope}[shift={(13.75,1.35)},scale=1.15]
\fill[capc] (0,0) circle (1.12);
\fill[white] (1.25,-0.15) -- (0.95,0.95) -- (-0.2,1.3) -- (-1.15,0.6) -- (-1.0,-0.7) -- (0.3,-1.1) -- cycle;
\fill[soft] (1.25,-0.15) -- (0.95,0.95) -- (-0.2,1.3) -- (-1.15,0.6) -- (-1.0,-0.7) -- (0.3,-1.1) -- cycle;
\foreach \p in {(1.25,-0.15),(0.95,0.95),(-0.2,1.3),(-1.15,0.6),(-1.0,-0.7),(0.3,-1.1)} {\draw[black!45,thin] (0,0) -- \p;}
\draw[acc,thin] (0,0) circle (1.12);
\draw[thick] (1.25,-0.15) -- (0.95,0.95) -- (-0.2,1.3) -- (-1.15,0.6) -- (-1.0,-0.7) -- (0.3,-1.1) -- cycle;
\node[dot] at (0,0) {};
\node[lab,above left=-2pt] at (0,0) {$x_i$};
\node[lab] at (0.42,0.62) {$C_i$};
\node[lab,anchor=north] at (0.05,-1.33) {shaded: $B_t(x_i)\setminus C_i$};
\end{scope}
\node[lab,font=\scriptsize\bfseries] at (13.75,4.15) {(c) One cell};
\end{tikzpicture}
\end{adjustbox}
\caption{(a) A side triangle in the coordinates \eqref{eq:tri}. Its side lies on the dashed line
$\xi=h$. The overhang is the part of the facility's disk that lies between the two rays bounding $T$ but
beyond its side. It need not be a whole cap. (b) A symmetric
triangle with corners at distance $R$ and half-angle $\beta=\arccos(h/R)$. For $t\le b<R$ its overhang
is the whole cap. (c) The side triangles of a cell. Their overhangs together form the part of the disk
outside the cell.}
\label{fig:triangle}
\end{figure}
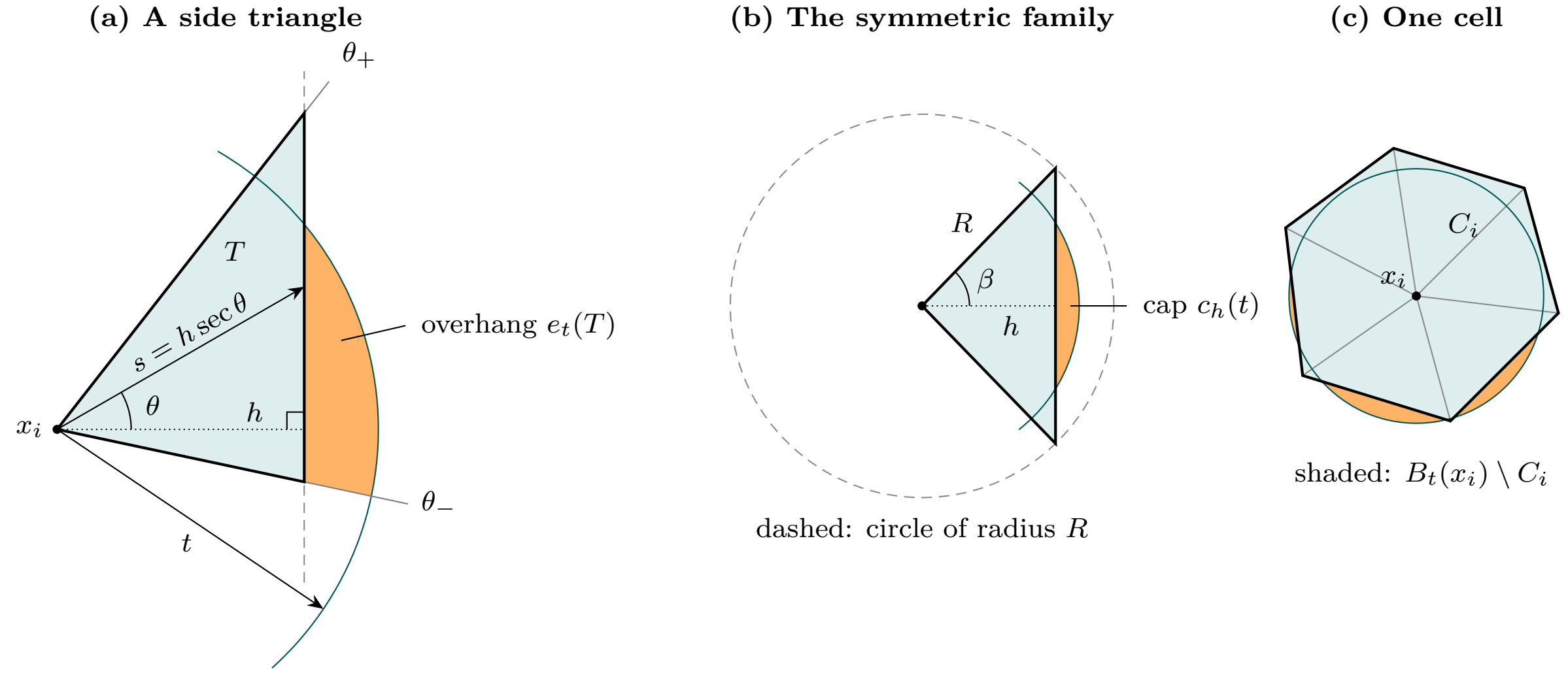

Join $x_i$ to the corners of $C_i$. This cuts the cell into triangles with common apex $x_i$, one for
each side (Figure~\ref{fig:triangle}c). We call them \emph{side triangles} and write $\cT$ for the set
of all side triangles of all cells. The service cost is additive over them:
\[
M_N(X)=\sum_{T\in\cT}M(T),\qquad M(T)=\int_T\norm{x-x_i}\dd x\quad(x_i\text{ the apex of }T).
\]
Besides service cost, we use the area of each triangle, its angle at the facility, and the disk area beyond its opposite side. Summing these quantities over the side triangles will give the four constraints on a layout and tour.

\paragraph{Coordinates for one triangle.}
Translate the apex of $T$ to the origin, and rotate the plane so that the side of $T$ lies on the
vertical line $\xi=h$ in local Cartesian coordinates $(\xi,\zeta)$. Here $h>0$ is the perpendicular distance from the facility to the line through the side.
These rigid motions preserve lengths and areas. In polar coordinates $(s,\theta)$ the line is
$s\cos\theta=h$, so
\begin{equation}\label{eq:tri}
T=\bigl\{(s,\theta):\ \theta_-\le\theta\le\theta_+,\ 0\le s\le h\sec\theta\bigr\},\qquad
-\tfrac{\pi}{2}<\theta_-<\theta_+<\tfrac{\pi}{2},
\end{equation}
where $\theta_\pm$ are the directions of the two other corners (Figure~\ref{fig:triangle}a). The three
numbers $h,\theta_-,\theta_+$ describe the triangle completely. In general, the triangle is not symmetric: the
two corners can be at different distances from the facility, and the foot of the perpendicular, at
$\theta=0$, can lie outside the side.

\paragraph{Area, angle and overhang.}
A ray from the facility in direction $\theta$ crosses $T$ in a segment of length $s=h\sec\theta$.
Integration in polar coordinates gives
\begin{equation}\label{eq:four}
M(T)=\int_{\theta_-}^{\theta_+}\frac{s^3}{3}\dd\theta,\qquad
\abs{T}=\int_{\theta_-}^{\theta_+}\frac{s^2}{2}\dd\theta,\qquad
\omega(T)=\int_{\theta_-}^{\theta_+}\dd\theta=\theta_+-\theta_-.
\end{equation}
Here $s=h\sec\theta$ and $\omega(T)$ is the angle at the facility. Where $t>s$, the radius-$t$ disk extends beyond the opposite side by area $(t^2-s^2)/2$ per unit angle. We call
\begin{equation}\label{eq:overhang}
e_t(T)=\int_{\theta_-}^{\theta_+}\frac{(t^2-s^2)_+}{2}\dd\theta
\end{equation}
the \emph{overhang} of the disk beyond $T$ (shaded in Figure~\ref{fig:triangle}a). Summed over all
side triangles of the layout, these quantities and the number of triangles satisfy
\begin{equation}\label{eq:sums}
\sum_{T\in\cT}\abs{T}=N,\qquad
\sum_{T\in\cT}\omega(T)=2\pi N,\qquad
\sum_{T\in\cT}e_t(T)=N\pi t^2-N+U_X(t),\qquad
\#\cT\le 6N .
\end{equation}
The first identity follows by adding the areas of triangles that tile the square. The second follows by adding the angles around each facility, which sum to $2\pi$. For the third, the overhangs around facility $x_i$ make up
$B_t(x_i)\setminus C_i$, the part of its disk outside its cell (Figure~\ref{fig:triangle}c). The parts inside, $C_i\cap B_t(x_i)$, overlap only on cell boundaries, which have area zero. Their union is the covered part of the square, since a customer is within $t$ of some facility exactly when it is within $t$ of its own. Their
total area is $N-U_X(t)$, and subtracting it from the $N$ disk areas gives the third sum. The bound on
the number of triangles is the following lemma, a consequence of Euler's formula (\appSide, see also
\citeGruber).

\Needspace{4\baselineskip}
\begin{restatable}[Side count]{lemma}{sideCount}\label{lem:sides}
For distinct interior facilities, the clipped Voronoi cells have at most $6N$ sides in total, where a
side shared by two cells is counted twice.
\end{restatable}

Total area and total angle in \eqref{eq:sums} are fixed. Total overhang depends on $U_X(t)$, and the triangle count is bounded by Lemma~\ref{lem:sides}. Combining the overhang identity with Lemma~\ref{lem:tour} gives
\begin{equation}\label{eq:floor}
\sum_{T\in\cT}e_t(T)\;\ge\;\sum_{i=1}^{N}o_t(\ell_i)-\pi t^2\qquad\text{for every } t\ge 0 .
\end{equation}

\subsection{The relaxed problem}\label{sec:relaxed}

Sections~\ref{sec:uncovered} and~\ref{sec:triangles} associate two lists with every layout and tour: the side triangles and the tour-edge lengths. The relaxed problem keeps their total area, total angle, triangle-count bound, and overhang constraint, as given in \eqref{eq:sums} and \eqref{eq:floor}.

A \emph{collection} is a finite list of triangles, with repetitions allowed. We include the empty list, although it cannot satisfy the area constraint below. Each triangle is given by numbers
$(h,\theta_-,\theta_+)$ as in \eqref{eq:tri}, and $M(T)$, $\abs T$, $\omega(T)$, $e_t(T)$ are defined by
\eqref{eq:four} and \eqref{eq:overhang}. The \emph{relaxed problem} has as unknowns a collection $\cC$
and $N$ numbers $\ell_1,\dots,\ell_N\ge 0$:
\begin{equation}\label{eq:relaxed}
\begin{array}{@{}lll@{}}
\underset{\cC,\,\boldsymbol\ell\in[0,\infty)^N}{\text{minimize}} & \displaystyle\sum_{T\in\cC}M(T)+r\sum_{i=1}^{N}\ell_i & \\[18pt]
\text{subject to} & \displaystyle\sum_{T\in\cC}\abs T=N & \text{[price $S$]}\\[12pt]
& \displaystyle\sum_{T\in\cC}\omega(T)=2\pi N & \text{[price $\lambda$]}\\[12pt]
& \#\cC\le 6N & \text{[price $-v\ge 0$]}\\[6pt]
& \displaystyle\sum_{T\in\cC}e_t(T)\ge\sum_{i=1}^{N}o_t(\ell_i)-\pi t^2\quad\text{for all } a\le t\le b
& \text{[price $\eta(t)\ge 0$]}
\end{array}
\end{equation}
The brackets indicate the multipliers used in Section~\ref{sec:prices}.

The optimal value is the infimum over collections and length lists satisfying the four constraints. Every layout $X$ of $N$ distinct interior facilities supplies a feasible choice: use its side triangles and the edge lengths of a closed tour through its facilities. Equations~\eqref{eq:sums} and~\eqref{eq:floor} verify all four constraints, and their relaxed objective equals $M_N(X)+rL$. Therefore
\begin{equation}\label{eq:relaxbound}
M_N(X)+rL\;\ge\;\text{the optimal value of \eqref{eq:relaxed}} .
\end{equation}
The relaxed problem discards most of the geometry. Its triangles need not fit together into cells that tile the square, its tour-edge lengths need not come from a tour through the facilities, and only the overlaps of consecutive disks enter. The overhang requirement is imposed for $a\le t\le b$. Section~\ref{sec:weight} explains this choice.

\Needspace{6\baselineskip}
\subsection{Lagrange multipliers and the lower bound}\label{sec:prices}

Attach multipliers $S,\lambda\in\R$ to the area and angle equalities in \eqref{eq:relaxed}, respectively, and $-v\ge0$ to its side-count inequality. For the inequalities indexed by radius, use a nonnegative integrable function $\eta$ on $[a,b]$. We also call these multipliers \emph{prices}. For a collection $\cC$ and a list $\boldsymbol\ell=(\ell_1,\ldots,\ell_N)$, the Lagrangian of \eqref{eq:relaxed} is
\[
\begin{aligned}
\mathcal L(\cC,\boldsymbol\ell\mid S,\lambda,v,\eta)
={}&\sum_{T\in\cC}M(T)+r\sum_{i=1}^N\ell_i
+S\left(N-\sum_{T\in\cC}\abs T\right)\\
&+\lambda\left(\sum_{T\in\cC}\omega(T)-2\pi N\right)
-v(\#\cC-6N)\\
&+\int_a^b\eta(t)\left[\sum_{i=1}^No_t(\ell_i)-\pi t^2
-\sum_{T\in\cC}e_t(T)\right]\dd t .
\end{aligned}
\]
For a feasible choice in \eqref{eq:relaxed}, the equality terms vanish, and the inequality terms are nonpositive. Thus $\mathcal L$ is no larger than the objective of \eqref{eq:relaxed}. All its integrals are finite since $0\le e_t(T),o_t(\ell)\le\pi t^2\le\pi b^2$ and $\eta$ is integrable.

Collect the terms associated with one triangle and one edge length:
\begin{equation}\label{eq:Gamma}
\Gamma(T)=M(T)-S\abs{T}+\lambda\,\omega(T)-\int_a^b\eta(t)\,e_t(T)\dd t ,
\end{equation}
\begin{equation}\label{eq:w}
w(\ell)=r\,\ell+\int_a^b\eta(t)\,o_t(\ell)\dd t .
\end{equation}
The same radius multiplier subtracts the overhang term from $\Gamma$ and adds the overlap term to $w$. Set $K=\pi\int_a^b\eta(t)t^2\dd t$. Grouping terms gives the identity in the second line below for every collection and nonnegative length list. The first inequality holds when they satisfy \eqref{eq:relaxed}:
\begin{equation}\label{eq:lagrange}
\begin{aligned}
\sum_{T\in\cC}M(T)+r\sum_{i=1}^N\ell_i
&\ge \mathcal L(\cC,\boldsymbol\ell \mid S,\lambda,v,\eta)\\
&=\sum_{T\in\cC}[\Gamma(T)-v]+\sum_{i=1}^Nw(\ell_i)
+N(S-2\pi\lambda+6v)-K .
\end{aligned}
\end{equation}
To obtain a lower bound, minimize $\mathcal L$ after removing the four aggregate constraints:
\[
g(S,\lambda,v,\eta)=\inf_{\cC,\,\boldsymbol\ell\in[0,\infty)^N}
\mathcal L(\cC,\boldsymbol\ell\mid S,\lambda,v,\eta).
\]
The collection may now be any finite list of triangles of the form \eqref{eq:tri}, including the empty list. If some triangle satisfies $\Gamma(T)-v<0$, repeating it makes $g=-\infty$. To obtain a finite lower bound we therefore require
\begin{enumerate}[label=(\Alph*),topsep=3pt,itemsep=1pt]
\item $\Gamma(T)\ge v$ for every triangle $T$, whatever the admissible values of $h$, $\theta_-$ and $\theta_+$.
\end{enumerate}
Under (A), each term in $\sum_{T\in\cC}[\Gamma(T)-v]$ is nonnegative and the empty collection makes this sum zero. Each of the $N$ edge terms can be minimized independently. Since $w(\ell)\ge0$ and $w(0)<\infty$, we obtain
\[
g(S,\lambda,v,\eta)=N(S-2\pi\lambda+6v)
+N\inf_{\ell\ge0}w(\ell)-K
\qquad\text{when (A) holds}.
\]
Every such value of $g$ is a lower bound on the relaxed optimum and hence on $M_N(X)+rL$. The Lagrange dual problem is to maximize $g$ over the multipliers defined above.

For the proof below, we use the bound obtained directly by dropping the nonnegative triangle terms from \eqref{eq:lagrange}:
\begin{equation}\label{eq:dualbound}
\sum_{T\in\cC}M(T)+r\sum_{i=1}^{N}\ell_i
\;\ge\;N(S-2\pi\lambda+6v)+\sum_{i=1}^{N}w(\ell_i)-K .
\end{equation}
Condition (A) is equivalent to the single-triangle bound
\[
M(T)-\int_a^b\eta(t)e_t(T)\dd t
\ge S\abs T-\lambda\omega(T)+v.
\]
Summing the right-hand side over the triangles gives $SN-2\pi N\lambda+v\#\cC$: their areas sum to $N$ and their angles sum to $2\pi N$. Since $\#\cC\le6N$ and $v\le0$, we have $v\#\cC\ge6Nv$. When $\eta=0$, the single-triangle bound reduces to $M(T)\ge S|T|-\lambda\omega(T)+v$, the form used in the argument for regular hexagons (\citeGruber; \citeXO, Lemma~1). The function $\eta$ incorporates the tour through the overhang constraint. The constant $K$ comes from the allowance $\pi t^2$ in Lemma~\ref{lem:tour}.

\subsection{Matching the reference cost}\label{sec:right}

To match the reference cost, we seek multipliers for which the right-hand side of \eqref{eq:dualbound} is $N\{\phi(d)+rd\}-K$.

In the relaxation, choose $N$ copies of the six side triangles of $H_d$ and set all $N$ length variables equal to $d$. This choice has total area $N$, total angle $2\pi N$, and $6N$ triangles. By \eqref{eq:ud}, its total overhang is $No_t(d)$ for $a\le t\le b$, which exceeds the required lower bound $No_t(d)-\pi t^2$ by $\pi t^2$. Thus it is feasible for \eqref{eq:relaxed}, with objective $N\{\phi(d)+rd\}$.

To make the triangle terms in \eqref{eq:lagrange} vanish for this choice, we require
\begin{enumerate}[label=(\Alph*),topsep=3pt,itemsep=1pt,start=2]
\item $\Gamma(T)=v$ for each of the six side triangles of $H_d$.
\end{enumerate}
Together with (A), this means that the reference triangles minimize $\Gamma$. Their total area is $1$, their total angle is $2\pi$, and their total overhang is $o_t(d)$. Summing (B) over the six triangles gives
\begin{equation}\label{eq:hexsum}
\phi(d)=S-2\pi\lambda+6v+\int_a^b\eta(t)\,o_t(d)\dd t .
\end{equation}

Section~\ref{sec:lb2} will choose the multipliers as functions of $d$, with $\eta$ given explicitly in \eqref{eq:eta}. For that choice, define
\begin{equation}\label{eq:rK}
\rs(d)=\int_a^b\eta(t)\sqrt{4t^2-d^2}\dd t,\qquad
K(d)=\pi\int_a^b\eta(t)\,t^2\dd t .
\end{equation}
Here $K(d)$ is the constant $K$ above for the chosen multiplier. Both numbers are finite and nonnegative since $\sqrt{4t^2-d^2}\le2b$ and $t^2\le b^2$ on $[a,b]$. When $d=d_0$, the interval has zero length, so $\rs(d_0)=K(d_0)=0$.

\begin{figure}[t]
\centering
\begin{adjustbox}{max width=\linewidth}
\begin{tikzpicture}
\draw[->] (0,0) -- (6.9,0) node[below,lab]{$\ell$};
\draw[->] (0,0) -- (0,3.6) node[left,lab]{$o_t(\ell)$};
\draw[dashed,thin] (1.500,1.583) -- (5.180,0);
\draw[acc,very thick] plot coordinates {(0.000,3.231) (0.125,3.138) (0.250,3.044) (0.375,2.951) (0.500,2.858) (0.625,2.765) (0.750,2.672) (0.875,2.579) (1.000,2.487) (1.125,2.396) (1.250,2.304) (1.375,2.213) (1.500,2.123) (1.625,2.034) (1.750,1.945) (1.875,1.856) (2.000,1.769) (2.125,1.682) (2.250,1.596) (2.375,1.511) (2.500,1.428) (2.625,1.345) (2.750,1.263) (2.875,1.183) (3.000,1.104) (3.125,1.026) (3.250,0.950) (3.375,0.875) (3.500,0.803) (3.625,0.731) (3.750,0.662) (3.875,0.595) (4.000,0.529) (4.125,0.466) (4.250,0.406) (4.375,0.348) (4.500,0.292) (4.625,0.240) (4.750,0.191) (4.875,0.146) (5.000,0.105) (5.125,0.068) (5.250,0.037) (5.375,0.013) (5.500,0.000) (6.500,0.000)};
\node[dot] at (4.500,0.292) {};
\draw[thin] (4.500,0) -- ++(0,-0.08) node[below,lab]{$d$};
\draw[thin] (5.500,0) -- ++(0,-0.08) node[below,lab]{$2t$};
\node[lab,acc,anchor=west] at (1.15,2.75) {overlap of two disks, centers $\ell$ apart};
\node[lab,anchor=north] at (2.1,0.62) {tangent at $d$};
\draw[thin] (2.25,0.62) -- (2.75,1.03);
\node[lab,anchor=west,text width=5.2cm] at (7.4,2.25)
{Every tour edge is compared with the same tangent, so the linear terms add up to $Nd-L$.};
\node[lab,anchor=west,text width=5.2cm] at (7.4,1.05)
{Weighting the slope $-\sqrt{4t^2-d^2}$ by $\eta(t)$ and integrating over $[a,b]$ gives $-\rs(d)$.};
\end{tikzpicture}
\end{adjustbox}
\caption{The tangent inequality \eqref{eq:tangent}, drawn for $d=0.9$ and $t=0.55$. Convexity keeps
the curve above the tangent for every tour-edge length, shorter or longer than $d$.}
\label{fig:tangent}
\end{figure}
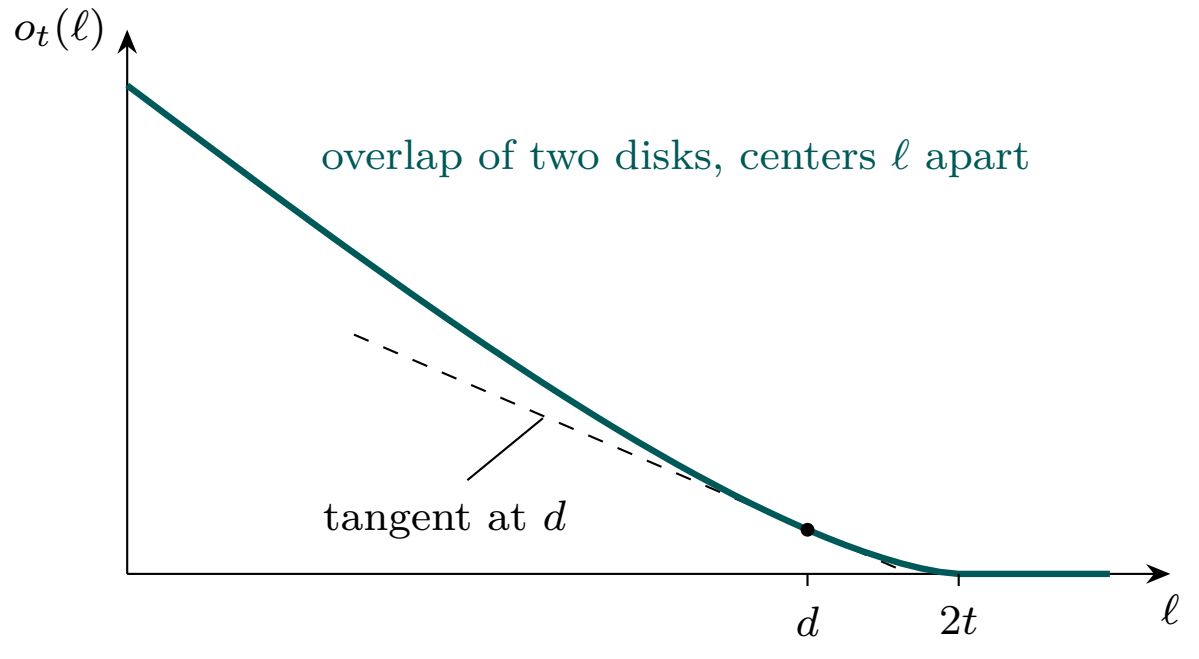

We also need $d$ to minimize the edge contribution $w$. Convexity of $o_t$ places its graph above its tangent at $d$ (Figure~\ref{fig:tangent}):
\begin{equation}\label{eq:tangent}
o_t(\ell)\;\ge\;o_t(d)-\sqrt{4t^2-d^2}\,(\ell-d),\qquad \ell\ge0,\quad a\le t\le b .
\end{equation}
Multiplying by $\eta(t)\ge0$, integrating and adding $r\ell$ yields
\begin{equation}\label{eq:wtangent}
w(\ell)\;\ge\;w(d)+\bigl(r-\rs(d)\bigr)(\ell-d),\qquad \ell\ge0 .
\end{equation}
Thus, if $r=\rs(d)$ then $w(\ell)\ge w(d)$ for every $\ell\ge0$. For a given $r$, Lemma~\ref{lem:spacing} will select $d$ with this property.

Assume (A), (B), and $r=\rs(d)$. The triangle terms vanish at the reference triangles, and the reference length minimizes every edge term. Using \eqref{eq:hexsum} to evaluate the dual function gives
\[
g(S,\lambda,v,\eta)=N\{\phi(d)+rd\}-K(d).
\]
The reference collection is feasible for the relaxation, so its optimal value is bounded on both sides:
\[
N\{\phi(d)+rd\}-K(d)
\ \le\ \text{optimal value of \eqref{eq:relaxed}}
\ \le\ N\{\phi(d)+rd\}.
\]
The gap between these two bounds is $K(d)$, the integral of the reference overhang excess $\pi t^2$ against $\eta(t)$. At $d=d_0$ and $r=0$, $a=b$ and $K(d_0)=0$, so the bounds coincide.

The remaining step is to construct multipliers satisfying (A) and (B):

\Needspace{4\baselineskip}
\begin{restatable}[Prices]{lemma}{multiplierPrices}\label{lem:prices}
For every $d\in(0,d_0]$ there are real numbers $S,\lambda$, a number $v<0$, and a nonnegative integrable function
$\eta$ on $[a,b]$, depending only on $d$, for which conditions (A) and (B) hold.
\end{restatable}

\section{Proof of Theorem 2: construction of the multipliers}\label{sec:lb2}

We prove conditions (A) and (B), then complete Theorem~\ref{thm:lower} in Section~\ref{sec:completion}. The reference triangles have side distance $a$ or $b$ and both outer vertices at distance $R$ from the facility. We first determine $S$ and $\lambda$ by requiring the first derivative of $\Gamma$ to vanish under each of two variations of a reference triangle. We then choose $\eta$ so that $\Gamma(T_h)$ is constant for the symmetric triangles with $a\le h\le b$. Section~\ref{sec:alltriangles} proves the uniform bound (A) for every admissible triangle.

\subsection{\texorpdfstring{Prices of area and angle}{Prices of area and angle}}\label{sec:constants}

For a triangle whose opposite side lies on $\xi=h$, the ray in direction $\theta$ ends at distance
$s=h\sec\theta$. Its service integral contributes $s^3/3$ per unit angle, its area contributes
$s^2/2$, and its disk overhang contributes $(t^2-s^2)_+/2$. Substituting these expressions into
\eqref{eq:Gamma} gives
\begin{equation}\label{eq:Psi}
\Gamma(T)=\int_{\theta_-}^{\theta_+}\Psi(h\sec\theta)\dd\theta,\qquad
\Psi(s)=\frac{s^3}{3}-\frac{S s^2}{2}+\lambda-\frac12\int_a^b\eta(t)\,(t^2-s^2)_+\dd t .
\end{equation}
$\Psi(s)$ is the net contribution per unit angle of a ray ending at distance $s$. Integrating it over the directions of a triangle gives $\Gamma(T)$.

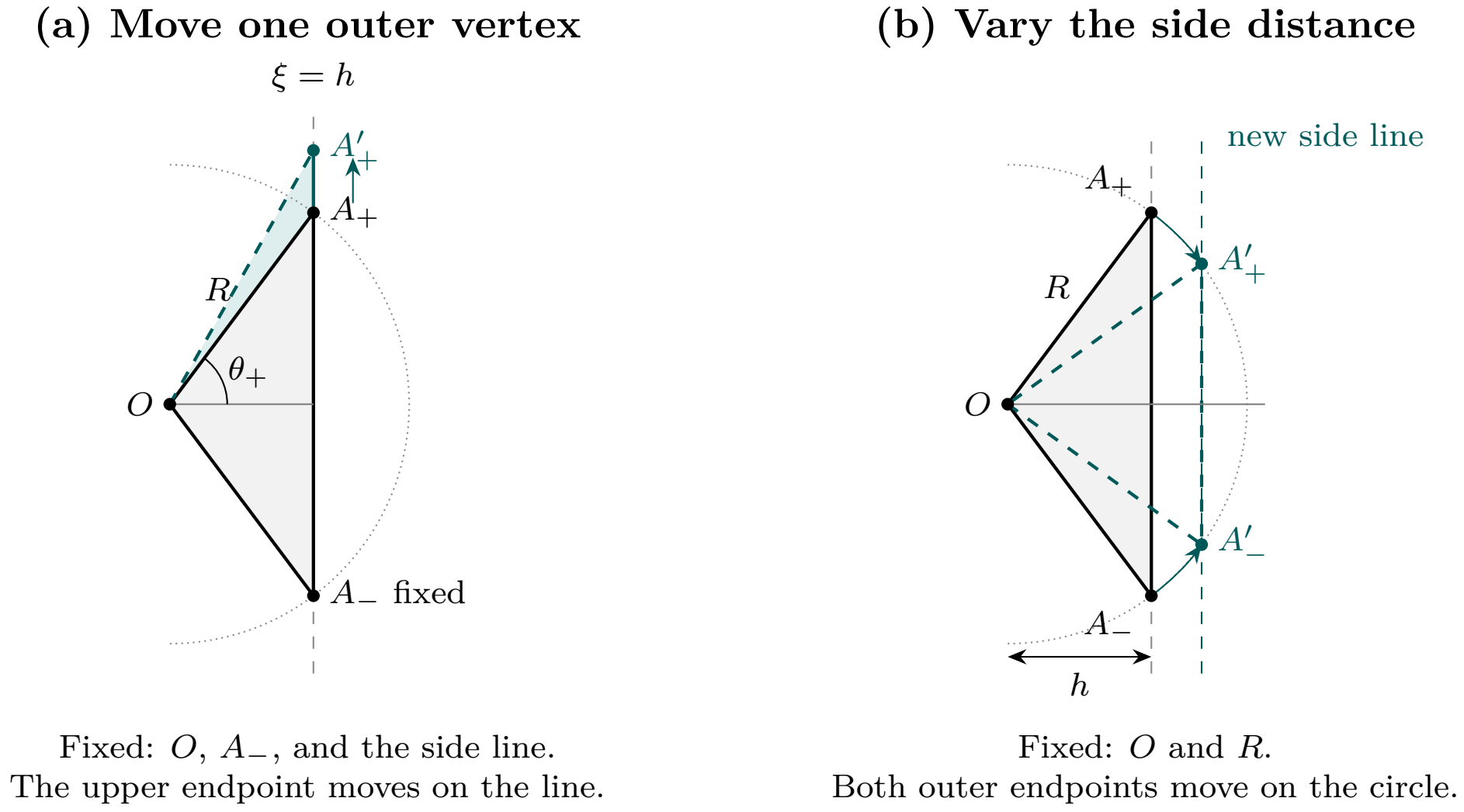
\begin{figure}[htbp]
\centering
\begin{adjustbox}{max width=\linewidth}
\begin{tikzpicture}[font=\scriptsize,scale=1.02]
\begin{scope}
\node[font=\small\bfseries] at (1.15,3.15) {(a) Move one outer vertex};
\fill[black!5] (0,0)--(1.2,-1.6)--(1.2,1.6)--cycle;
\fill[soft] (0,0)--(1.2,1.6)--(1.2,2.12)--cycle;
\draw[black!45,densely dotted] (0,-2) arc[start angle=-90,end angle=90,radius=2];
\draw[black!45,dashed] (1.2,-2.25)--(1.2,2.48);
\node[anchor=south] at (1.2,2.48) {$\xi=h$};
\draw[thick] (0,0)--(1.2,-1.6)--(1.2,1.6)--cycle;
\draw[acc,thick,dashed] (0,0)--(1.2,2.12);
\draw[acc,thick] (1.2,1.6)--(1.2,2.12);
\draw[acc,->] (1.53,1.68)--(1.53,2.06);
\draw[black!50] (0,0)--(1.2,0);
\draw (0.48,0) arc[start angle=0,end angle=53.1301,radius=.48];
\node at (.66,.27) {$\theta_+$};
\fill (0,0) circle(1.5pt) node[left]{$O$};
\fill (1.2,-1.6) circle(1.5pt) node[right]{$A_-$ fixed};
\fill (1.2,1.6) circle(1.5pt) node[right]{$A_+$};
\fill[acc] (1.2,2.12) circle(1.5pt) node[right,acc]{$A_+'$};
\node[anchor=east] at (.64,.96) {$R$};
\node[align=center] at (1.15,-3.05) {Fixed: $O$, $A_-$, and the side line.\\The upper endpoint moves on the line.};
\end{scope}
\begin{scope}[shift={(7,0)}]
\node[font=\small\bfseries] at (1.15,3.15) {(b) Vary the side distance};
\fill[black!5] (0,0)--(1.2,-1.6)--(1.2,1.6)--cycle;
\draw[black!45,densely dotted] (0,-2) arc[start angle=-90,end angle=90,radius=2];
\draw[black!45,dashed] (1.2,-2.25)--(1.2,2.25);
\draw[acc,dashed] (1.62,-2.25)--(1.62,2.25);
\draw[thick] (0,0)--(1.2,-1.6)--(1.2,1.6)--cycle;
\draw[acc,thick,dashed] (0,0)--(1.62,-1.17286)--(1.62,1.17286)--cycle;
\draw[acc,->] (1.2,1.6) arc[start angle=53.1301,end angle=35.9041,radius=2];
\draw[acc,->] (1.2,-1.6) arc[start angle=-53.1301,end angle=-35.9041,radius=2];
\draw[black!50] (0,0)--(2.15,0);
\fill (0,0) circle(1.5pt) node[left]{$O$};
\fill (1.2,1.6) circle(1.5pt) node[above left]{$A_+$};
\fill (1.2,-1.6) circle(1.5pt) node[below left]{$A_-$};
\fill[acc] (1.62,1.17286) circle(1.5pt) node[right,acc]{$A_+'$};
\fill[acc] (1.62,-1.17286) circle(1.5pt) node[right,acc]{$A_-'$};
\node[anchor=east] at (.65,.98) {$R$};
\draw[<->] (0,-2.11)--(1.2,-2.11) node[midway,below]{$h$};
\node[acc,anchor=west] at (1.7,2.25) {new side line};
\node[align=center] at (1.15,-3.05) {Fixed: $O$ and $R$.\\Both outer endpoints move on the circle.};
\end{scope}
\end{tikzpicture}
\end{adjustbox}
\caption{The two variations used to choose $\lambda$ and $S$. The solid triangle is the reference position, and the dashed outline shows a perturbation. (a) Only the upper outer vertex moves along the fixed side line. This gives the necessary condition $\Psi(R)=0$. (b) The line moves parallel to itself, and both endpoints stay at distance $R$ from $O$. At $h=b$, this gives $V'(b)=0$.}
\label{fig:variations}
\end{figure}

\paragraph{Moving one outer vertex determines $\lambda$.}
Take one reference triangle and put its facility at $O=(0,0)$, with opposite side on the line
$\xi=h$, where $h=a$ or $h=b$. Its outer vertices are
$A_-=(h,-\sqrt{R^2-h^2})$ and $A_+=(h,\sqrt{R^2-h^2})$.
Keep $O$, $A_-$ and the line $\xi=h$ fixed, and slide $A_+$ a small distance up or down that line
(Figure~\ref{fig:variations}a). The upper endpoint need not remain on the circle during this variation.
Only the upper angular endpoint $\theta_+$ changes.

By \eqref{eq:Psi}, the derivative of $\Gamma$ with respect to $\theta_+$ is
$\Psi(h\sec\theta_+)$. At the reference triangle, $h\sec\theta_+=R$, so this derivative is
$\Psi(R)$. Conditions (A) and (B) require the reference triangle to minimize $\Gamma$. Since the endpoint can move in either direction, this derivative must vanish. Thus $\Psi(R)=0$. Since $R>b$, the overhang term in \eqref{eq:Psi} vanishes at $R$, giving

\begin{equation}\label{eq:lambda}
\lambda=\frac{S R^2}{2}-\frac{R^3}{3}=\int_0^R(S-u)\,u\dd u .
\end{equation}
Every reference endpoint is at distance $R$, so this choice of $\lambda$ satisfies the same necessary condition for all six triangles.

\paragraph{Varying the symmetric triangle determines $S$.}
Keep the facility $O$ and the radius $R$ fixed, and move the opposite side parallel to itself. At each distance $0<h<R$, place its two endpoints at
\[
A_-(h)=(h,-\sqrt{R^2-h^2}),\qquad A_+(h)=(h,\sqrt{R^2-h^2}).
\]
As $h$ varies, both endpoints move along the circle of radius $R$ (Figure~\ref{fig:variations}b). Denote this symmetric triangle by $T_h$. The reference triangles are
copies of $T_a$ and $T_b$. The coefficients $S$ and $\lambda$ are held fixed when $h$ varies.

Define
\begin{equation}\label{eq:V}
V(h)=J_R(h)-S\,h\sqrt{R^2-h^2}+2\lambda\arccos(h/R),\qquad 0<h<R .
\end{equation}
By \lemHex, the three terms are respectively the service integral of $T_h$, its area multiplied by $-S$, and its angle multiplied by $\lambda$. Thus $V(h)=M(T_h)-S|T_h|+\lambda\omega(T_h)$ is the value before subtracting the weighted overhang.

For $t\le b<R$, the overhang of $T_h$ is the cap $c_h(t)$. The cap derivative \eqref{eq:capderiv} gives, for $a<h<b$,
\[
\frac{\mathrm d}{\mathrm dh}\int_a^b\eta(t)c_h(t)\dd t
=-2\int_h^b\eta(t)\sqrt{t^2-h^2}\dd t .
\]
Its absolute value is at most $2b\int_h^b\eta(t)\dd t$, which tends to zero as $h\uparrow b$ since $\eta$ is integrable. For $h\ge b$, all caps in this integral vanish. If $a=b$, the integral is identically zero. In either case, the derivative of $\Gamma(T_h)$ at $h=b$ is $V'(b)$. Since (A) and (B) require $T_b$ to minimize $\Gamma$, we must have $V'(b)=0$.

Differentiating the three terms in \eqref{eq:V} and substituting \eqref{eq:lambda} gives
\[
V'(h)=(R-2S)\sqrt{R^2-h^2}+h^2\arcosh(R/h)
\]
(the differentiation is also recorded in \appV). Setting $h=b$ gives
$(2S-R)\sqrt{R^2-b^2}=b^2\arcosh(R/b)$, and hence
\begin{equation}\label{eq:S}
S=\frac12\Bigl(R+\frac{b^2\arcosh(R/b)}{\sqrt{R^2-b^2}}\Bigr).
\end{equation}

\Needspace{4\baselineskip}
\begin{restatable}[Minimum of $V$]{lemma}{minimumOfV}\label{lem:V}
With $S$ and $\lambda$ as in \eqref{eq:S} and \eqref{eq:lambda}, we have $R/2<S<R$. The function $V$
decreases on $(0,b)$, increases on $(b,R)$, and extends continuously to $V(R)=0$. Its minimum is
$V(b)<0$.
\end{restatable}

The proof is in \appV. We set
\[
v=V(b).
\]
Each of the four far-side triangles of $H_d$ has net cost $v$: its weighted overhang is zero. Lemma~\ref{lem:V} gives $v<0$, as required in \eqref{eq:lagrange}. The price $\lambda$ is positive as well: the
estimate $\arcosh(R/b)=\int_b^R(u^2-b^2)^{-1/2}\dd u\ge\sqrt{R^2-b^2}/R$ gives $S\ge\tfrac12R(1+b^2/R^2)$. Since $b^2/R^2=(1+d^4/4)^{-1}\ge3/4$, we get $S\ge7R/8$ and $\lambda\ge5R^3/48$.

\subsection{The price of overhang}\label{sec:weight}

For $a<t\le b$, a disk of radius $t$ extends beyond a near side at distance $a$ but not beyond a far side at distance $b$. Subtracting the weighted overhang in \eqref{eq:Gamma} can therefore lower the value of $T_a$ without changing the value of $T_b$.

The same overhang also connects service to the tour. For an actual layout, \eqref{eq:sums} gives $\sum_T e_t(T)=N\pi t^2-N+U_X(t)$. Lemma~\ref{lem:tour} bounds this sum below by $\sum_i o_t(\ell_i)-\pi t^2$. We use radii in $[a,b]$, where the reference cell's overhang equals the overlap $o_t(d)$. For $t<a$, the overlap and its derivative at $d=2a$ vanish. These radii cannot penalize a small reduction of tour length. For $t>b$, the reference disk extends beyond the far sides as well, giving extra overhang in every cell and a gap proportional to $N$ in the reference comparison.

\paragraph{Net cost on the symmetric family.}
For $T_h$, the disk overhang is the full cap $c_h(t)$ when $t\le b<R$
(Figure~\ref{fig:triangle}b). Restricting \eqref{eq:Gamma} to this family therefore gives
\begin{equation}\label{eq:F}
F(h):=V(h)-\int_a^b\eta(t)\,c_h(t)\dd t .
\end{equation}
Thus $F(h)=\Gamma(T_h)$ is the net cost on the symmetric family, and $V(h)$ is its value before the weighted cap is subtracted. For an arbitrary triangle, the corresponding value is $\Gamma(T)$ and its angular integrand is $\Psi$.

The difference $V(h)-V(b)$ is the amount available for subtraction without falling below the
proposed bound $v=V(b)$. For $t\in[a,b]$ we have $c_b(t)=0$, so $F(b)=V(b)$ for every $\eta$.
Condition (B) requires the near-side triangle $T_a$ and the far-side triangle $T_b$ to have the same net cost $v$. Thus $F(a)$ must also equal $V(b)$, so the weighted cap at $h=a$ must equal $V(a)-V(b)$. On the rest of the symmetric
family, the desired lower bound requires $F(h)\ge V(b)$.

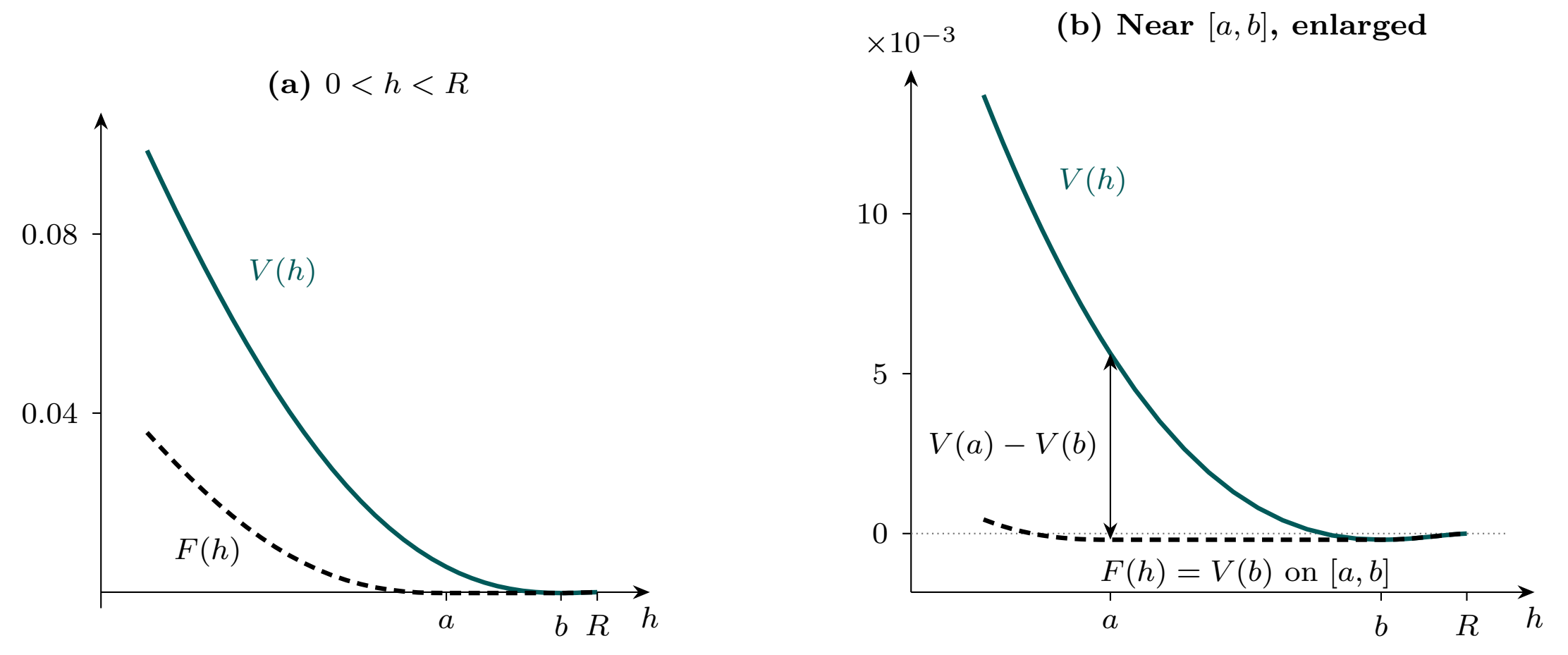
\begin{figure}[t]
\centering
\begin{adjustbox}{max width=\linewidth}
\begin{tikzpicture}

\draw[->] (0,0) -- (5.15,0) node[below,lab]{$h$};
\draw[->] (0,-0.15) -- (0,4.5);
\draw[acc,very thick] plot coordinates {(0.432,4.144) (0.566,3.862) (0.699,3.586) (0.833,3.318) (0.967,3.057) (1.101,2.805) (1.234,2.562) (1.368,2.329) (1.502,2.106) (1.635,1.893) (1.769,1.691) (1.903,1.501) (2.037,1.322) (2.170,1.154) (2.304,0.998) (2.438,0.854) (2.571,0.722) (2.705,0.602) (2.839,0.493) (2.973,0.396) (3.106,0.311) (3.240,0.237) (3.360,0.179) (3.479,0.131) (3.599,0.090) (3.718,0.057) (3.838,0.032) (3.957,0.013) (4.077,0.001) (4.196,-0.006) (4.316,-0.008) (4.353,-0.008) (4.391,-0.007) (4.429,-0.007) (4.467,-0.006) (4.505,-0.004) (4.543,-0.003) (4.580,-0.002) (4.618,-0.001) (4.656,-0.000)};
\draw[very thick,densely dashed] plot coordinates {(0.432,1.498) (0.566,1.352) (0.699,1.211) (0.833,1.076) (0.967,0.948) (1.101,0.827) (1.234,0.714) (1.368,0.609) (1.502,0.512) (1.635,0.424) (1.769,0.344) (1.903,0.274) (2.037,0.211) (2.170,0.157) (2.304,0.112) (2.438,0.075) (2.571,0.045) (2.705,0.023) (2.839,0.007) (2.973,-0.003) (3.106,-0.007) (3.240,-0.008) (3.360,-0.008) (3.479,-0.008) (3.599,-0.008) (3.718,-0.008) (3.838,-0.008) (3.957,-0.008) (4.077,-0.008) (4.196,-0.008) (4.316,-0.008) (4.353,-0.008) (4.391,-0.007) (4.429,-0.007) (4.467,-0.006) (4.505,-0.004) (4.543,-0.003) (4.580,-0.002) (4.618,-0.001) (4.656,-0.000)};
\foreach \x/\l in {3.240/a,4.316/b,4.656/R} {\draw[thin] (\x,0) -- ++(0,-0.08) node[below,lab]{$\l$};}
\foreach \y/\l in {1.68/0.04,3.36/0.08} {\draw[thin] (0,\y) -- ++(-0.08,0) node[left,lab]{$\l$};}
\node[lab,acc,anchor=west] at (1.25,3.0) {$V(h)$};
\node[lab,anchor=west] at (0.55,0.38) {$F(h)$};
\node[lab,font=\scriptsize\bfseries] at (2.5,4.75) {(a) $0<h<R$};

\begin{scope}[shift={(7.6,0.55)}]
\draw[->] (0,-0.55) -- (5.85,-0.55) node[below,lab]{$h$};
\draw[->] (0,-0.55) -- (0,4.35);
\draw[black!55,thin,densely dotted] (0,0) -- (5.6,0);
\draw[acc,very thick] plot coordinates {(0.680,4.114) (0.772,3.886) (0.863,3.665) (0.955,3.451) (1.046,3.244) (1.138,3.044) (1.229,2.851) (1.321,2.665) (1.412,2.486) (1.504,2.313) (1.595,2.147) (1.687,1.988) (1.778,1.836) (1.870,1.690) (2.101,1.351) (2.332,1.053) (2.563,0.794) (2.793,0.574) (3.024,0.390) (3.255,0.241) (3.486,0.125) (3.717,0.040) (3.948,-0.017) (4.179,-0.048) (4.410,-0.058) (4.483,-0.057) (4.556,-0.055) (4.629,-0.051) (4.702,-0.046) (4.775,-0.039) (4.848,-0.032) (4.921,-0.025) (4.994,-0.017) (5.067,-0.010) (5.140,-0.004) (5.214,-0.000)};
\draw[very thick,densely dashed] plot coordinates {(0.680,0.132) (0.772,0.098) (0.863,0.068) (0.955,0.041) (1.046,0.019) (1.138,-0.001) (1.229,-0.017) (1.321,-0.030) (1.412,-0.040) (1.504,-0.048) (1.595,-0.053) (1.687,-0.056) (1.778,-0.058) (1.870,-0.058) (2.101,-0.058) (2.332,-0.058) (2.563,-0.058) (2.793,-0.058) (3.024,-0.058) (3.255,-0.058) (3.486,-0.058) (3.717,-0.058) (3.948,-0.058) (4.179,-0.058) (4.410,-0.058) (4.483,-0.057) (4.556,-0.055) (4.629,-0.051) (4.702,-0.046) (4.775,-0.039) (4.848,-0.032) (4.921,-0.025) (4.994,-0.017) (5.067,-0.010) (5.140,-0.004) (5.214,-0.000)};
\draw[<->,thin] (1.870,-0.058) -- (1.870,1.690);
\node[lab,anchor=east,inner sep=1pt] at ($(1.870,0.82)+(-0.07,0)$) {$V(a)-V(b)$};
\foreach \x/\l in {1.870/a,4.410/b,5.214/R} {\draw[thin] (\x,-0.55) -- ++(0,-0.08) node[below,lab]{$\l$};}
\foreach \y/\l in {0/0,1.5/5,3.0/10} {\draw[thin] (0,\y) -- ++(-0.08,0) node[left,lab]{$\l$};}
\node[lab,anchor=south] at (0,4.35) {$\times 10^{-3}$};
\node[lab,acc,anchor=west] at (1.25,3.3) {$V(h)$};
\node[lab,anchor=north] at (3.14,-0.1) {$F(h)=V(b)$ on $[a,b]$};
\node[lab,font=\scriptsize\bfseries] at (3.1,4.75) {(b) Near $[a,b]$, enlarged};
\end{scope}
\end{tikzpicture}
\end{adjustbox}
\caption{The functions $V$ (solid) and $F$ (dashed) on the symmetric family, computed for $d=0.9$
($a=0.45$, $b\approx0.599$, $R\approx0.647$). Before the overhang subtraction, $V(a)$ exceeds $V(b)\approx-0.19\times10^{-3}$ by
$V(a)-V(b)\approx5.8\times10^{-3}$. The weighted caps remove this excess, making
$F(h)=\Gamma(T_h)$ equal to $V(b)$ on $[a,b]$ and at least $V(b)$ throughout $0<h<R$.}
\label{fig:VF}
\end{figure}

\paragraph{The condition on $\eta$.}
To select $\eta$, require the net cost $F(h)=\Gamma(T_h)$ to be constant on $[a,b]$. Equivalently, require the weighted cap to equal $V(h)-V(b)$ throughout this interval:
\begin{equation}\label{eq:flat}
\int_a^b\eta(t)\,c_h(t)\dd t=V(h)-V(b),\qquad\text{that is,}\quad F(h)=V(b),\qquad a\le h\le b .
\end{equation}
Agreement at $h=a$ alone would give (B), since $F(b)=V(b)$ for every $\eta$. We impose \eqref{eq:flat} on the whole interval since the proof of (A) in Section~\ref{sec:alltriangles} uses that $F$ is nondecreasing on $[a,R)$. Since $F(a)=F(b)$, a nondecreasing $F$ is constant on $[a,b]$.

\paragraph{Solving for $\eta$.}
For $a\le h\le b$, only radii $t>h$ contribute to the weighted cap, so the integral in \eqref{eq:flat} starts at $h$. Differentiating with \eqref{eq:capderiv} gives the slope equation below. Conversely, integrating that equation from $h$ to $b$ recovers \eqref{eq:flat}, since the cap integral and $V(h)-V(b)$ both vanish at $h=b$. Thus \eqref{eq:flat} is equivalent to
\begin{equation}\label{eq:slope}
2\int_h^b\eta(t)\sqrt{t^2-h^2}\dd t=-V'(h),\qquad a\le h\le b .
\end{equation}
The left-hand side is the weighted cap area lost per unit of outward movement. The right-hand side is the decrease in $V(h)$ before that cap area is subtracted. The two effects cancel, so
$F(h)=\Gamma(T_h)$ remains constant between $a$ and $b$. For $a<h<b$, a further differentiation
gives the Abel equation \citep{GorenfloVessella1991}
\begin{equation}\label{eq:abel}
\int_h^b\frac{\eta(t)}{\sqrt{t^2-h^2}}\dd t=\frac{V''(h)}{2h},\qquad a<h<b.
\end{equation}
In the variables $h^2$ and $t^2$, this is Abel's equation in standard form. Lemma~\ref{lem:eta} gives its solution explicitly, and \appEta{} verifies directly that this function satisfies \eqref{eq:slope}.

\Needspace{4\baselineskip}
\begin{restatable}[The weight]{lemma}{overhangWeight}\label{lem:eta}
Let $q(y)=V''(\sqrt y)/(2\sqrt y)$, that is,
\[
q(y)=\arcosh\!\Bigl(\frac{R}{\sqrt y}\Bigr)-\frac{R-S}{\sqrt{R^2-y}},\qquad 0<y\le b^2,
\]
and for $0<t<b$ define
\begin{equation}\label{eq:eta}
\eta(t)=\frac{2t}{\pi}\biggl[\frac{q(b^2)}{\sqrt{b^2-t^2}}-\int_{t^2}^{b^2}\frac{q'(y)}{\sqrt{y-t^2}}\dd y\biggr].
\end{equation}
Set $\eta(b)=0$ by convention. This endpoint value does not affect the integrals. Then $\eta(t)\ge 0$ for $0<t<b$, $\int_a^b\eta(t)\dd t<\infty$, and
$2\int_h^b\eta(t)\sqrt{t^2-h^2}\dd t=-V'(h)$ for all $0<h\le b$.
\end{restatable}

The proof is in \appEta. Formula~\eqref{eq:eta} is the Abel inverse of \eqref{eq:abel}: it gives the weight whose cap derivative cancels $V'(h)$ in \eqref{eq:slope}. The endpoint term and the integral involving $-q'$ are both nonnegative, since $q(b^2)>0$ and $q'<0$. The singularity at $b$ is integrable. Only $\eta$ on $[a,b]$ is used in \eqref{eq:Gamma}. Its extension to $(0,a)$ lets us compare the truncated cap integral with the full slope identity when $h<a$.

\paragraph{Consequences for the symmetric family.}
Use this function $\eta$ from now on. Integrating the slope identity from $h$ to $b$ gives \eqref{eq:flat}, so $F(h)=V(b)$ on $[a,b]$ and the reference triangles satisfy (B). For $b<h<R$, all caps in \eqref{eq:F} vanish, so $F(h)=V(h)\ge V(b)$ by Lemma~\ref{lem:V}. For $0<h<a$, differentiate \eqref{eq:F} and subtract the full slope identity in Lemma~\ref{lem:eta}:
\[
F'(h)=V'(h)+2\int_a^b\eta(t)\sqrt{t^2-h^2}\dd t
=-2\int_h^a\eta(t)\sqrt{t^2-h^2}\dd t\le0.
\]
Thus $F$ decreases toward $F(a)=V(b)$ on $(0,a)$. For $a\le h\le b$, differentiating \eqref{eq:F} and using \eqref{eq:slope} gives
\[
F'(h)=V'(h)+2\int_h^b\eta(t)\sqrt{t^2-h^2}\dd t=0.
\]
This calculation includes $h=a$ and $h=b$, since $|\partial_h c_h(t)|\le2b$ and $\eta$ is integrable. For $b<h<R$, the caps vanish and $F'(h)=V'(h)>0$ by the sign calculation in \appV. If $a=b$, then $F=V$ and $V'(b)=0$. Together these cases give
\[
F(h)\ge V(b)\quad(0<h<R),\qquad
F(h)<0,\quad F'(h)\ge0\quad(a\le h<R),\qquad F(R)=0,
\]
where the endpoint value is defined by continuity. At $R$, formula~\eqref{eq:app-Vprime} gives $F'(h)=V'(h)\to0$ as $h\uparrow R$, so the mean value theorem gives the left derivative $F'_-(R)=0$.

At $h=a$, equations~\eqref{eq:rK} and \eqref{eq:slope}, with $d=2a$, give
\[
\rs(d)=2\int_a^b\eta(t)\sqrt{t^2-a^2}\dd t=-V'(a).
\]
Thus $\rs(d)$ is the rate at which $V(h)$ decreases as $h$ increases from $a$, with $R$, $S$ and $\lambda$ held fixed.

\subsection{Arbitrary triangles}\label{sec:alltriangles}

It remains to prove condition (A), $\Gamma(T)\ge v$, for arbitrary $h>0$ and angular endpoints allowed by \eqref{eq:tri}. For fixed $h$, equation~\eqref{eq:Psi} bounds $\Gamma(T)$ below by the integral over all directions where $\Psi(h\sec\theta)<0$. We determine those directions from the symmetric-family values $F$.

For a symmetric triangle, the substitution $s=h\sec\theta$
gives
\begin{equation}\label{eq:FPsi}
F(h)=\int_{-\arccos(h/R)}^{\arccos(h/R)}\Psi(h\sec\theta)\dd\theta
=2h\int_h^R\frac{\Psi(s)}{s\sqrt{s^2-h^2}}\dd s,\qquad 0<h<R,
\end{equation}
which is an Abel transform of $\Psi$. We use its inverse below to determine the sign of $\Psi$ on $(0,R)$.

\Needspace{4\baselineskip}
\begin{restatable}[Side-triangle inequality]{lemma}{sideTriangleInequality}\label{lem:triangle}
Let $S$, $\lambda$, $\eta$ be as in \eqref{eq:S}, \eqref{eq:lambda}, \eqref{eq:eta}, and $\Psi$ as in
\eqref{eq:Psi}. For every $h>0$ and every $-\pi/2<\theta_-<\theta_+<\pi/2$,
\begin{equation}\label{eq:triineq}
\int_{\theta_-}^{\theta_+}\Psi(h\sec\theta)\dd\theta\;\ge\;V(b).
\end{equation}
Equality holds when $a\le h\le b$ and $\theta_\pm=\pm\arccos(h/R)$, in particular for the six side
triangles of $H_d$.
\end{restatable}

Appendix~\ref{app:triangle} supplies the regularity and endpoint details for the inversion.

\begin{figure}[t]
\centering
\begin{adjustbox}{max width=\linewidth}
\begin{tikzpicture}
\fill[black!22] (1.477,-0.07) rectangle (6.020,0.07);
\draw[->] (0,0) -- (7.45,0) node[right,lab]{$s$};
\draw[->] (0,-1.25) -- (0,4.45);
\node[lab,anchor=south] at (0,4.45) {$\times10^{-3}$};
\foreach \y/\l in {0/0,2.4/1} {\draw[thin] (0,\y) -- ++(-0.08,0) node[left,lab]{$\l$};}
\foreach \x/\l in {1.477/\rho,1.890/a,5.027/b,6.020/R}
{\draw[black!50,thin,densely dotted] (\x,0) -- (\x,-1.25); \node[lab,anchor=north] at (\x,-1.25) {$\l$};}
\draw[very thick,densely dashed] plot coordinates {(3.558,4.200) (3.656,3.702) (3.754,3.228) (3.852,2.779) (3.950,2.353) (4.048,1.953) (4.146,1.579) (4.244,1.230) (4.342,0.908) (4.440,0.613) (4.538,0.346) (4.635,0.106) (4.733,-0.105) (4.831,-0.287) (4.929,-0.439) (5.027,-0.562)};
\draw[acc,very thick] plot coordinates {(0.000,3.620) (0.111,3.215) (0.222,2.829) (0.334,2.463) (0.445,2.117) (0.556,1.792) (0.667,1.489) (0.778,1.208) (0.889,0.950) (1.001,0.716) (1.112,0.507) (1.223,0.323) (1.334,0.164) (1.445,0.033) (1.556,-0.071) (1.668,-0.147) (1.779,-0.194) (1.890,-0.212) (2.075,-0.217) (2.259,-0.222) (2.444,-0.227) (2.628,-0.233) (2.813,-0.240) (2.997,-0.248) (3.182,-0.256) (3.366,-0.266) (3.551,-0.276) (3.735,-0.289) (3.920,-0.304) (4.104,-0.321) (4.289,-0.343) (4.474,-0.371) (4.658,-0.407) (4.843,-0.461) (5.027,-0.562) (5.117,-0.647) (5.208,-0.707) (5.298,-0.740) (5.388,-0.746) (5.479,-0.725) (5.569,-0.676) (5.659,-0.599) (5.749,-0.493) (5.840,-0.358) (5.930,-0.194) (6.020,0.000) (6.121,0.253) (6.222,0.544) (6.324,0.873) (6.425,1.242) (6.526,1.651) (6.627,2.100) (6.728,2.590) (6.829,3.122) (6.930,3.695)};
\node[lab,anchor=east,inner sep=1pt] (r0) at (4.1,0.6) {$\rho_0$};
\draw[thin] (r0.east) -- ($(4.683,0)+(-0.03,0.05)$);
\node[lab,acc,anchor=west,align=left] at (1.05,1.55) {$\Psi(s)$, after\\ overhang subtraction};
\node[lab,anchor=west] at (3.95,3.6) {without it ($\eta=0$)};
\node[lab,anchor=west,text width=4.6cm] at (8.0,3.05)
{Without the overhang subtraction, the angular contribution is negative only for lengths between $\rho_0$ and $R$.
At $s=a$, this unweighted contribution is positive.};
\node[lab,anchor=west,text width=4.6cm] at (8.0,1.15)
{The subtraction lowers the angular contribution for shorter rays and widens that range to $(\rho,R)$, with
$\rho<a$.};
\end{tikzpicture}
\end{adjustbox}
\caption{The angular integrand $\Psi(s)$ for radial endpoint distance $s$, computed for $d=0.9$
($\rho\approx0.430$, $a=0.45$, $\rho_0\approx0.583$, $b\approx0.599$, $R\approx0.647$). The dashed curve is $\Psi$ with $\eta=0$, and $\rho_0$ is its zero below $R$. The two curves
coincide for $s\ge b$. After the weighted overhang subtraction, $\Psi$ is negative exactly on the shaded interval
$(\rho,R)$, as in \eqref{eq:sign}.}
\label{fig:psi}
\end{figure}
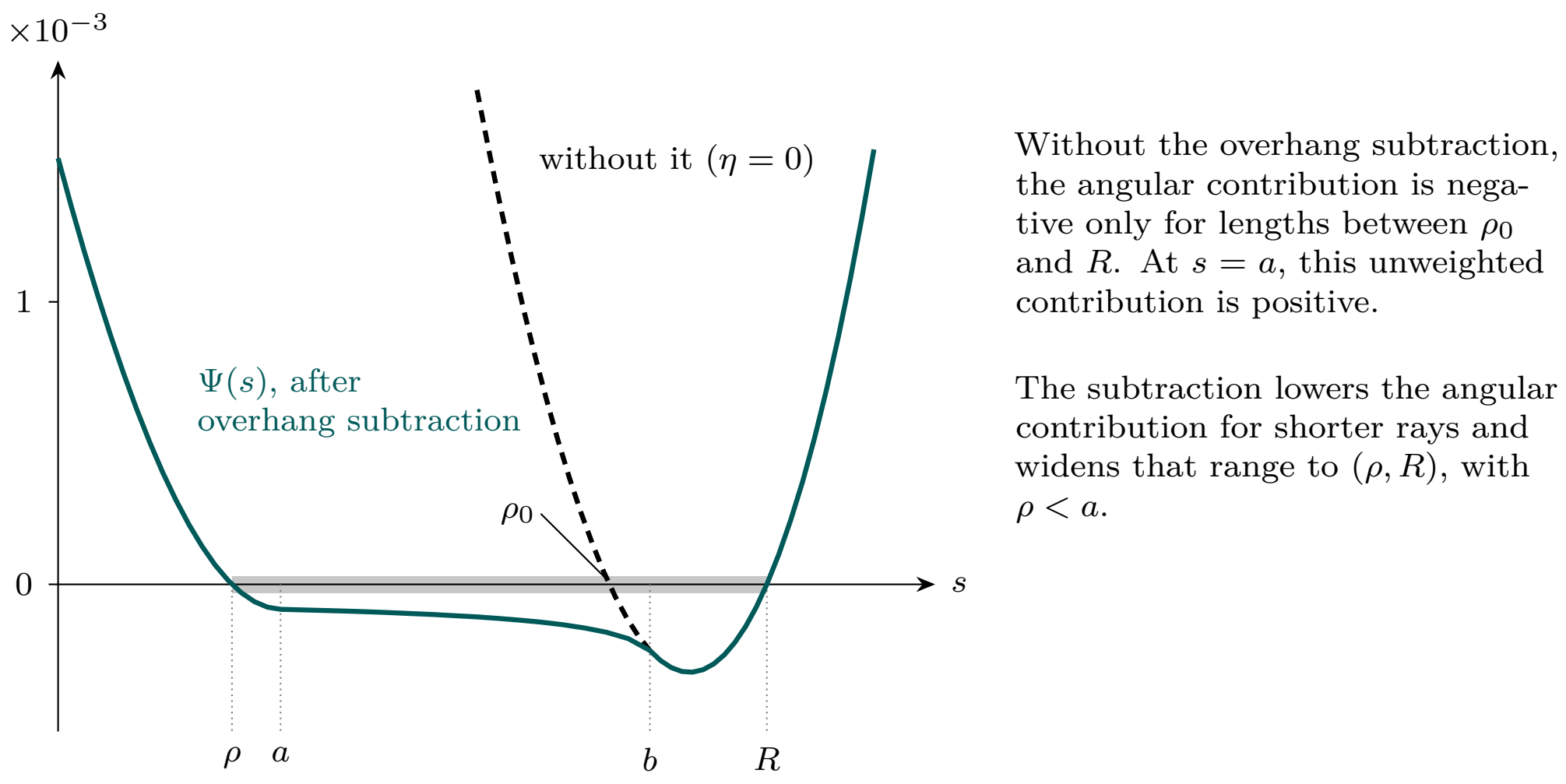

\paragraph{Step 1: the sign of $\Psi$.}
The function $\Psi$ is negative on a single interval (Figure~\ref{fig:psi}): there is a radius $\rho$ with $0\le\rho<a$ such that
\begin{equation}\label{eq:sign}
\Psi\ge 0\ \text{on }(0,\rho),\qquad \Psi<0\ \text{on }(\rho,R),\qquad \Psi>0\ \text{on }(R,\infty).
\end{equation}
For $s>R>b$, we have $(t^2-s^2)_+=0$ for every $t\in[a,b]$, so the integral term vanishes. Then $\Psi(R)=0$ and $\Psi'(s)=s(s-S)>0$ give the positive sign.
For $a\le s<R$, use the relation \eqref{eq:FPsi} between $F$ and $\Psi$. With
$x=h^2$ and $y=s^2$, it becomes the Abel transform
\[
\frac{F(\sqrt x)}{\sqrt x}
=\int_x^{R^2}\frac{\Psi(\sqrt y)/y}{\sqrt{y-x}}\dd y .
\]
Its inverse, using $F(R)=0$, is
\[
\Psi(s)=-\frac{s^2}{\pi}\int_s^R
\frac{F'(u)/u-F(u)/u^2}{\sqrt{u^2-s^2}}\dd u,\qquad 0<s<R .
\]
For $a\le s<R$, the symmetric-family bounds give $F'(u)\ge0$ and $F(u)<0$ for $s<u<R$. The integrand is therefore positive, and the leading minus sign gives $\Psi(s)<0$.

For $0<s<a$, direct differentiation gives
$\Psi'(s)=s\bigl(s-S+\int_a^b\eta(t)\dd t\bigr)$. The factor in parentheses increases with $s$,
so $\Psi$ can only decrease and then increase on this interval. Since $\Psi(a)<0$, there is at most
one zero in $(0,a)$, crossed from positive to negative. Denote it by $\rho$, or set $\rho=0$ when
there is no such zero. This proves the full sign pattern \eqref{eq:sign}.

\begin{figure}[t]
\centering
\begin{adjustbox}{max width=\linewidth}
\begin{tikzpicture}[scale=2.0,rayneg/.style={acc,line width=0.55pt},raypos/.style={orange!85!black,line width=0.55pt,densely dashed}]
\begin{scope}
\clip (-0.45,-1.36) rectangle (3.0,1.5);
\fill[soft] (0,0) -- (0.75,-0.6614) -- (0.75,0.6614) -- cycle;
\draw[raypos] (0,0) -- (0.75,-1.2990);
\draw[raypos] (0,0) -- (0.75,-1.0323);
\draw[raypos] (0,0) -- (0.75,-0.8330);
\draw[raypos] (0,0) -- (0.75,-0.6753);
\draw[rayneg] (0,0) -- (0.75,-0.5449);
\draw[rayneg] (0,0) -- (0.75,-0.4330);
\draw[rayneg] (0,0) -- (0.75,-0.3339);
\draw[rayneg] (0,0) -- (0.75,-0.2437);
\draw[rayneg] (0,0) -- (0.75,-0.1594);
\draw[rayneg] (0,0) -- (0.75,-0.0788);
\draw[rayneg] (0,0) -- (0.75,0.0000);
\draw[rayneg] (0,0) -- (0.75,0.0788);
\draw[rayneg] (0,0) -- (0.75,0.1594);
\draw[rayneg] (0,0) -- (0.75,0.2437);
\draw[rayneg] (0,0) -- (0.75,0.3339);
\draw[rayneg] (0,0) -- (0.75,0.4330);
\draw[rayneg] (0,0) -- (0.75,0.5449);
\draw[raypos] (0,0) -- (0.75,0.6753);
\draw[raypos] (0,0) -- (0.75,0.8330);
\draw[raypos] (0,0) -- (0.75,1.0323);
\draw[raypos] (0,0) -- (0.75,1.2990);
\draw[black!45,thin,densely dotted] (-80:1) arc[start angle=-80,end angle=80,radius=1];
\draw[black!45,thin,densely dotted] (-80:0.62) arc[start angle=-80,end angle=80,radius=0.62];
\draw[black!60,thin] (0.75,-1.34) -- (0.75,1.34);
\draw[very thick] (0,0) -- (0.75,-0.2437) -- (0.75,1.1119) -- cycle;
\node[dot] at (0,0) {};
\node[lab,left] at (0,0) {$x_i$};
\node[lab] at (-76:1.08) {$R$};
\node[lab] at (-76:0.54) {$\rho$};
\node[lab,anchor=west,text width=3.3cm] at (0.84,0) {rays shorter than $R$: $\Psi<0$. These form the symmetric triangle (shaded)};
\node[lab,anchor=west] at (0.82,1.0) {longer than $R$: $\Psi>0$};
\node[lab,anchor=west] at (0.82,-1.0) {longer than $R$: $\Psi>0$};
\node[lab,anchor=west,text width=2.5cm] at (1.6,0.66) {another triangle on the same line (thick outline)};
\draw[thin] (1.58,0.68) -- (0.76,0.9);
\node[lab,font=\scriptsize\bfseries] at (0.9,1.42) {(a) $\rho\le h<R$};
\end{scope}
\begin{scope}[shift={(3.7,0)}]
\clip (-0.25,-1.36) rectangle (2.7,1.5);
\fill[soft] (0,0) -- (0.4,0.4739) -- (0.4,0.9165) -- cycle;
\fill[soft] (0,0) -- (0.4,-0.4739) -- (0.4,-0.9165) -- cycle;
\draw[raypos] (0,0) -- (0.4,-1.2686);
\draw[raypos] (0,0) -- (0.4,-0.9657);
\draw[rayneg] (0,0) -- (0.4,-0.7684);
\draw[rayneg] (0,0) -- (0.4,-0.6279);
\draw[rayneg] (0,0) -- (0.4,-0.5213);
\draw[raypos] (0,0) -- (0.4,-0.4365);
\draw[raypos] (0,0) -- (0.4,-0.3665);
\draw[raypos] (0,0) -- (0.4,-0.3069);
\draw[raypos] (0,0) -- (0.4,-0.2548);
\draw[raypos] (0,0) -- (0.4,-0.2082);
\draw[raypos] (0,0) -- (0.4,-0.1657);
\draw[raypos] (0,0) -- (0.4,-0.1261);
\draw[raypos] (0,0) -- (0.4,-0.0887);
\draw[raypos] (0,0) -- (0.4,-0.0527);
\draw[raypos] (0,0) -- (0.4,-0.0175);
\draw[raypos] (0,0) -- (0.4,0.0175);
\draw[raypos] (0,0) -- (0.4,0.0527);
\draw[raypos] (0,0) -- (0.4,0.0887);
\draw[raypos] (0,0) -- (0.4,0.1261);
\draw[raypos] (0,0) -- (0.4,0.1657);
\draw[raypos] (0,0) -- (0.4,0.2082);
\draw[raypos] (0,0) -- (0.4,0.2548);
\draw[raypos] (0,0) -- (0.4,0.3069);
\draw[raypos] (0,0) -- (0.4,0.3665);
\draw[raypos] (0,0) -- (0.4,0.4365);
\draw[rayneg] (0,0) -- (0.4,0.5213);
\draw[rayneg] (0,0) -- (0.4,0.6279);
\draw[rayneg] (0,0) -- (0.4,0.7684);
\draw[raypos] (0,0) -- (0.4,0.9657);
\draw[raypos] (0,0) -- (0.4,1.2686);
\draw[black!45,thin,densely dotted] (-85:1) arc[start angle=-85,end angle=85,radius=1];
\draw[black!45,thin,densely dotted] (-85:0.62) arc[start angle=-85,end angle=85,radius=0.62];
\draw[black!60,thin] (0.4,-1.34) -- (0.4,1.34);
\node[dot] at (0,0) {};
\node[lab] at (-80:1.08) {$R$};
\node[lab] at (-27:0.72) {$\rho$};
\node[lab,anchor=west] at (0.5,0.70) {between $\rho$ and $R$: $\Psi<0$};
\node[lab,anchor=west] at (0.5,-0.70) {between $\rho$ and $R$: $\Psi<0$};
\node[lab,anchor=west] at (0.46,0.0) {shorter than $\rho$: $\Psi\ge 0$};
\node[lab,anchor=west] at (0.46,1.15) {longer than $R$: $\Psi>0$};
\node[lab,font=\scriptsize\bfseries] at (0.8,1.42) {(b) $h<\rho$};
\end{scope}
\end{tikzpicture}
\end{adjustbox}
\caption{Rays from a facility to a line at distance $h$. The angular contribution $\Psi(s)$ depends only on the radial endpoint distance $s$: it is negative for solid rays and
nonnegative for dashed rays. (a) For $\rho\le h<R$ the rays with negative angular contribution are those of the
symmetric triangle with corners at distance $R$. Any other triangle with its side on the line, such as
the one outlined, omits some of them or contains longer rays, so its net cost is at least $F(h)$.
(b) For $h<\rho$ the rays near the foot of the perpendicular are shorter than $\rho$, and the rays with
negative angular contribution form two bands. Their combined contribution is at least $F(\rho)$. The drawing is schematic.}
\label{fig:signs}
\end{figure}
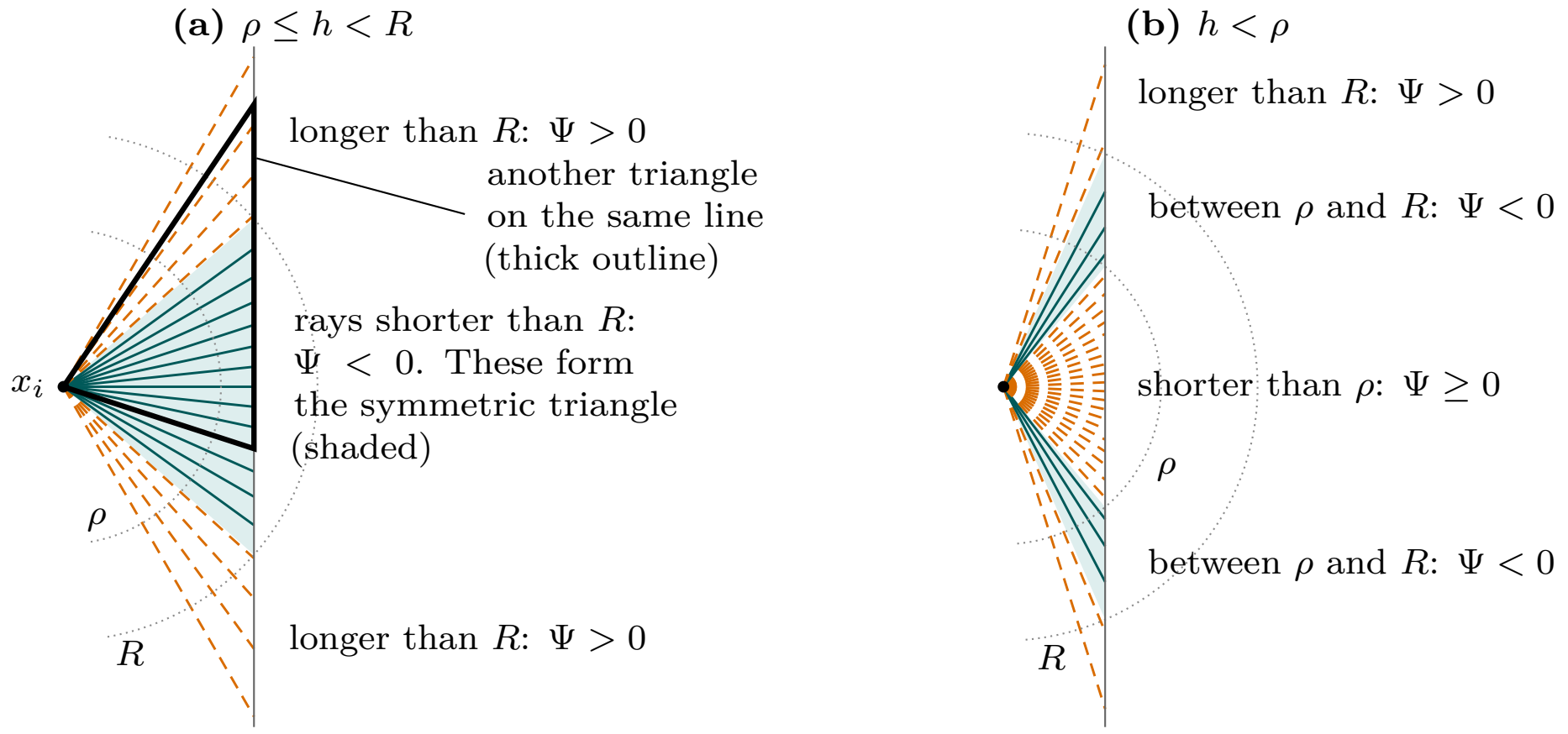

\paragraph{Step 2: an arbitrary interval of directions.}
Fix $h$. The integral in \eqref{eq:triineq} is at least the integral of $\Psi(h\sec\theta)$ over all
directions $\theta\in(-\pi/2,\pi/2)$ in which it is negative: leaving out directions with negative
values, or including directions with nonnegative values, can only increase an integral. By \eqref{eq:sign},
the negative directions are those with $\rho<h\sec\theta<R$ (Figure~\ref{fig:signs}).
\begin{itemize}[topsep=3pt,itemsep=1pt]
\item If $h\ge R$, there are no directions with negative integrand, so the integral is at least $0>V(b)$.
\item If $\rho\le h<R$, the directions with negative integrand form the interval $\abs{\theta}<\arccos(h/R)$, except possibly for the single zero direction $\theta=0$ when $h=\rho$. This is the angular interval of the symmetric triangle $T_h$, and its integral is $F(h)\ge V(b)$.
\item If $0<h<\rho$, the directions with negative integrand form two intervals symmetric about $\theta=0$. With $s=h\sec\theta$, their combined contribution is
\[
2\int_\rho^R\Psi(s)\frac{h}{s\sqrt{s^2-h^2}}\dd s.
\]
For fixed $s>\rho$, the factor $h/(s\sqrt{s^2-h^2})$ increases with $h$ on $(0,\rho)$. Since $\Psi(s)<0$ for $\rho<s<R$, replacing $h$ by $\rho$ makes the integral no larger. Its resulting value is $F(\rho)\ge V(b)$.
\end{itemize}
In the last case, the two intervals do not describe a single triangle. Their combined contribution is used only as a lower bound. For the reference triangles, the interval is symmetric, and the integral equals $F(a)=F(b)=V(b)$, which gives the equality cases.

\multiplierPrices*
\begin{proof}
Take $S$, $\lambda$ and $\eta$ from \eqref{eq:S}, \eqref{eq:lambda} and \eqref{eq:eta}, and $v=V(b)$.
Then $v<0$ by Lemma~\ref{lem:V}, and $\eta$ is nonnegative and integrable on $[a,b]$ by
Lemma~\ref{lem:eta}. By \eqref{eq:Psi}, the net cost of a triangle is the integral in
\eqref{eq:triineq}. Lemma~\ref{lem:triangle} therefore states condition (A), and its equality cases
include condition (B).
\end{proof}

\subsection{Completion of the proof of Theorem 2}\label{sec:completion}
Fix the multipliers from Lemma~\ref{lem:prices} and the quantities $\rs(d)$ and $K(d)$ in \eqref{eq:rK}. We first apply their lower bound to a layout with interior facilities.

\Needspace{4\baselineskip}
\begin{restatable}[Service--routing trade-off]{lemma}{serviceRoutingTradeoff}\label{lem:tradeoff}
For every layout $X$ of $N$ distinct facilities in the interior of $Q_N$ and every closed tour of
length $L$, inequalities \eqref{eq:tradeoff} and \eqref{eq:tradeoff-total} hold.
\end{restatable}

\begin{proof}
Consider the relaxed problem with routing weight $r=\rs(d)$. The actual side triangles and tour-edge lengths satisfy \eqref{eq:relaxed}, and (A) holds by Lemma~\ref{lem:prices}, so \eqref{eq:dualbound} applies with $K=K(d)$. By \eqref{eq:wtangent} with $r=\rs(d)$, each $w(\ell_i)$ is at least $w(d)$. Hence
\[
M_N(X)+\rs(d)\,L\;\ge\;N\,(S-2\pi\lambda+6v)+N\,w(d)-K(d)=N\{\phi(d)+\rs(d)\,d\}-K(d),
\]
where the equality follows from \eqref{eq:hexsum} and the definition \eqref{eq:w} of $w$. This is
\eqref{eq:tradeoff-total}. Subtracting $\rs(d)L$ gives \eqref{eq:tradeoff}.
\end{proof}

The coefficient $\rs(d)$ is the weighted marginal increase in disk overlap when a tour edge of length $d$ is shortened. The allowance $K(d)$ is the term $\pi t^2$ from the tour constraint integrated against the same weight. Appendix~\ref{app:coverage} gives a related direct consequence of the side-triangle bound: a weighted excess in uncovered area forces higher service cost.

\Needspace{4\baselineskip}
\begin{restatable}[Boundary facilities]{lemma}{boundaryFacilities}\label{lem:boundary}
Inequalities \eqref{eq:tradeoff} and \eqref{eq:tradeoff-total} hold for every layout of $N$ distinct
facilities in $Q_N$, including facilities on its sides and corners, with the same $\rs(d)$ and $K(d)$
and for every closed tour.
\end{restatable}

The proof (\appBdry) moves every facility slightly toward the center $z$ of the square, to
$(1-\varepsilon)x_i+\varepsilon z$. The facilities become interior and stay distinct, the tour
length becomes $(1-\varepsilon)L$, and the service cost changes by at most
$N\varepsilon\sqrt{N/2}$. Lemma~\ref{lem:tradeoff} applies to the moved layout, and
$\varepsilon\to 0$ gives the claim.

\Needspace{4\baselineskip}
\begin{restatable}[Optimal spacing]{lemma}{optimalSpacing}\label{lem:spacing}
For $0<d\le d_0$, we have 
\[
\rs(d)=-\phi'(d),
\]
where the derivative at $d_0$ is taken from the left. Thus, for every $r\ge 0$, a spacing $d_r$
that minimizes $\phi(d)+rd$ over $(0,d_0]$ satisfies $\rs(d_r)=r$. For $r=0$, we have $d_r=d_0$.

Moreover, $\phi$ is strictly convex on $(0,d_0]$, so this minimizing spacing is unique.
\end{restatable}

Appendix~\ref{app:spacing} differentiates $\phi(d)=2J_R(a)+4J_R(b)$ with respect to $d$. Subtracting $S$ times the derivative of the cell-area identity and adding $\lambda$ times the derivative of the angle identity leaves $\phi'(d)$ unchanged. The choice of $\lambda$ cancels the terms containing $R'(d)$, and $V'(b)=0$ cancels the term containing $b'(d)$. The remaining term is $2V'(a)a'(d)=V'(a)=-\rs(d)$.

For $r>0$, the left derivative of $\phi(d)+rd$ at $d_0$ is $r>0$, so the minimizer is interior and satisfies $\phi'(d_r)+r=0$. For $r=0$, the function $\phi$ decreases up to $d_0$. In either case, the minimizing spacing satisfies $\rs(d_r)=r$. Strict convexity makes this choice unique.

\universalLowerBound*
\begin{proof}
Fix $r\ge 0$ and let $d_r$ be a minimizing spacing in \eqCsix, which exists by \lemSpacing. Apply
Lemmas~\ref{lem:tradeoff} and~\ref{lem:boundary} with $d=d_r$. Since $\rs(d_r)=r$ by
Lemma~\ref{lem:spacing}, inequality \eqref{eq:tradeoff-total} reads
\begin{equation}\label{eq:final}
M_N(X)+rL\;\ge\;N\{\phi(d_r)+r\,d_r\}-K(d_r)=N\,C_6(r)-K_r,\qquad K_r:=K(d_r).
\end{equation}
The same calculation and the feasible reference collection in Section~\ref{sec:right} also give
\[
NC_6(r)-K_r\le\text{optimal value of \eqref{eq:relaxed}}\le NC_6(r).
\]
The relaxation gives the lower bound. Theorem~\ref{thm:upper} provides layouts and closed tours that attain the same limiting cost per facility.

The spacing $d_r$, and hence $K_r$, depends on $r$ only. For $r=0$ we have $d_r=d_0$ and $a=b$, so
$K_0=0$.
\end{proof}

\section{Numerical implications of the bounds}\label{sec:numerics}

We now evaluate the minimizing spacing and the constants in the bounds. The calculations also compare $C_6(r)$ with the earlier relaxed bound and illustrate the gap for finite networks.

\subsection{Optimal spacing and cell shape}\label{sec:numericaldesign}
Table~\ref{tab:values} lists the optimal spacing $d_r$, cell height $2R$, aspect ratio $2R/d_r$, limiting cost $C_6(r)$, and constants $K_r$ and $D_r$ for several routing weights. The ``saving'' column gives the percentage reduction from $\phi(d_0)+rd_0$, the cost per facility of staggered rows of regular hexagons. For $r=0.5$, the optimal cell is $2.8$ times as high as it is wide, and its cost is $17\%$ below that of regular hexagons.

We compute the spacing from $\rs(d_r)=r$, and we check the entries at 50-digit working precision. We evaluate the constant $K_r$ both using the regularized weight formula and independently by interchanging the integrals in \eqref{eq:eta}.

\begin{table}[!ht]
\centering
\caption{Optimal spacing, cell shape, limiting cost and constants. Lengths are in units in which the facility density is one.}\label{tab:values}
\begin{tabular}{@{}rrrrrrrr@{}}
\toprule
$r$ & $d_r$ & $2R$ & $2R/d_r$ & $C_6(r)$ & Saving (\%) & $K_r$ & $D_r$\\
\midrule
0& 1.075 & 1.241 & 1.15 & 0.3772 & 0.0& 0 & 8.1\\
0.05 & 0.958 & 1.263 & 1.32 & 0.4277 & 0.7& 0.10& 8.2\\
0.1& 0.885 & 1.303 & 1.47 & 0.4737 & 2.3& 0.16& 8.6\\
0.25 & 0.744 & 1.447 & 1.95 & 0.5948 & 7.9& 0.31& 10.4\\
0.5& 0.610 & 1.697 & 2.78 & 0.7624 & 16.6 & 0.59& 14.2\\
1& 0.470 & 2.153 & 4.58 & 1.0282 & 29.2 & 1.28& 23.0\\
2& 0.346 & 2.901 & 8.39 & 1.4271 & 43.5 & 3.18& 41.9\\
5& 0.223 & 4.496 & 20.2 & 2.2403 & 61.0 & 11.9& 101.0\\
\bottomrule
\end{tabular}
\end{table}

\FloatBarrier
\subsection{Comparison with the earlier relaxed bound}\label{sec:earlierbound}
The idealized region of \citet[Section~3.2, equation~(17)]{XieOuyang2015} is a disk of radius $R_\alpha$ between two parallel lines at distance $R_\alpha\cos\alpha$ from its center. At unit area,
\[
R_\alpha=(\sin2\alpha+\pi-2\alpha)^{-1/2},\qquad 0\le\alpha<\pi/2.
\]
Omitting setup cost and using our routing weight gives its optimized cost
\[
\begin{aligned}
C_{\mathrm{XO}}(r)=\min_{0\le\alpha<\pi/2}\bigg\{
&\frac{2R_\alpha^3}{3}\left[\sin\alpha\cos\alpha
+\cos^3\alpha\log(\sec\alpha+\tan\alpha)+\pi-2\alpha\right]\\
&+2rR_\alpha\cos\alpha\bigg\}.
\end{aligned}
\]
The service term integrates distance over this unit-area region and $2R_\alpha\cos\alpha$ is its straight-route spacing. This is their relaxed geometry in the present normalization, before facility-density optimization. The relative difference $100[C_6(r)-C_{\mathrm{XO}}(r)]/C_6(r)$ is approximately $0.28376\%$, $0.08180\%$, $0.001496\%$ and $0.0001743\%$ at $r=0,0.05,0.5,1$, respectively. At $r=0$, the relaxed region is a disk and $C_{\mathrm{XO}}(0)=2/(3\sqrt\pi)$. At these four routing weights, the computed value $C_{\mathrm{XO}}(r)$ is below $C_6(r)$. 

\subsection{Finite-size performance}\label{sec:finiteperformance}
The constants in Theorem~\ref{thm:upper} give uniform finite-size bounds. They generally grow with $r$, although $D_r$ is not monotone near zero: $D_0\simeq8.08290$, whereas $D_{0.01}\simeq8.07034$. At $r=1$, requiring $D_r/\sqrt N\le0.1C_6(r)$ gives approximately $N\ge5.0\times10^4$. Evaluating the construction directly can give smaller gap bounds.

For $r=1$, take the unrounded minimizing spacing $d_r$ underlying Table~\ref{tab:values} and set $m=\lfloor\sqrt N/d_r\rfloor$, $n=\lfloor\sqrt N d_r\rfloor$ and $q=N-mn$. Let $M_{\mathrm{row}}$ and $L_{\mathrm{row}}$ be the service integral and closed tour length of the row-site set in Appendix~\ref{app:upperconstruction}. Insert the remaining $q$ sites near an interior row site as in that appendix, using a radius at most $1/(2N)$. Note that service cost does not increase, and the added tour length is at most $q/N$, which gives the feasible-cost bound
\[
U_N=M_{\mathrm{row}}+r(L_{\mathrm{row}}+q/N).
\]
\Needspace{6\baselineskip}
Compute $M_{\mathrm{row}}$ over the finite row-site set's clipped Voronoi cells, with an independent one-dimensional quadrature as a check. For the resulting $N$-site layout and tour, Theorem~\ref{thm:lower} gives
\[
0\le\frac{M_N(X)+rL(X,\sigma)-A_N(r)}{N}
\le\frac{U_N}{N}-C_6(r)+\frac{K_r}{N}.
\]
The right-hand side is the last column of Table~\ref{tab:finite}.

\begin{table}[!ht]
\centering
\caption{Computed finite-instance bounds at $r=1$. The last column is an upper bound on the constructed layout's excess cost per facility. Values are rounded floating-point evaluations, not interval-arithmetic certificates.}\label{tab:finite}
\begin{tabular}{@{}rrrrr@{}}
\toprule
$N$ & row sites $mn$ & $C_6-K_r/N$ & $U_N/N$ & $U_N/N-C_6+K_r/N$\\
\midrule
25 & 20 & 0.97712 & 1.21612 & 0.23900\\
100 & 84 & 1.01542 & 1.25287 & 0.23745\\
400 & 378 & 1.02500 & 1.14122 & 0.11622\\
1600 & 1530 & 1.02739 & 1.10654 & 0.07915\\
\bottomrule
\end{tabular}
\end{table}

\par\noindent Rounding the row and column counts explains the nonmonotonicity in $U_N/N$. These gap bounds are smaller than $D_r/\sqrt N+K_r/N$.
\FloatBarrier
\section{Conclusion}\label{sec:conclusion}
This paper proves the asymptotic optimality of staggered-row 
hexagonal layouts for the joint facility location and replenishment 
routing problem under uniform demand and Euclidean distance. 
For every fixed routing weight, we derive upper and lower bounds 
with the same leading term, showing that the optimality gap per 
facility vanishes at a rate of $O(N^{-1/2})$ as the service region 
and the number of facilities increase at a fixed density. We 
also characterize the unique optimal spacing within the 
staggered-row family. The main methodological contribution is 
a geometric relaxation that connects customer travel and 
replenishment routing through disk overlaps at consecutive tour 
stops. Combining the planar side-count inequality with Lagrangian 
duality and Abel integral equations, we obtain a lower bound 
that applies to arbitrary facility locations and visiting orders, 
without imposing restrictions on the shapes or areas of service 
regions. Several directions may be interesting for further investigation, including the extension to heterogeneous customer demand, capacitated 
vehicles, and multiple replenishment tours. Another question is 
whether geometric configurations outside the staggered-row family 
can attain the same limiting cost. 

\medskip
\section*{AI Use Disclosure}
The authors acknowledge the use of OpenAI's ChatGPT and Codex for manuscript editing and supporting Lean formalization. The authors are responsible for all mathematical arguments, numerical results, references, and the final manuscript.


\label{or:references-start}

\newpage
\Needspace{10\baselineskip}
\appendix
\numberwithin{equation}{section}
\section*{Appendices: Complete Proofs}
\label{ec:start}
The appendices provide the remaining geometric and analytic proofs. The proofs of Lemmas~\ref{lem:tour} and~\ref{lem:tradeoff} and Theorem~\ref{thm:lower} appear only in the main text. Appendix~\ref{app:coverage} derives a direct comparison between service cost and uncovered area.

\section{Upper-bound proofs}\label{app:upper}
\subsection{The unit-area hexagon}\label{app:hexagon}
\unitAreaHexagon*
\begin{proof}
Comparing distance to the origin with distance to the six sites $(\pm d,0)$ and $(\pm d/2,\pm1/d)$ gives the half-planes in \eqref{eq:hexagon}. Let $P$ be their intersection. Its vertices are $(0,\pm R)$ and $(\pm d/2,\pm(1/d-R))$, all at distance $R$ since $d^2/4+(1/d-R)^2=R^2$. Hence $P\subseteq B_R(0)$.

For a point $(x,y)\in P$, the origin is the nearest site in row $0$, and the nearest sites in rows $\pm1$ are among the four listed above. Every other row is at vertical distance at least $2/d-R>R$ from this point, since $dR=1/2+d^4/8\le2/3$. No other lattice site can therefore be closer than the origin, so $P=H_d$. The side distances are half the corresponding neighbor distances, giving $a,a,b,b,b,b$, and direct integration of \eqref{eq:hexagon} gives $|H_d|=2dR-d^4/4=1$.

For a symmetric side triangle let $\beta=\arccos(h/R)$. Polar integration gives
\[
\int_T\|x\|\dd x=\frac{h^3}{3}\int_{-\beta}^{\beta}\sec^3\theta\dd\theta
=\frac{h^3}{3}\{\sec\beta\tan\beta+\log(\sec\beta+\tan\beta)\}=J_R(h).
\]
Adding the six triangle integrals proves \eqref{eq:phiformula}.
\end{proof}

\subsection{Existence of a minimizing spacing}\label{app:existence}
\minimizingSpacingExists*
\begin{proof}
Compare $H_d$ with the centered rectangle of width $d$ and height $1/d$. Both have area one, so their two set differences have equal area. On the part of $H_d$ outside the rectangle, $|y|\ge1/(2d)$. On the part of the rectangle outside $H_d$, $|y|\le1/(2d)$. Thus the integral of $|y|$ over $H_d$ is at least its integral over the rectangle. Using $\sqrt{x^2+y^2}\ge|y|$ gives
\begin{equation}\label{eq:vertical}
\phi(d)\ge\int_{H_d}|y|\dd x\dd y
\ge d\int_{-1/(2d)}^{1/(2d)}|y|\dd y=\frac1{4d}.
\end{equation}
Thus $\phi(d)+rd\to\infty$ as $d\downarrow0$. Its formula is continuous, so its minimum is attained on a closed subinterval of $(0,d_0]$. At the regular endpoint,
\[
\phi(d_0)=\left(\frac{d_0}{2}\right)^3\left(\frac43+\log3\right).
\]
\end{proof}

\subsection{Finite-square upper bound}\label{app:upperconstruction}
\finiteSquareUpperBound*
\begin{proof}
Write $s=\sqrt N$, $m=\lfloor s/d\rfloor$, and $n=\lfloor sd\rfloor$. The size assumption gives $s\ge\max\{d,d^{-1}\}$, so $m,n\ge1$. Use the $mn$ interior sites
\[
x_{jk}=\left(d\left(j+\frac14+\frac{k\bmod2}{2}\right),\frac{k+1/2}{d}\right),
\qquad 0\le j<m,\quad0\le k<n.
\]
Let $X_{\mathrm{row}}$ be this finite set and write $M_{\mathrm{row}}=\int_{Q_N}\dist(x,X_{\mathrm{row}})\dd x$. Rounding leaves $q=N-mn\le(d+d^{-1})s$ facilities to be added.

Let $W$ be the union of the cells of these $mn$ facilities in the translated infinite lattice. These full hexagons have disjoint interiors, so $|W|=mn=N-q$. The part of the square not covered by these cells has area
\[
\begin{aligned}
|Q_N\setminus W|
&=|Q_N|-|Q_N\cap W|\\
&=N-\bigl(|W|-|W\setminus Q_N|\bigr)
=q+|W\setminus Q_N|.
\end{aligned}
\]
The deficit in total cell area is $q$. Any cell area outside the square adds to the uncovered part inside it.
\Needspace{5\baselineskip}
Since $md\le s$ and $n/d\le s$, the cell formula gives
\[
W\subset\left[-\frac d4,s+\frac d4\right]\times\left[-\frac{d^3}8,s+\frac{d^3}8\right].
\]
The part of this rectangle outside $Q_N$ has area $(d/2+d^3/4)s+d^4/8$. Using $s\ge d$ and $d^4\le4/3$, we obtain
\[
\begin{aligned}
|Q_N\setminus W|&\le(d+d^{-1})s+\frac{ds}{2}+\frac{d^3s}{4}+\frac{d^4}{8}\\
&\le\frac32(d+d^{-1})s,
\end{aligned}
\]
where $d^3s/4+d^4/8\le3d^3s/8\le s/(2d)$.

Every customer can reach one of these rows within vertical distance $3/(2d)$. Indeed, consecutive rows are $1/d$ apart, the lower gap is $1/(2d)$, and the upper gap is $s-n/d+1/(2d)<3/(2d)$. In that row the horizontal distance to a facility is at most $7d/4$: the left gap is at most $3d/4$ and the right gap is less than $s-md+3d/4<7d/4$. The nearest-facility distance is therefore at most $\frac74\sqrt{d^2+d^{-2}}$. On $Q_N\cap W$, assign customers to the centers of their full cells. On the remaining region, use the distance bound above. Then
\[
\begin{aligned}
M_{\mathrm{row}}&\le mn\phi(d)+\frac74\sqrt{d^2+d^{-2}}\,|Q_N\setminus W|\\
&\le N\phi(d)+\frac{21\sqrt2}{8}(d^2+d^{-2})s\\
&\le N\phi(d)+4(d^2+d^{-2})s.
\end{aligned}
\]
Traverse alternate rows in opposite directions. Each within-row edge has length $d$, and each row connection has length $\sqrt{d^{-2}+d^2/4}$. A closing edge has length at most $\sqrt2s$. Since $n\le sd$, the resulting route has length
\[
\begin{aligned}
n(m-1)d+(n-1)\sqrt{d^{-2}+d^2/4}+\sqrt2s
\le mnd+\{\sqrt{1+d^4/4}+\sqrt2\}s.
\end{aligned}
\]
Place the remaining $q$ distinct sites in a ball of radius $\varepsilon\le1/(2N)$ around an interior row site $A$. Choose the ball to lie inside the square and to contain no other row site, and choose every new site different from $A$. Inserting a site $z$ after $A$ in the route replaces $AB$ by $Az,zB$ and adds at most $|Az|+|zB|-|AB|\le2\varepsilon$. If $mn=1$, start with the zero-length loop $AA$. All insertions add at most $2q\varepsilon\le1\le s$, and none increases service cost, so $M_N(X)\le M_{\mathrm{row}}$. The final tour satisfies
\[
\begin{aligned}
L(X,\sigma)&\le Nd+\{\sqrt{1+d^4/4}+\sqrt2+1\}s\\
&\le Nd+\{2/\sqrt3+\sqrt2+1\}s\le Nd+4s.
\end{aligned}
\]
The resulting layout $X$ and closed visiting order $\sigma$ satisfy \eqref{eq:upperparts}, and $L_N(X)\le L(X,\sigma)$.

For the given $r$, choose a minimizer $d$ of $\phi(d)+rd$ using Lemma~\ref{lem:existence}. Adding the weighted routing bound gives
\[
A_N(r)\le N\{\phi(d)+rd\}+4(d^2+d^{-2}+r)\sqrt N,
\]
which is \eqref{eq:upper}. The spacing, and hence $D_r$, is independent of $N$. An alternative constant involving only $r$ follows by putting $U=\phi(d_0)+rd_0$. The inequalities
\[
\frac1{4d}\le\phi(d)\le\phi(d)+rd=C_6(r)\le U
\]
give $d^{-2}\le16U^2$. Therefore $D_r\le4d_0^2+64U^2+4r$, and this explicit upper estimate can replace $D_r$ in the cost bound for $N\ge\max\{3,\lceil16U^2\rceil\}$. Dividing by $N$ also gives $\limsup_{N\to\infty}A_N(r)/N\le C_6(r)$.
\end{proof}

\section{Local lower-bound proofs}\label{app:local}
\subsection{At most \texorpdfstring{$6N$}{6N} side triangles}\label{app:side}
\sideCount*
\begin{proof}
Use the planar graph formed by the clipped Voronoi boundaries and the square boundary.
Remove degree-two vertices lying on straight segments, keeping the four square corners.
If a cell side contains a vertex of this graph, split the side at that vertex. The resulting number of pieces
is at least the number of original cell sides. Every cell is convex with nonempty interior,
so the graph is a connected subdivision of the square into $N$ faces.
Let $V_i,V_b$ count interior and noncorner boundary vertices, and let $E_i$ count interior edges. There are $V_b+4$ boundary edges, one between each consecutive pair of boundary vertices. Euler's formula gives $E_i=V_i+N-1$. Interior and noncorner boundary vertices have degree at least three, and corners have degree at least two. Thus
\[
2(E_i+V_b+4)\ge3V_i+3V_b+8.
\]
Substitute $E_i=V_i+N-1$ into this inequality to obtain
\[
V_i+V_b\le2N-2.
\]
Counting each interior edge twice and each boundary edge once bounds the total side count $k$:
\[
k\le2E_i+V_b+4=2V_i+2N+V_b+2\le6N-2-V_b\le6N.
\]
The degree bounds allow more than three cells to meet at an interior vertex and allow an
interior edge to end at a corner. Neither degeneracy changes the estimate.
\end{proof}

\subsection{Cap and overlap derivatives}\label{app:overlap}
For $0\le h<t$, a vertical slice gives $c_h(t)=\int_h^t2\sqrt{t^2-u^2}\dd u$. Differentiating this identity and using \eqref{eq:overlap} gives the cap derivative \eqref{eq:capderiv} and the overlap derivatives:
\begin{equation}\label{eq:app-overlapderiv}
\partial_hc_h(t)=-2\sqrt{t^2-h^2},\qquad
o_t'(\ell)=-\sqrt{4t^2-\ell^2},\qquad
o_t''(\ell)=\frac{\ell}{\sqrt{4t^2-\ell^2}}>0.
\end{equation}
The cap derivative applies for $0<h<t$, and the last two formulas apply for $0<\ell<2t$. Since the slice integrand tends to zero at $u=t$, the cap derivative extends continuously through $h=t>0$ with value zero, agreeing with the branch $c_h(t)=0$ for $h\ge t$. When $\ell\ge2t$, the disks have no overlap of positive area, so $o_t(\ell)=0$. As $\ell$ increases to $2t$ from below, both $o_t(\ell)$ and $o_t'(\ell)$ tend to zero, matching this constant branch. The derivative is increasing on $(0,2t)$ and remains zero on $(2t,\infty)$, so $o_t$ is convex on $[0,\infty)$. When $t=0$, both disks have area zero and $o_0(\ell)=0$ for every $\ell\ge0$, which is also convex.

\subsection{The minimum of \texorpdfstring{$V$}{V}}\label{app:V}
\minimumOfV*
\begin{proof}
Differentiating \eqref{eq:V} gives
\begin{align}
V'(h)&=(R-2S)\sqrt{R^2-h^2}+h^2\arcosh(R/h),\label{eq:app-Vprime}\\
V''(h)&=2h\left(\arcosh(R/h)-\frac{R-S}{\sqrt{R^2-h^2}}\right).\label{eq:app-Vsecond}
\end{align}
To determine the sign of $V'(h)$, put $p(h)=h^2\arcosh(R/h)/\sqrt{R^2-h^2}$ and $u=\sqrt{1-h^2/R^2}$. Thus $0<u<1$, $h^2=R^2(1-u^2)$ and $\sqrt{R^2-h^2}=Ru$. Also,
\[
\arcosh(R/h)=\arcosh\!\left(\frac1{\sqrt{1-u^2}}\right)
=\frac12\log\frac{1+u}{1-u}=\artanh u.
\]
Substituting these expressions and differentiating gives
\[
\frac{p(h)}R=(u^{-1}-u)\artanh u,\qquad
\frac{\mathrm d}{\mathrm du}\{(u^{-1}-u)\artanh u\}
=\frac{u-(1+u^2)\artanh u}{u^2}<0,
\]
using $\artanh u=\int_0^u(1-v^2)^{-1}\dd v>u$ for $0<u<1$. As $h$ increases, $u$ decreases, so $p$ increases strictly from 0 to $R$. The choice $2S=R+p(b)$ gives $R/2<S<R$ and
\[
V'(h)=\sqrt{R^2-h^2}\{p(h)-p(b)\}.
\]
This is negative before $b$ and positive after $b$. Every term of $V$ tends to zero at $R$, so $V(R)=0$ and the minimum $V(b)$ is negative.
\end{proof}

\subsection{Integrability and slope cancellation}\label{app:eta}
\overhangWeight*
\begin{proof}
\textbf{1. Positivity and integrability.} From the definition of $q$, we have $V''(h)=2h\,q(h^2)$ and
\[
q'(y)=-\frac{R}{2y\sqrt{R^2-y}}-\frac{R-S}{2(R^2-y)^{3/2}}<0,\qquad0<y\le b^2,
\]
where $S<R$ by Lemma~\ref{lem:V}. To check the sign of $q(b^2)$, put $u=\sqrt{1-b^2/R^2}\in(0,1)$. Substituting \eqref{eq:S} gives
\[
q(b^2)=\frac{(1+u^2)\artanh u-u}{2u^2}>0,
\]
using $\artanh u=\int_0^u\frac{\dd v}{1-v^2}>u$. Thus both terms in brackets in \eqref{eq:eta} are positive, and $\eta(t)\ge0$ for $0<t<b$.

Fix $0<a<b$. The derivative $q'$ is continuous on $[a^2,b^2]$, so
\[
M:=\max_{y\in[a^2,b^2]}|q'(y)|<\infty.
\]
For $a\le t<b$,
\[
0\le-\int_{t^2}^{b^2}\frac{q'(y)}{\sqrt{y-t^2}}\dd y
\le M\int_{t^2}^{b^2}\frac{\dd y}{\sqrt{y-t^2}}=2M\sqrt{b^2-t^2}.
\]
Substituting this bound into \eqref{eq:eta} yields
\[
0\le\eta(t)\le\frac{2t}{\pi}\left[\frac{q(b^2)}{\sqrt{b^2-t^2}}+2M\sqrt{b^2-t^2}\right].
\]
Integrating the two terms gives
\[
\int_a^b\eta(t)\dd t\le\frac{2q(b^2)}\pi\sqrt{b^2-a^2}
+\frac{4M}{3\pi}(b^2-a^2)^{3/2}<\infty.
\]
For $a=b$, the integral is zero. The value assigned to $\eta(b)$ does not affect any of these integrals.

\textbf{2. Evaluate the square-root transform.} For $x<y$, the substitution $v=x+(y-x)\sin^2\theta$ gives
\begin{equation}\label{eq:app-beta}
\int_x^y\frac{\dd v}{\sqrt{(y-v)(v-x)}}=2\int_0^{\pi/2}\dd\theta=\pi.
\end{equation}
Fix $0<h<b$. Insert \eqref{eq:eta} and set $v=t^2$:
\[
\begin{aligned}
\int_h^b\frac{\eta(t)}{\sqrt{t^2-h^2}}\dd t
&=\frac{q(b^2)}{\pi}\int_{h^2}^{b^2}\frac{\dd v}{\sqrt{(b^2-v)(v-h^2)}}\\
&\quad+\frac1\pi\int_{h^2}^{b^2}\int_v^{b^2}\frac{-q'(y)}{\sqrt{(y-v)(v-h^2)}}\dd y\dd v.
\end{aligned}
\]
The first integral equals $\pi$ by \eqref{eq:app-beta}. For the double integral, the region is $h^2<v<y<b^2$. Its integrand is nonnegative, so Tonelli's theorem \citep[Theorem~1.7.18]{Tao2011} permits reversing the order. Applying \eqref{eq:app-beta} again gives
\[
\begin{aligned}
\int_h^b\frac{\eta(t)}{\sqrt{t^2-h^2}}\dd t
&=q(b^2)+\frac1\pi\int_{h^2}^{b^2}(-q'(y))
\left[\int_{h^2}^y\frac{\dd v}{\sqrt{(y-v)(v-h^2)}}\right]\dd y\\
&=q(b^2)-\int_{h^2}^{b^2}q'(y)\dd y=q(h^2).
\end{aligned}
\]
In particular, the transform is finite for every $0<h<b$. Its limit as $h\uparrow b$ is
$q(b^2)>0$. At $h=b$, the integration interval is empty, so the integral itself is zero.
The slope identity proved next, unlike this inverse-square-root transform, does extend to $h=b$.

\textbf{3. Integrate to cancel the slope.} For $h<t$, direct integration gives
\[
2\sqrt{t^2-h^2}=\int_h^t\frac{2u}{\sqrt{t^2-u^2}}\dd u.
\]
Using this identity and Tonelli's theorem for the nonnegative integrand, we obtain
\[
\begin{aligned}
2\int_h^b\eta(t)\sqrt{t^2-h^2}\dd t
&=\int_h^b\eta(t)\left[\int_h^t\frac{2u}{\sqrt{t^2-u^2}}\dd u\right]\dd t\\
&=\int_h^b2u\left[\int_u^b\frac{\eta(t)}{\sqrt{t^2-u^2}}\dd t\right]\dd u\\
&=\int_h^b2u\,q(u^2)\dd u
=\int_h^b V''(u)\dd u\\
&=V'(b)-V'(h)=-V'(h).
\end{aligned}
\]
The third line uses the identity from Step 2, and the last line uses $V'(b)=0$. At $h=b$, both sides of the slope identity are zero.
\end{proof}

\subsection{Side-triangle inequality}\label{app:triangle}
\sideTriangleInequality*
\begin{proof}
Fix the reference spacing $d$ and its parameters $a$, $b$ and $R$ from \eqref{eq:abR}. The triangle being bounded may have any admissible side distance $h$ and endpoints $\theta_-,\theta_+$. Its outer vertices need not lie on the circle of radius $R$ centered at the facility.

\textbf{1. Bounds for $F$.} For $0<h<R$, put
\begin{equation}\label{eq:app-Fintegral}
F(h)=\int_{-\arccos(h/R)}^{\arccos(h/R)}\Psi(h\sec\theta)\dd\theta.
\end{equation}
Write $\beta=\arccos(h/R)$. The first three terms of $\Psi(s)$ are $s^3/3-Ss^2/2+\lambda$. Their angular integral is
\[
\int_{-\beta}^{\beta}
\left(\frac{h^3\sec^3\theta}{3}-\frac{Sh^2\sec^2\theta}{2}+\lambda\right)\dd\theta
=J_R(h)-Sh\sqrt{R^2-h^2}+2\lambda\beta=V(h).
\]
For $h<t<R$, set $\alpha=\arccos(h/t)$. The remaining integral in $\Psi$ contributes a cap, since
\[
\frac12\int_{-\alpha}^{\alpha}(t^2-h^2\sec^2\theta)\dd\theta
=t^2\alpha-h^2\tan\alpha=c_h(t).
\]
For $t\le h$, the positive-part integrand in \eqref{eq:Psi} is zero. Thus
\begin{equation}\label{eq:app-Fcap}
F(h)=V(h)-\int_a^b\eta(t)c_h(t)\dd t.
\end{equation}
The bound $|\partial_hc_h(t)|\le2b$ and integrability of $\eta$ permit differentiation under the integral. Using \eqref{eq:capderiv}, for $h<b$ we obtain
\[
F'(h)=V'(h)+2\int_{\max\{a,h\}}^b\eta(t)\sqrt{t^2-h^2}\dd t.
\]
The slope identity in Lemma~\ref{lem:eta} states that
\[
V'(h)=-2\int_h^b\eta(t)\sqrt{t^2-h^2}\dd t,\qquad 0<h<b.
\]
If $a\le h<b$, substitution into the preceding formula gives $F'(h)=0$. If $0<h<a$, split the integral from $h$ to $b$ at $a$. Its part over $[a,b]$ cancels, leaving
\[
F'(h)=-2\int_h^a\eta(t)\sqrt{t^2-h^2}\dd t.
\]
For $h\ge b$, every cap in the correction is zero. At $b$, the cap derivative also vanishes, so $F'(b)=V'(b)=0$. For $b<h<R$, \appV{} gives $F'(h)=V'(h)=\sqrt{R^2-h^2}\,[p(h)-p(b)]>0$, since the function $p$ defined there is strictly increasing. Together with $F(b)=V(b)$, these identities give
\begin{equation}\label{eq:app-Fproperties}
\begin{cases}
F'(h)=-2\displaystyle\int_h^a\eta(t)\sqrt{t^2-h^2}\dd t<0,&0<h<a,\\[4pt]
F(h)=V(b),\quad F'(h)=0,&a\le h\le b,\\[2pt]
F(h)=V(h)<0,\quad F'(h)>0,&b<h<R.
\end{cases}
\end{equation}
Hence $F(h)\ge V(b)$ for $0<h<R$ and $F(R)=0$ by continuity (Figure~\ref{fig:VF}). The same reasoning applies when $a=b$.

\textbf{2. Identify the negative radial interval.} We use the following Abel inversion formula \citep{GorenfloVessella1991}:
\begin{equation}\label{eq:app-inversion}
\Psi(s)=-\frac{s^2}{\pi}\int_s^R
\frac{F'(h)/h-F(h)/h^2}{\sqrt{h^2-s^2}}\dd h,\qquad0<s<R.
\end{equation}
We verify the formula for this $\Psi$, including the changes of integration order, differentiation, and endpoint conditions. Set $T=R^2$, $f(y)=\Psi(\sqrt y)/y$, and $g(x)=F(\sqrt x)/\sqrt x$. Changing variables $s=h\sec\theta$ in \eqref{eq:app-Fintegral} yields $g(x)=\int_x^T f(y)/\sqrt{y-x}\dd y$. Apply this transform twice and use \eqref{eq:app-beta} to obtain
\[
\int_z^T\frac{g(x)}{\sqrt{x-z}}\dd x=\pi\int_z^T f(y)\dd y,\qquad0<z<T.
\]
The interchange is valid by absolute integrability: $f$ is bounded on $[z,T]$ and the inner beta integral is $\pi$.
For $s>0$, differentiation of \eqref{eq:Psi} gives
\[
\Psi'(s)=s(s-S)+s\int_{\max\{a,s\}}^b\eta(t)\dd t,
\]
where the last integral is read as zero when $s\ge b$.
The tail integral is continuous since $\eta\in L^1([a,b])$, including at $a$ and $b$.
Thus $\Psi$ has a continuous first derivative on every interval bounded away from zero,
and $\Psi(R)=0$. The transform can also be written as $g(x)=2\int_0^{\sqrt{T-x}}f(x+u^2)\dd u$. Since $f(T)=0$ and $f'$ is bounded, differentiating this representation gives no endpoint term:
\[
g'(x)=2\int_0^{\sqrt{T-x}}f'(x+u^2)\dd u,\qquad
|g'(x)|\le2\sqrt{T-x}\sup_{y\in[z,T]}|f'(y)|\quad(z\le x<T).
\]
Thus $g(T)=0$ and $g'$ extends continuously to $g'(T)=0$ and is bounded on $[z,T]$. Integration by parts, followed by differentiation with respect to $z$, gives
\[
\pi\int_z^T f(y)\dd y=-2\int_z^Tg'(x)\sqrt{x-z}\dd x,\qquad
-\pi f(z)=\int_z^T\frac{g'(x)}{\sqrt{x-z}}\dd x.
\]
The boundary terms vanish and the remaining square-root singularity is integrable. Finally, substitute $x=h^2$, $z=s^2$, and $2h g'(h^2)=F'(h)/h-F(h)/h^2$ to get \eqref{eq:app-inversion}.

For $a\le s<R$, \eqref{eq:app-Fproperties} gives $F<0$ and $F'\ge0$ on $(s,R)$, so the integral in \eqref{eq:app-inversion} is positive and $\Psi(s)<0$. For $0<s<a$, direct differentiation gives
\[
\Psi'(s)=s\left(s-S+\int_a^b\eta(t)\dd t\right).
\]
The bracket increases strictly with $s$. The only possible positive turning point is therefore a minimum. Since $\Psi(a)<0$, any increasing portion below $a$ stays negative. Any zero there is unique and crossed from positive to negative. For $s>R>b$, we have $(t^2-s^2)_+=0$ for every $t\in[a,b]$. The integral term in \eqref{eq:Psi} is therefore zero, and $\Psi'(s)=s(s-S)>0$, with $\Psi(R)=0$. Thus, for some $\rho\in[0,a)$,
\begin{equation}\label{eq:app-sign}
\Psi\ge0\text{ on }(0,\rho),\qquad
\Psi<0\text{ on }(\rho,R),\qquad
\Psi>0\text{ on }(R,\infty).
\end{equation}
If $\rho>0$, continuity gives $\Psi(\rho)=0$. For $\rho=0$, the first interval is empty.

\textbf{3. Sum the negative angular contributions.} Fix $I=[\theta_-,\theta_+]$. An integral over $I$ is bounded below by the sum of all negative contributions over $(-\pi/2,\pi/2)$: removing negative pieces or adding nonnegative pieces can only increase it. These negative contributions are integrable and occur away from $\pm\pi/2$. If $h\ge R$, the integrand is nonnegative in every direction, so the integral is at least $0>V(b)$.

If $\rho\le h<R$, the directions with negative integrand form the symmetric interval $(-\beta,\beta)$, with $\beta=\arccos(h/R)$, apart from a possible zero at $\theta=0$. Its integral is $F(h)$, so
\[
\int_I\Psi(h\sec\theta)\dd\theta\ge F(h)\ge V(b).
\]
If $0<h<\rho$, the directions with negative integrand form two intervals. Substitution $s=h\sec\theta$ on the positive half-axis and evenness give
\[
\int_I\Psi(h\sec\theta)\dd\theta\ge2\int_\rho^R\Psi(s)\frac{h}{s\sqrt{s^2-h^2}}\dd s.
\]
For fixed $s>\rho$, the factor $h/(s\sqrt{s^2-h^2})$ multiplying $\Psi(s)$ increases with $h$ on $0<h<s$, since
\[
\frac{\partial}{\partial h}\frac{h}{s\sqrt{s^2-h^2}}
=\frac{s}{(s^2-h^2)^{3/2}}>0.
\]
With $h<\rho$ and $\Psi(s)<0$, multiplication reverses the comparison. Hence
\[
2\int_\rho^R\Psi(s)\frac{h}{s\sqrt{s^2-h^2}}\dd s
\ge2\int_\rho^R\Psi(s)\frac{\rho}{s\sqrt{s^2-\rho^2}}\dd s
=F(\rho)\ge V(b).
\]
The comparison integral is finite by continuity of $\Psi$ and integrability of $(s-\rho)^{-1/2}$ at $\rho>0$. Figure~\ref{fig:signs} shows both cases. This proves \eqref{eq:triineq} for every admissible angular interval. For the symmetric triangles in the lemma, the integral equals $F(h)=V(b)$ by \eqref{eq:app-Fproperties}.
\end{proof}

\section{Global lower-bound proofs}\label{app:global}
\subsection{Service cost and uncovered area}\label{app:coverage}

The side-triangle bound shows how a weighted excess in uncovered area forces an increase in service cost. Fix $d$ and its reference parameters and weight. Every layout of distinct interior facilities satisfies
\begin{equation}\label{eq:coverage}
M_N(X)-N\phi(d)\ge\int_a^b\eta(t)[U_X(t)-Nu_d(t)]\dd t,
\qquad u_d(t)=1-\pi t^2+o_t(d).
\end{equation}

\begin{proof}
\textbf{1. Sum the overhangs.} For fixed $t$, use $(t^2-s^2)_+=t^2-\min\{t^2,s^2\}$ in each side triangle. The angles at each facility sum to $2\pi$, so the $t^2$ terms contribute $\pi t^2$ per cell. The subtracted part is the area of that cell within its radius-$t$ disk. Since the cells partition $Q_N$, these in-cell areas sum to $N-U_X(t)$. Therefore
\begin{equation}\label{eq:app-areas}
\sum_{\text{triangles}}\frac12\int(t^2-(h\sec\theta)^2)_+\dd\theta
=N\pi t^2-N+U_X(t).
\end{equation}
Figure~\ref{fig:triangle} illustrates the decomposition. Nonnegativity permits the weighted sum and the angular integrals to be exchanged by Tonelli's theorem.

\textbf{2. Add the local service bounds.} For one side triangle $T$, expanding \eqref{eq:Psi} in \eqref{eq:triineq} and isolating its service cost gives
\[
\begin{aligned}
\int_T\|x\|\dd x
&\ge S|T|+\left(\frac{R^3}3-\frac{SR^2}2\right)(\theta_+-\theta_-)+V(b)\\
&\quad+\frac12\int_a^b\eta(t)\int_{\theta_-}^{\theta_+}(t^2-(h\sec\theta)^2)_+\dd\theta\dd t.
\end{aligned}
\]
The coordinates here are relative to the triangle's facility. Sum over the $k$ triangles. Service costs sum to $M_N(X)$, areas to $N$, and angular widths to $2\pi N$. Lemma~\ref{lem:sides} and $V(b)<0$ give $kV(b)\ge6NV(b)$. Using \eqref{eq:app-areas}, we obtain
\begin{equation}\label{eq:app-sum}
\begin{aligned}
M_N(X)&\ge N\left[S+\pi\left(\frac{2R^3}3-SR^2\right)+6V(b)\right]\\
&\quad+\int_a^b\eta(t)\{N\pi t^2-N+U_X(t)\}\dd t.
\end{aligned}
\end{equation}
\textbf{3. Subtract the reference equality.} The reference hexagon has area one, total angle $2\pi$, and six triangles attaining equality. For $a\le t\le b$, only the two sides at distance $a=d/2$ contribute caps, whose total area is $2c_a(t)=o_t(d)$. The other four caps have area zero. Summing the six equality cases gives
\begin{equation}\label{eq:app-reference}
\phi(d)=S+\pi\left(\frac{2R^3}3-SR^2\right)+6V(b)+\int_a^b\eta(t)o_t(d)\dd t.
\end{equation}
Subtract $N$ times this equality from \eqref{eq:app-sum}. The remaining integral is exactly the uncovered-area difference in \eqref{eq:coverage}.
\end{proof}

\subsection{Extension to boundary facilities}\label{app:boundary}
\boundaryFacilities*
\begin{proof}
Let $z$ be the square's center and set $x_i^\varepsilon=(1-\varepsilon)x_i+\varepsilon z$, $0<\varepsilon<1$, with $X^\varepsilon=\{x_1^\varepsilon,\ldots,x_N^\varepsilon\}$. These sites are interior and distinct: their pairwise differences are multiplied by $1-\varepsilon$. Each site moves at most $\varepsilon\sqrt{N/2}$. The triangle inequality, applied to each site and then minimized over sites, gives
\[
|\dist(x,X^\varepsilon)-\dist(x,X)|\le\varepsilon\sqrt{N/2}.
\]
Integrate over the same square of area $N$. For the unchanged visiting order, every tour edge is multiplied by $1-\varepsilon$, so
\[
|M_N(X^\varepsilon)-M_N(X)|\le N\varepsilon\sqrt{N/2}\to0,
\qquad L^\varepsilon=(1-\varepsilon)L\to L.
\]
With $N$ and $d$ fixed, apply Lemma~\ref{lem:tradeoff} to $X^\varepsilon$ and let $\varepsilon\downarrow0$. The reference parameters and constants are unchanged.
\end{proof}

\subsection{Optimal spacing and routing weight}\label{app:spacing}
\optimalSpacing*
\begin{proof}
For a symmetric triangle, write
\[
A(h,R)=h\sqrt{R^2-h^2},\qquad
\Theta(h,R)=2\arccos(h/R),\qquad
\lambda=\frac{SR^2}{2}-\frac{R^3}{3}.
\]
Fix a reference spacing $d$. In the partial derivatives below, $S$ and $\lambda$ are constants. Set $P(h,R)=J_R(h)-SA(h,R)+\lambda\Theta(h,R)$. At the reference radius, $\partial_hP(h,R)=V'(h)$. Direct differentiation gives
\[
\partial_RJ_R(h)=\frac{2hR^2}{3\sqrt{R^2-h^2}},\quad
\partial_RA(h,R)=\frac{hR}{\sqrt{R^2-h^2}},\quad
\partial_R\Theta(h,R)=\frac{2h}{R\sqrt{R^2-h^2}}.
\]
Substituting the chosen value of $\lambda$ gives, at the reference radius,
\[
\partial_RP(h,R)=\frac{h}{\sqrt{R^2-h^2}}
\left(\frac{2R^2}{3}-SR+\frac{2\lambda}{R}\right)=0.
\]
The six triangles satisfy
\[
2A(a,R)+4A(b,R)=1,\qquad
2\Theta(a,R)+4\Theta(b,R)=2\pi.
\]
Differentiate $\phi(d)=2J_R(a)+4J_R(b)$, subtract $S$ times the derivative of the first identity, and add $\lambda$ times the derivative of the second. Both identities have zero derivative, so
\[
\begin{aligned}
\phi'(d)&=2\partial_hP(a,R)a'(d)+4\partial_hP(b,R)b'(d)\\
&\quad+\{2\partial_RP(a,R)+4\partial_RP(b,R)\}R'(d)\\
&=2V'(a)\frac12+4V'(b)b'(d)=V'(a).
\end{aligned}
\]
We used $a'(d)=1/2$ and $V'(b)=0$. Since $d=2a$, the definition \eqref{eq:rK} and the slope identity \eqref{eq:slope} at $h=a$ give
\[
r_*(d)=2\int_a^b\eta(t)\sqrt{t^2-a^2}\dd t=-V'(a)=-\phi'(d).
\]
The identity at $d_0$ is understood from the left.

Lemma~\ref{lem:existence} gives a positive minimizer of $\phi(d)+rd$. For $r>0$, the left derivative at $d_0$ is $r>0$, so a minimizer lies in $(0,d_0)$ and satisfies $\phi'(d)+r=0$. For $r=0$, Lemma~\ref{lem:V} gives $\phi'(d)=V'(a)<0$ for $d<d_0$, and $d_0$ is a minimizer. In both cases a minimizing spacing satisfies $r_*(d)=r$.

The map $r_*(d)=-V'(a)$ is continuous on $(0,d_0]$. At $d_0$ we have $a=b$ and $r_*(d_0)=0$. As $d\downarrow0$,
\[
R=\frac1{2d}\left(1+\frac{d^4}4\right),\qquad
\frac bR=\frac1{\sqrt{1+d^4/4}}\to1,\qquad
\frac aR=\frac{d^2}{1+d^4/4}\to0.
\]
Integrating $u/R\le1\le u/b$ against $(u^2-b^2)^{-1/2}\dd u$ on $(b,R)$ gives
\[
\frac{\sqrt{R^2-b^2}}{R}\le\arcosh(R/b)\le\frac{\sqrt{R^2-b^2}}{b}.
\]
Substitution in \eqref{eq:S} yields
\[
\frac12\left(1+\frac{b^2}{R^2}\right)\le\frac SR
\le\frac12\left(1+\frac bR\right)\longrightarrow1.
\]
By \eqref{eq:app-Vprime},
\[
\frac{r_*(d)}{R^2}
=\left(\frac{2S}R-1\right)\sqrt{1-(a/R)^2}-(a/R)^2\arcosh(R/a)\longrightarrow1,
\]
where $z^2\arcosh(1/z)\le z^2\log(2/z)\to0$. Hence $r_*(d)\sim1/(4d^2)$.

\Needspace{6\baselineskip}
\paragraph{Strict convexity.}
Put $x=d^2/2$, so $0<x\le1/\sqrt3$. Substituting the geometric formulas into $\rs(d)=-V'(a)$ gives
\[
\rs(d)=G(x):=\frac{(1-x^2)\sqrt{1+x^2}}{8x^2}\operatorname{arsinh}x
+\frac{x}{2}\log x.
\]
Here $\arcosh(R/a)=-\log x$ and $\arcosh(R/b)=\operatorname{arsinh}x$.
Differentiation yields
\[
G'(x)=\frac{\log x}{2}+\frac38+\frac1{8x^2}
-\frac{x^4+x^2+2}{8x^3\sqrt{1+x^2}}\operatorname{arsinh}x.
\]
Since $\operatorname{arsinh}x=\int_0^x(1+t^2)^{-1/2}\dd t\ge x/\sqrt{1+x^2}$,
\[
G'(x)\le\frac{\log x}{2}
+\frac{2x^4+3x^2-1}{8x^2(1+x^2)}
\le\frac1{16}-\frac{\log3}{4}<0.
\]
For the second inequality, $x^2\le1/3$ implies
$(3x^2-1)(x^2+2)\le0$, which bounds the rational term by $1/16$. Also, $\log x\le-\tfrac12\log3$.
Thus $\phi''(d)=-dG'(d^2/2)>0$ on $(0,d_0)$, and continuity extends strict convexity to $(0,d_0]$. The minimizer of $\phi(d)+rd$ is thus unique.
\end{proof}


\end{document}